\documentclass[reqno, 11pt, letterpaper]{amsart}

  \title[A 0-transition calculus with applications to gas giants]{A transition calculus for 0-pseudodifferential operators with applications to gas giant geometry}
  \author{Andy Yin}
  \date{November 2, 2026}

    \usepackage{amssymb}
  \usepackage{latexsym}
  \usepackage{amsfonts}
  \usepackage{amsthm}
  \usepackage{amsmath}
  \allowdisplaybreaks
  \usepackage{xcolor}
  \usepackage{bbm}
  \usepackage{quiver}
  \usepackage{mathrsfs}
  \usepackage{xparse}
  \usepackage{titletoc}
  \usepackage{adjustbox}
  \usepackage{microtype}
  \usepackage[hidelinks]{hyperref}
  \usepackage[backend=biber,style=alphabetic,sorting=nyt,citestyle=alphabetic]{biblatex}
  \DeclareFieldFormat{postnote}{#1}

    \NewDocumentCommand{\DeclareOptionalSubscriptCommand}{m m}{%
  \NewDocumentCommand{#1}{o o o}{%
    {\IfValueTF{##3}{##3{\mathrm{#2}}}{\mathrm{#2}}}%
    \IfValueT{##1}{_{##1}}%
    \IfValueT{##2}{^{##2}}%
  }%
}

    \def\R{{\mathbb R}}
    
    \def\Z{{\mathbb Z}}
    \def\C{{\mathbb C}}
    \def\N{{\mathbb N}}
    
    \newcommand{\gas}{\mathrm{0,tr}}
    \newcommand{\btr}{\mathrm{b,tr}}
    \newcommand{\sctr}{\mathrm{sc,tr}}
    \newcommand{\dbsctr}{\mathrm{2b,2sc,tr}}
    \newcommand{\ztrtext}{0-transition}
    \newcommand{\btrtext}{b-transition}

    \NewDocumentCommand{\sgsp}{m m}{#1^{1}_{#2}}
    \NewDocumentCommand{\dbsp}{m m}{#1^{2}_{#2}}
    \NewDocumentCommand{\trsp}{m m}{#1^{3}_{#2}}
    \NewDocumentCommand{\anysp}{m m o}{%
  #1_{#2\IfValueT{#3}{,#3}}%
}

    \NewDocumentCommand{\anydiag}{m o o}{%
  {\IfValueTF{#3}{#3{\mathrm{d}}}{\mathrm{d}}}_{#1}%
  \IfValueT{#2}{^{#2}}%
}
    
    \newcommand{\diaglift}{\mathrm{d}_{12}}

    \RenewDocumentCommand{\bf}{m o o}{%
  {\IfValueTF{#3}{#3{\mathrm{b}}}{\mathrm{b}}}_{#1}%
  \IfValueT{#2}{^{#2}}%
}
    \DeclareOptionalSubscriptCommand{\ff}{ff}
    \DeclareOptionalSubscriptCommand{\bff}{bff}
    \DeclareOptionalSubscriptCommand{\tff}{tff}
    \DeclareOptionalSubscriptCommand{\tf}{tf}
    
    \DeclareOptionalSubscriptCommand{\zff}{zff}
    \DeclareOptionalSubscriptCommand{\uff}{zff} %Deprecated
    \DeclareOptionalSubscriptCommand{\ztff}{ztff} %Used only if necessary to disambiguate zero-transition front face from b-transition front face
    \DeclareOptionalSubscriptCommand{\btff}{btff}
    
    \DeclareOptionalSubscriptCommand{\scff}{scff}
    \DeclareOptionalSubscriptCommand{\sctff}{sctff}

  \newtheorem{theorem}{Theorem}[section]
  \newtheorem{lemma}[theorem]{Lemma}
  \newtheorem{proposition}[theorem]{Proposition}
  \newtheorem{corollary}[theorem]{Corollary}

  \theoremstyle{definition}
  \newtheorem{definition}[theorem]{Definition}

  \usepackage{graphicx}

\begin{document}

  \begin{abstract}
      A gas giant Riemannian metric $g$ on a compact $n$-dimensional manifold-with-boundary $M$ blows up at the boundary $\partial M$ at the rate $1/x^\alpha$, for $\alpha \in (0,2)$, where $x \in C^\infty(M)$ is a boundary-defining function. Such metrics arise in the modeling of acoustic waves in gaseous planets. In \cite{dHIKM_gas-giants_2024}, it was shown that, if $\alpha > 2/n$, the Laplace-Beltrami operator $\Delta_g$ is essentially self-adjoint and has discrete spectrum. Moreover, its spectrum is stable under truncations of $M$. The proof leverages the 0-calculus of Mazzeo-Melrose \cite{mazzeo-melrose_meromorphic-extension-resolvent_1987}. We extend this theory by introducing a \emph{0-transition calculus} for studying families of 0-operators. We demonstrate a parametrix construction for $\Delta_g$ that is stable under truncations. Using this, we study the behavior of the resolvent of $\Delta_g$ under truncation and extend the characterization of the spectrum obtained in \cite{dHIKM_gas-giants_2024}. In particular, we prove that (on appropriate functional-analytic domains) simple eigenvalues admit asymptotic expansions in the truncation parameter $\epsilon$ and that the corresponding eigenfunctions admit joint asymptotic expansions in $(\epsilon, x)$.
  \end{abstract}
  \maketitle

  \section{Introduction}

Let $M$ be a compact $n$-dimensional manifold with boundary, where $n \geq 3$, and let $x \in C^\infty(M)$ be a boundary-defining function, meaning that $x \geq 0$, $x|_{\partial M} = 0$, and $dx|_{\partial M}$ is nonvanishing. A \emph{gas giant metric} on $M$ is of the form $g = \overline{g}/x^\alpha$, where $\overline{g}$ is a smooth Riemannian metric on $M$ up to $\partial M$ and $\alpha \in (0,2)$ \cite{dHIKM_gas-giants_2024}. The excluded extremes $\alpha = 0, 2$ correspond, respectively, to a smooth Riemannian metric and an asymptotically hyperbolic metric. Unlike an asymptotically hyperbolic metric, a gas giant metric is incomplete and its boundary lies at finite distance.   

The term `gas giant metric' was introduced by \cite{dHIKM_gas-giants_2024}. Such metrics are derived from seismological models of gaseous planets \cite{prat_asymptotic_2016} \cite{dewberry_constraining_2021}. Recently, they have spurred work in spectral geometry and geometric inverse problems \cite{dHIKM_gas-giants_2024} \cite{de_hoop_characterization_2025} \cite{ilmavirta_tomography_2026}.

In \cite{dHIKM_gas-giants_2024}, it was shown that, if $\alpha > 2/n$, the Laplace-Beltrami operator $\Delta_g$ is essentially self-adjoint and has discrete spectrum. It was also shown that the spectrum is stable under truncations. Truncation is a ubiquitous technique in the modeling of gaseous planets to produce numerically tractable models \cite{dewberry_constraining_2021}. In our context, truncation means the following: let $M_\epsilon = \{ p \in M : x(p) \geq \epsilon \}$ (after rescaling $x$ if necessary, so that this set is nonempty for $\epsilon \in [0,1)$) and let $\Delta_\epsilon = \Delta_g|_{M_\epsilon}$. In \cite{dHIKM_gas-giants_2024}, it was shown that, if $\alpha >  2/n$, the Dirichlet eigenvectors and eigenvalues of $\Delta_\epsilon$ converge to the eigenvectors and eigenvalues of $\Delta_g$ on its domain of self-adjointness.

To prove these results, \cite{dHIKM_gas-giants_2024} leveraged the fact that the rescaled operator $x^{2-\alpha} \Delta_g$ is an elliptic \emph{0-differential operator}, for which there exists a robust theory \cite{mazzeo-melrose_meromorphic-extension-resolvent_1987} \cite{mazzeo_differential-edge-operators-I_1991} \cite{hintz_0-calculus-mazzeo-melrose_2023}. Indeed, this operator admits a parametrix (approximate inverse, modulo compact error) lying in the calculus of \emph{0-pseudodifferential operators}. This parametrix construction plays a key role in the characterization of $\mathrm{spec} (\Delta_g)$. The authors of \cite{dHIKM_gas-giants_2024} remarked that a more precise characterization of the truncated spectra could be obtained using a uniform parametrix construction for the family $\{ \Delta_\epsilon \}$.

In this paper, we demonstrate such a construction. We introduce a calculus of \emph{\ztrtext{} pseudodifferential operators} and show that $\Delta_\epsilon$ admits a family of parametrices that lies in this calculus. In particular, the parametrices admit well-behaved asymptotic expansions. Using this, we study the uniform behavior of the resolvent $(\Delta_\epsilon - \lambda)^{-1}$. As a result, we show that simple eigenvalues admit asymptotic expansions in $\epsilon$ and also characterize the joint asymptotics of the corresponding eigenvectors in $(\epsilon, x)$.

The use of transition calculi to study families of degenerating operators, such as the Laplacian or Dirac operators associated to a degenerating family of metrics, is ubiquitous. Examples of this include \cite{mcdonald_laplacian-cone-singularities_1990} \cite{mazzeo-melrose_analytic-surgery_1995} \cite{guillarmou-hassell_low-energy-resolvent-asymptotically-conic_2007} \cite{albin-rochon-sher_fibered-cusps_2021} \cite{anghel_dirac-spectrum-degenerating-riemann-surfaces_2024} \cite{albin_cheegermuller_2022}.

There are several motivations for using kernel-based methods and a transition calculus to study the spectrum of $\Delta_\epsilon$. Obtaining asymptotic expansions for the kernel of the resolvent yields asymptotic expansions of eigenvalues and eigenvectors with relative simplicity. Moreover, asymptotic expansions of geometric operators tend to be fruitful sources of invariants in spectral geometry. Also, although we work exclusively with weighted $L^2$ Sobolev spaces, pseudodifferential operators have good mapping properties on a variety of spaces, including $L^p$ spaces and H\"{o}lder spaces.

Another advantage of developing a transition calculus is that our techniques could be generalized to study other families of operators besides the Laplacian. For example, the wave operator associated to a gas giant metric plays an important role in gas giant modeling \cite{prat_asymptotic_2016} \cite{dewberry_constraining_2021}. Moreover, although the spectrum of $\Delta_g$ is no longer discrete in the asymptotically hyperbolic case (corresponding to $\alpha = 2$) \cite{mazzeo_hodge_1988}, our \ztrtext{} calculus should be equally applicable to the study of families of metrics degenerating to an asymptotically hyperbolic metric, because $\Delta_g$ is a genuine $0$-operator in this case.

\section*{On artificial intelligence}

This paper is the outcome of a project that commenced in April 2024. Artificial intelligence models began to be used in the research process from January 2026, including Google Gemini (versions 3-3.5) and ChatGPT (versions GPT-5.5 to GPT-6.1). Wolfram Mathematica was used for computer algebra computations, with agentic programming assistance from OpenAI Codex (GPT-5.5 to 6.1). An AI model (GPT-6 to 6.1) generated the figures in this paper, with the exception of Figures \ref{fig_commutative-diagram-double-and-triple} and \ref{fig_triple-space-submanifolds-poset}. Some proofreading was done by AI models. No writing was done by AI models.

\section{\ztrtext{} differential operators}

Let $M$ be an $n$-dimensional compact manifold with boundary. Let $(x, y)$ be local coordinates on $M$, where $x$ is a boundary-defining function and $y$ restricts to a coordinate system on $\partial M$. We define a ring of \emph{\ztrtext{} differential operators} on $M$. Such an operator is really a family of differential operators parametrized by $\epsilon \in [0,1)$, which includes, up to a prefactor, the family of Laplacians associated to truncations of a gas giant.

Let $\sgsp{M}{\gas} = ([0,1) \times M)_b = [[0,1) \times M; \{0\} \times \partial M]$. We will frequently use the notation $\anysp{M}{0,1}$ to denote $[0,1) \times M$ and $\anysp{M}{01}$ to denote $\sgsp{M}{\gas}$. The index $0$ in $M_{01}$ refers to the factor of $[0,1)$ and the index $1$ refers to the copy of $M$. Comma-separated indices represent a Cartesian product, whereas concatenation of indices without commas indicates a blow-up of a product. This notation will be useful when it comes time to construct blow-ups involving multiple copies of $M$ and it becomes necessary to index them.

Let the boundary hypersurfaces of $\anysp{M}{01}$ be $\bf{0}$ (the proper lift of $\{0\} \times M$), $\bf{1}$ (the proper lift of $[0,1) \times \partial M$), and $\tf[1]$ (the front face). Let $\beta^{01}_{0,1} : \sgsp{M}{\gas} \to [0,1) \times M$ be the blow-down map and $\beta^{0,1}_0 : [0,1) \times M \to [0,1), \beta^{0,1}_1 : [0,1) \times M \to M$ be projections. Then define $\beta^{01}_0 : \anysp{M}{01} \to [0,1)$ by $\beta^{01}_0 = \beta^{0,1}_0 \circ \beta^{01}_{0,1}$. Let $\mathcal{V}_\gas(\anysp{M}{01})$ be the space of vector fields on $\anysp{M}{01}$ tangent to the fibres of $\beta^{01}_0$, vanishing at $\bf{1}$, and tangent at $\tf[1]$ to the fibres of the projection $\tf[1] \to \partial M$. This is generated by lifts along $\beta^{01}_1$ of vector fields from $\mathcal{V}_0(M)$ (which is generated by $x \partial_x$ and $x \partial_y$).

We define the ring of \ztrtext{} differential operators $\mathrm{Diff}^*_\gas(M)$ to be generated by $\mathcal{V}_\gas(\anysp{M}{01})$ with coefficients in $C^\infty(\anysp{M}{01})$. Equivalently, it’s generated over $C^\infty(\anysp{M}{01})$ by $x \partial_x$ and $x \partial_y$. Define $\mathrm{Diff}^m_\gas(M)$ to consist of the operators of order $m$.

\section{The \ztrtext{} double space}
The \ztrtext{} double space $\dbsp{M}{\gas}$ is constructed as a blow-up of $[0,1) \times M^2$. Label the boundary hypersurfaces of $[0,1) \times M^2$ as $\bf{0} = \{0\} \times M^2, \bf{1} = [0,1) \times \partial M \times M, \bf{2} = [0,1) \times M \times \partial M$. Let $(x_i, y_i)$ be a coordinate system for the $i$th copy of $M$ and let $\epsilon$ parametrize $[0,1)$. Then $(\epsilon, x_1, x_2, y_1, y_2)$ is a coordinate system on $[0,1) \times M^2$, wherein $\bf{i} = \{ x_i = 0 \}$ and $\bf{0} = \{ \epsilon = 0 \}$.

We will frequently use the shorthand $\anysp{M}{0,1,2}$ to denote $[0,1) \times M^2$ and $\anysp{M}{012}$ to denote $\dbsp{M}{\gas}$. The index $0$ refers to the factor of $[0,1)$ and the indices $1,2$ refer to copies of $M$. This notation is helpful whenever it is important to distinguish the two copies of $M$.

Define the following submanifolds of $\anysp{M}{0,1,2}$:
\begin{align}
    W^{12} &= \{ \epsilon = x_1 = x_2 = |y_1-y_2| = 0\} \nonumber \\
    C^{12} &= \{ \epsilon = x_1 = x_2 = 0 \} \nonumber \\
    S^{12} &= \{ x_1 = x_2 = |y_1-y_2| = 0 \} \nonumber \\
    B^{12} &= \{ x_1 = x_2 = 0 \} \nonumber \\
    C^1 &= \{ \epsilon = x_1 = 0 \} \nonumber \\
    C^2 &= \{ \epsilon = x_2 = 0 \} \label{eqn_double-space-blow-up-centers}
\end{align}

We also frequently abuse notation and use the same symbols to denote the lifts of these submanifolds to blow-ups of $\anysp{M}{0,1,2}$. Each of these is a $p$-submanifold of $\anysp{M}{0,1,2}$. Now define $\anysp{M}{012} = [\anysp{M}{0,1,2};W^{12};C^{12};S^{12};B^{12};C^1;C^2]$.

%-- Updated 1 Sep 26 to incorporate new face bff
Denote the front faces of these blow-ups, in order, by $\tff$, $\btff$, $\zff$, $\bff$, $\tf[1]$, and $\tf[2]$. The corner structure of $\anysp{M}{012}$ can be represented by an abstract simplicial complex (family of sets closed under taking subsets) on the boundary hypersurfaces of $\anysp{M}{012}$ (three lifts of boundary hypersurfaces present in $\anysp{M}{0,1,2}$ and six produced by blowing up). This simplicial complex consists of all sets of boundary hypersurfaces whose intersection is nonempty. The complex has maximal sets $\{ \bf{0}, \tff$, $\btff, \tf[1] \}$, $\{ \bf{0}, \tff, \btff, \tf[2] \}$, $\{ \bf{1}, \tff, \btff, \tf[1] \}$, $\{ \bf{1}, \bff, \tff, \btff \}$, $\{ \bf{2}, \bff, \tff, \btff \}$, $\{ \bf{2}, \tff, \btff, \tf[2] \}$, $\{ \bf{1}, \bff, \tff, \zff \}$, and $\{ \bf{2}, \bff, \tff, \zff \}$.

\subsection{Regions of $\dbsp{M}{\gas}$}

To understand the space $\anysp{M}{012} = \dbsp{M}{\gas}$, it is helpful to separate it into three regions. Observe that the first blow-up space in the sequence of blow-ups, $[\anysp{M}{0,1,2}; W^{12}]$, has four boundary hypersurfaces: $\bf{i}$ for $i = 0,1,2$, in addition to $\tff$. The intersection of all four of these boundary hypersurfaces is nonempty. That is, the intersection complex of $[\anysp{M}{0,1,2}; W^{12}]$ is a $3$-simplex on $\{ \bf{0}, \bf{1}, \bf{2}, \tff\}$.

Also note that the diagonal $\{ x_1 = x_2, y_1 = y_2 \}$ lifts to a $p$-submanifold of $[\anysp{M}{0,1,2}; W^{12}]$. We denote the lifted diagonal by $\anydiag{12}$; we will abuse notation by using the same notation to denote the lift of the diagonal to $\anysp{M}{012}$. Note that, in $[\anysp{M}{0,1,2}; W^{12}]$, the intersection $\bf{0} \cap \bf{1} \cap \bf{2} \cap \anydiag{12}$ is empty.

\begin{figure}
    \centering
    \includegraphics[width=0.5\linewidth]{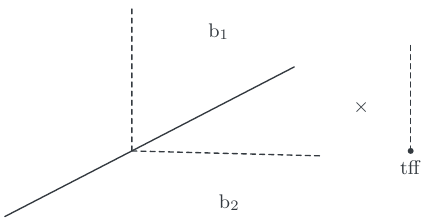}
    \caption{The corner structure of the upper region of $[\anysp{M}{0,1,2}; W^{12}]$.}
    \label{fig_double-space-upper-region-blowndown}

    \centering
    \includegraphics[width=0.5\linewidth]{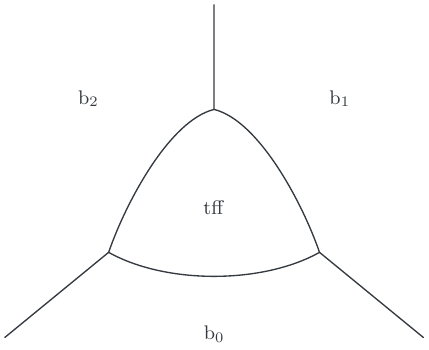}
    \caption{The corner structure of the middle region of $[\anysp{M}{0,1,2}; W^{12}]$.}
    \label{fig_double-space-middle-region-blowndown}

    \centering
    \includegraphics[width=0.5\linewidth]{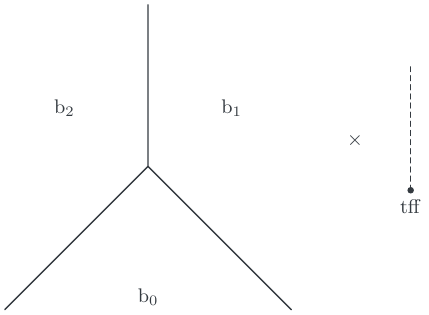}
    \caption{The corner structure of the lower region of $[\anysp{M}{0,1,2}; W^{12}]$.}
    \label{fig_double-space-lower-region-blowndown}
\end{figure}

Now define $U_\mathrm{up} = [\anysp{M}{0,1,2}; W^{12}] \setminus \bf{0}$, $U_\mathrm{mid} = [\anysp{M}{0,1,2}; W^{12}] \setminus \big( (\bf{0} \cup \anydiag{12}) \cap (\bf{1} \cap \bf{2}) \big) = [\anysp{M}{0,1,2}; W^{12}] \setminus (C^{12} \cup S^{12})$, and $U_\mathrm{low} = [\anysp{M}{0,1,2}; W^{12}] \setminus \anydiag{12}$.  Because $\bf{0} \cap \bf{1} \cap \bf{2} \cap \anydiag{12}$ is empty, $U_\mathrm{up} \cup U_\mathrm{mid} \cup U_\mathrm{low} = [\anysp{M}{0,1,2}; W^{12}]$. Figures \ref{fig_double-space-upper-region-blowndown}, \ref{fig_double-space-middle-region-blowndown}, and \ref{fig_double-space-lower-region-blowndown} depict the corner structures of these regions.

\begin{figure}
    \centering
    \includegraphics[width=0.5\linewidth]{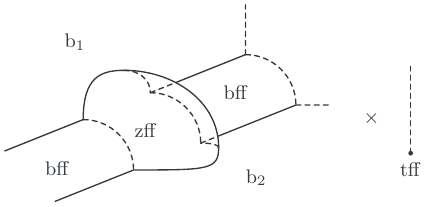}
    \caption{The corner structure of the upper region of $\anysp{M}{012} = \dbsp{M}{\gas}$.}
    \label{fig_double-space-upper-region}

    \centering
    \includegraphics[width=0.5\linewidth]{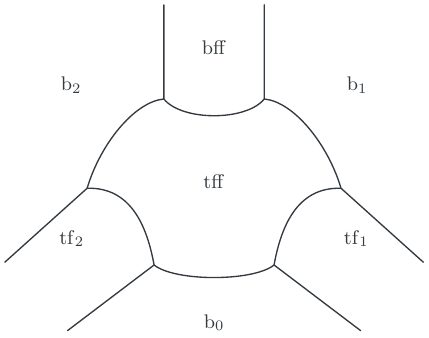}
    \caption{The corner structure of the middle region of $\anysp{M}{012} = \dbsp{M}{\gas}$.}
    \label{fig_double-space-middle-region}

    \centering
    \includegraphics[width=0.5\linewidth]{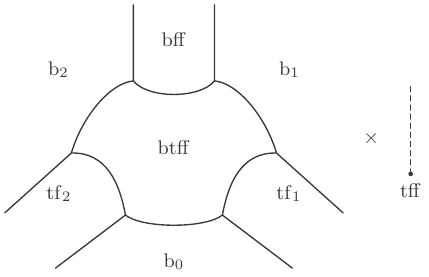}
    \caption{The corner structure of the lower region of $\anysp{M}{012} = \dbsp{M}{\gas}$.}
    \label{fig_double-space-lower-region}
\end{figure}

Lift these regions to $\anysp{M}{012}$. We shall use the same notation for the lifts and, in fact, we shall always assume unless specified otherwise that $U_\ast$ refers to a region of $\anysp{M}{012}$. The corner structures of these regions are depicted in Figures \ref{fig_double-space-upper-region}, \ref{fig_double-space-middle-region}, and \ref{fig_double-space-lower-region}.

Let us describe the structure of each region in more detail and introduce some useful projective coordinate systems. Recall that $(\epsilon, x_1, x_2, y_1, y_2)$ is a coordinate system on $\anysp{M}{0,1,2} = [0,1) \times M^2$. Define $s = x_1/x_2$, $\rho = \epsilon/x_2$, $r = x_1/\epsilon$. Also define $s' = s^{-1}$, $\rho' = \rho^{-1}$, and $r' = r^{-1}$. Also define $Y = y_1 - y_2$, $u = Y/x_2$, $v = Y/x_1$, and $w = Y/\epsilon$.

\begin{figure}
    \centering
    \includegraphics[width=0.5\linewidth]{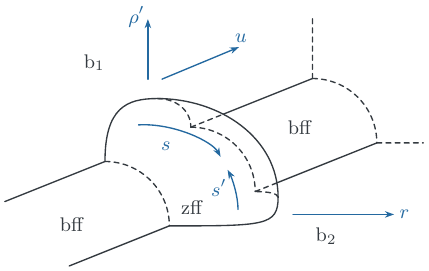}
    \caption{Local coordinate systems on the upper region of $\anysp{M}{012}$. This diagram suppresses a factor of $\R^{n}_1$ parametrized by $(\epsilon, y_2)$.}
    \label{fig_double-space-upper-region-coords}
\end{figure}

Away from $\bf{2}$ and $\bff$, $U_\mathrm{up}$ admits a coordinate chart $(\epsilon, s, u, \rho', y_2)$. Away from $\bf{1}$ and $\bff$, $U_\mathrm{up}$ admits a chart $(\epsilon, s', v, r, y_2)$. (See the diagram in Figure \ref{fig_double-space-upper-region-coords}.) The upper region can be identified with the product of an interval with an `overblown' $0$-double space \cite{hintz_0-calculus-mazzeo-melrose_2023}. Let $\tilde{M} = \R^n_1 = [0,\infty) \times \R^{n-1}$ be a local model for $M$, parametrized by $(\tilde{x}, \tilde{y})$. (Other tilded coordinates are defined in terms of $\tilde{x}, \tilde{y}$ in the same way their non-tilded counterparts are defined in terms of $x, y$.) The $0$-double space $\dbsp{\tilde{M}}{0'}$ is, away from its b-front face, covered by coordinate charts $(\tilde{s}, \tilde{u}, \tilde{x}_2, \tilde{y}_2)$ and $(\tilde{s}', \tilde{v}, \tilde{x}_1, \tilde{y}_2)$. Let $\tilde{\epsilon}$ parametrize the unit interval $[0,1)$. The mapping defined equivalently by $(\tilde{\epsilon}, \tilde{s}, \tilde{u}, \tilde{x}_2, \tilde{y}_2) = (\epsilon, s, u, \rho', y_2)$ or by $(\tilde{\epsilon}, \tilde{s}', \tilde{v}, \tilde{x}_1, \tilde{y}_2) = (\epsilon, s', v, r, y_2)$ is a (locally defined) map from $U_\mathrm{up}$ to $[0,1) \times \dbsp{\tilde{M}}{0'}$ that is a diffeomorphism onto its image, in the sense of manifolds with corners.

\begin{figure}
    \centering
    \includegraphics[width=0.5\linewidth]{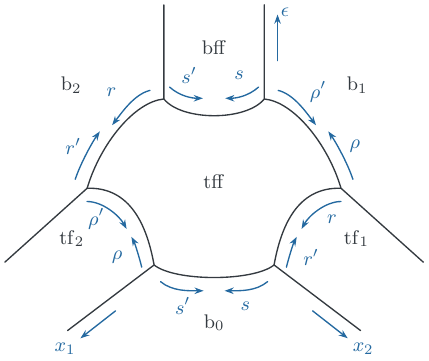}
    \caption{Local coordinate systems on the middle region of $\anysp{M}{012}$. This diagram suppresses a factor of $\R^{2(n-1)}$ parametrized by $(u,y_2)$, $(v,y_2)$, or $(w,y_2)$.}
    \label{fig_double-space-middle-region-coords}
\end{figure}

The middle region $U_\mathrm{mid} \subset \anysp{M}{012} = \dbsp{M}{\gas}$ is covered by six coordinate systems. The chart $(\rho, r, u, x_2, y_2)$ is smooth up to $\tff \cup \bf{1} \cup \tf[1]$. The chart $(\rho, s', v, x_1, y_2)$ is smooth up to $\tff \cup \bf{0} \cup \tf[2]$. The chart $(r', s, u, x_2, y_2)$ is smooth up to $\tff \cup \bf{0} \cup \tf[1]$. The chart $(r', \rho', v, x_1, y_2)$ is smooth up to $\tff \cup \bf{2} \cup \tf[2]$. The chart $(s', r, w, \epsilon, y_2)$ is smooth up to $\tff \cup \bff \cup \bf{2}$. The chart $(\rho', s, w, \epsilon, y_2)$ is smooth up to $\tff \cup \bff \cup \bf{1}$. See Figure \ref{fig_double-space-middle-region-coords} for a diagram of these coordinates.

\begin{figure}
    \centering
    \includegraphics[width=0.5\linewidth]{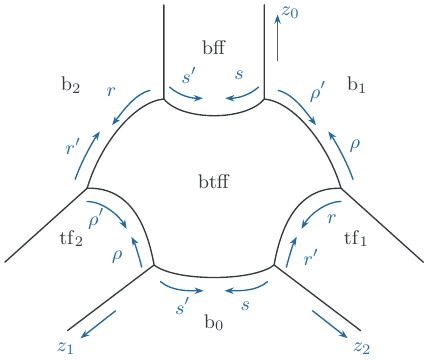}
    \caption{Local coordinate systems on the lower region of $\anysp{M}{012}$. This diagram suppresses a factor of $\R^{2(n-1)}$ parametrized by $(Y, y_2)$.}
    \label{fig_double-space-lower-region-coords}
\end{figure}

The lower region $U_\mathrm{low} \subset \anysp{M}{012} = \dbsp{M}{\gas}$ is covered by six coordinate systems. Let $(z_0, z_1, z_2) = (\epsilon, x_1, x_2)/|Y|$. The charts are $(\rho, r, z_2, Y, y_2), (r', s, z_2, Y, y_2)$, $(\rho, s', z_1, Y, y_2)$, $(r', \rho', z_1, Y, y_2)$, $(s', r, z_0, Y, y_2)$, and $(\rho', s, z_0, Y, y_2)$. (See Figure \ref{fig_double-space-lower-region-coords}.)

\subsection{Boundary hypersurfaces}

We comment here on the boundary hypersurfaces of $\anysp{M}{012}$. The boundary hypersurfaces $\zff, \tff$, and $\bf{0}$ intersect $\diaglift$. The face $\zff$ is referred to as the `$0$-front face'. From the local identification of $\dbsp{M}{\gas}$ with $[0,1) \times \dbsp{\tilde{M}}{0'}$, it is immediate that $\zff$ can be locally identified with $[0,1) \times \mathrm{ff}( \dbsp{\tilde{M}}{0'} )$, but in fact it is canonically diffeomorphic to $[0,1) \times \mathrm{ff}( \dbsp{M}{0'} )$. The face $\tff$ is referred to as the `transition front face' or `0-transition front face', as it describes the transition of a \ztrtext{} operator from the $\epsilon > 0$ regime to the $\epsilon = 0$ regime. Locally in the upper region, $\tff \cong \dbsp{\tilde{M}}{0'}$. The face $\bf{0}$, namely the lift of $\{\epsilon = 0\}$ to $\anysp{M}{012}$, is referred to as the `lower face' or `bottom face'. Locally, it is diffeomorphic to the overblown double space $\dbsp{\tilde{M}}{0'}$ \cite{hintz_0-calculus-mazzeo-melrose_2023}.

\begin{figure}
    \centering
    \includegraphics[width=0.5\linewidth]{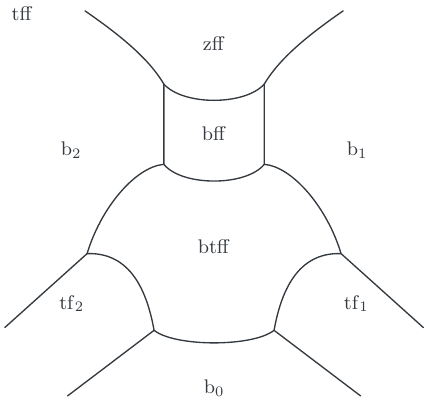}
    \caption{The corner structure of $\tff$.}
    \label{fig_tff-corner-structure}
\end{figure}

Of these front faces, the transition front face is most interesting to us as it captures the transitional behavior of an operator as $\epsilon \to 0$. It intersects every other boundary hypersurface of $\dbsp{M}{\gas}$. Figure \ref{fig_tff-corner-structure} depicts its corner structure.

\subsection{Some notation}

Let $\beta^{012}_{0,1,2} : \dbsp{M}{\gas} = \anysp{M}{012} \to [0,1) \times M^2 = \anysp{M}{0,1,2}$ be the blow-down map. In general, given indices $\sigma, \sigma'$, the symbol $\beta^{\sigma}_{\sigma'}$ will be used to denote a natural map $M_{\sigma} \to M_{\sigma'}$. There is usually only one reasonable choice for this map, built from blow-downs and projections.

Let $\lambda$ be a nonsingular, nondegenerate half-density on $[0,1) \times M^2$. Locally, we could set $\lambda = (d\epsilon dx_1 dx_2 dy_1 dy_2)^{1/2}$. Where there is no chance for misinterpretation, we will also use $\lambda$ to denote its own lift along $\beta^{012}_{0,1,2}$.

Given a boundary hypersurface $\mathrm{h}$, write $\varrho_\mathrm{h}$ for a boundary-defining function for $\mathrm{h}$.

\section{The \ztrtext{} calculus}

A \ztrtext{} operator is a family of pseudodifferential operators parametrized by $\epsilon$. It can be represented by a half-density distribution on $\anysp{M}{012}$ that is conormal with respect to $\diaglift$ and polyhomogeneous at the boundary hypersurfaces of $\anysp{M}{012}$.

It is natural to normalize \ztrtext{} operators with respect to a certain half-density $\sigma$. Define $\sigma = \varrho_{\uff}^{-n} \varrho_{\tff}^{-n} \left(\beta^{012}_{0,1,2}\right)^* \lambda$. The reason is that the normalization of the kernel of the identity operator with respect to $\sigma$ is a product of delta functions without degenerate or singular factors. Define the `\ztrtext{} half-density bundle' ${}^\gas\Omega^{1/2}(\anysp{M}{012})$ to be trivialized by $\sigma$.

Recall that $I^*(X,Y;E)$ denotes the space of conormal distributions on the manifold $X$ with respect to the submanifold $Y$, taking values in the bundle $E$ \cite{hormander_partial-differential-operators-III_2007}.

\begin{definition}

The small \ztrtext{} calculus is
\begin{align*}
    \Psi^*_\gas \left( M\right) = \mathcal{A}^{\mathcal{E}_0}_\mathrm{phg}(\dbsp{M}{\gas}) I^*\left( \dbsp{M}{\gas}, \diaglift; {}^\gas\Omega^{1/2} \right)
\end{align*}
The index family $\mathcal{E}_0$ on $\dbsp{M}{\gas}$ is generated by $(0,0)$ at the faces intersecting the diagonal -- namely $\zff, \tff$, and $\bf{0}$ -- and $\emptyset$ at other faces.

\end{definition}
The calculus is filtered by the spaces $\Psi^m_\gas \left( M\right) = \mathcal{A}^{\mathcal{E}_0}_\mathrm{phg}(\dbsp{M}{\gas}) I^m\left( \dbsp{M}{\gas}, \diaglift; {}^\gas\Omega^{1/2} \right)$, where the order of conormal distributions is defined so that operators have the expected order, i.e. so that an $m$th order differential operator lies in $\Psi^m_\gas$.

A differential operator $P = \sum_{j \gamma} a_{j, \gamma}(\epsilon,x,y) (x\partial_x)^j (x\partial_y)^\gamma \in \mathrm{Diff}^*_\gas(M)$ is identified with a kernel (unique only up to a conformal factor)
\begin{align}
    K_P = \sum_{j+|\gamma|\leq m} a_{j, \gamma}(\epsilon,x_1,y_1) (x_1\partial_{x_1})^j (x_1\partial_{y_1})^\gamma \delta(x_1-x_2) \delta(y_1-y_2) \lambda \label{eqn_kernel-gas-diff-op}
\end{align}

As a conormal distribution, an element $Q$ of the small calculus has a principal symbol $\sigma(Q)$, valued in the following symbol space:
\begin{align*}
    S^{\{m\}}_\gas(M) := 
S^{\{m\}}(N^* \diaglift) \otimes \Omega_\mathrm{fibre}(N^* \diaglift) \otimes {}^\gas\Omega^{1/2}(\anysp{M}{012})|_{\diaglift}
\end{align*}
for some $m$. Here, $\diaglift \subset \anysp{M}{012}$ is the lifted diagonal. The superscript $\{m\}$ indicates a quotient by symbols of order one lower.

The order of the principal symbol gives a filtration of $\Psi^*_\gas \left( M\right)$ by subsets $\Psi^m_\gas \left( M\right)$ consisting of the operators of order $m$.

In fact, this symbol space can be identified with symbolic sections of a pullback of the $0$-cotangent bundle from $M$ to $\sgsp{M}{\gas}$.

\begin{lemma}
\label{lemma_symbol-space}

\begin{align*}
    S^{\{m\}}_\gas(M) &\cong S^{\{m\}}(N^* \diaglift) \\
    &\cong S^{\{m\}} \left( \left( \beta^{01}_1 \right)^*\left({}^0T^*M\right) \right)
\end{align*}
where $\beta^{01}_1 : \sgsp{M}{\gas} \to M$ is the blow-down followed by a projection.

\end{lemma}

\begin{proof}
    First observe that there is a natural diffeomorphism $\diaglift \cong \sgsp{M}{\gas}$. To show this, observe that the diagonal in $[0,1) \times M^2$ is naturally diffeomorphic to $[0,1) \times M$. The diffeomorphism identifies $\diaglift \cap W^{12}$ with $\{0\} \times \partial M$. Therefore, the lift of $\diaglift$ to $[ [0,1) \times M^2 ; W^{12}]$ is identified with $\sgsp{M}{\gas}$. Besides $W^{12}$, all other submanifolds defined in \eqref{eqn_double-space-blow-up-centers} are either disjoint from the lifted diagonal or lift to one of its boundary hypersurfaces, meaning their blow-up does not change $\diaglift$. Thus $\diaglift \subset \dbsp{M}{\gas}$ is identified with $\sgsp{M}{\gas}$.

    The rest of the argument is standard; see \cite[4.3]{mazzeo-melrose_analytic-surgery_1995}.
\end{proof}

The quantization map $\mathrm{Op}$ that produces a small-calculus operator (modulo lower-order operators) from a principal symbol defines an isomorphism $S^{\{m\}}_\gas(M) \to \Psi^{\{m\}}_\gas(M)$ (where again curly braces indicate a quotient by operators of order one lower).

Now we will define the full \ztrtext{} calculus.

\begin{definition}
    Suppose $\mathcal{E}$ is an index family on $\dbsp{M}{\gas}$. Then the set of residual \ztrtext{} pseudodifferential with index family $\mathcal{E}$ is
    \begin{align*}
        \Psi^{-\infty, \mathcal{E}}_\gas \left( M \right) = \mathcal{A}^\mathcal{E}_\mathrm{phg}\left( \dbsp{M}{\gas}; {}^\gas \Omega^{1/2} \right).
    \end{align*}
    If $m \in \R$, then the set of \ztrtext{} pseudodifferential operators of order $m$ with index family $\mathcal{E}$ is
    \begin{align*}
        \Psi^{m, \mathcal{E}}_\gas \left( M \right) = \mathcal{A}^\mathcal{E}_\mathrm{phg}(\dbsp{M}{\gas}) I^m\left( \dbsp{M}{\gas}, \diaglift; {}^\gas\Omega^{1/2} \right)
    \end{align*}
\end{definition}

Let $\beta^{012}_{0,i} : \dbsp{M}{\gas} \to [0,1) \times M$ be defined by $\beta^{012}_{0,i} = \beta^{0,1,2}_{0,i} \circ \beta^{012}_{0,1,2}$, where $\beta^{0,1,2}_{0,i} : [0,1) \times M^2 \to \anysp{M}{0,i} = [0,1) \times M$ is a projection. Let $\beta^{012}_i : \dbsp{M}{\gas} \to M$ be defined by $\beta^{012}_i = \beta^{0,i}_i \circ \beta^{012}_{0,i}$, where $\beta^{0,i}_i : M_{0,i} = [0,1) \times M \to M$ is a projection.

By lifting a half-density on $M$ to $\dbsp{M}{\gas}$ along $\beta^{012}_2$, multiplying by a kernel, and pushing forward along $\beta^{012}_{0,1}$, a \ztrtext{} pseudodifferential operator maps half-densities on $M$ to half-densities on $M \times [0,1)$. However, $\beta^{012}_{0,1}$ is not a b-fibration.

Instead, consider the map $\beta^{012}_{01} : \dbsp{M}{\gas} \to \sgsp{M}{\gas}$. This is the unique b-map such that $\beta^{012}_{0,1}$ factors as $\beta^{01}_{0,1} \circ \beta^{012}_{01}$. It is a b-fibration; its lifting matrix is
\begin{align}
    e^{012}_{01} = \begin{pmatrix}
        1 & 0 & 0 & 0 & 0 & 0 & 0 & 0 & 1 \\
        0 & 1 & 0 & 0 & 0 & 1 & 1 & 0 & 0 \\
        0 & 0 & 0 & 1 & 1 & 0 & 0 & 1 & 0
    \end{pmatrix} \label{eqn_beta-012-01-lifting-matrix}
\end{align}
where the rows correspond to the faces $\bf{0}$, $\bf{1}$, and $\tf[1]$ of $\sgsp{M}{\gas}$ and the columns correspond to the faces $\bf{0}$, $\bf{1}$, $\bf{2}$, $\tff$, $\btff$, $\zff$, $\bff$, $\tf[1]$, and $\tf[2]$ of $\dbsp{M}{\gas}$. It is thus more appropriate, a priori, to define a \ztrtext{} operator as a map from half-densities on $M$ to half-densities on $\sgsp{M}{\gas}$.

\begin{definition}
    A kernel $K_Q \in \Psi^*_{\gas}$ defines an operator $Q: \dot{C}^\infty(M, \Omega^{1/2}) \to \dot{C}^\infty(M_{01}, \Omega^{1/2})$ given by
    \begin{align}
        Qu &:= (\beta^{012}_{01})_* \left( K_Q \left(\beta^{012}_2\right)^*u \right) \label{eqn_op-mapping-M-to-single-space}
    \end{align}
    for all $u \in C^\infty(M, \Omega^{1/2}(M))$.
    
    We choose to normalise $K_Q$ with respect to the density $\sigma$ and write $K_Q = \mathcal{K}_Q \sigma$, where is a scalar distribution. Then $Qu = (\beta^{012}_{01})_* \left( \mathcal{K}_Q \sigma \left(\beta^{012}_2\right)^*u \right)$.
\end{definition}

Of course, we can also identify $Q$ with a family of operators $Q_\epsilon : \dot{C}^\infty(M, \Omega^{1/2}) \to \dot{C}^\infty(M, \Omega^{1/2})$ by setting $Q_{\epsilon_0} u = (Qu)/(d\epsilon)^{1/2}|_{\epsilon = \epsilon_0}$ for $\epsilon_0 \in (0,1)$.

We can also identify $Q$ as an operator $Q: \dot{C}^\infty(\sgsp{M}{\gas}, \Omega^{1/2}) \to \dot{C}^\infty(\sgsp{M}{\gas}, \Omega^{1/2})$ given by
\begin{align*}
    Qu &:= (\beta^{012}_{01})_* \left( K_Q \left(\beta^{012}_{02} \right)^*u \right)
\end{align*}
where $\beta^{012}_{02} : \dbsp{M}{\gas} \to \sgsp{M}{\gas}$ is defined analogously to $\beta^{012}_{02}$ by exchanging the roles of the first and second copies of $M$.

\section{Mapping properties}

From the index family of a \ztrtext{} operator, we can read off its mapping properties. We introduce some notation for index families. Suppose $E$ is an index set and $s = \min(\{ \mathrm{Re}(z) : (z,p) \in E \text{ for some } p \in \N_0 \})$. For $c \in \R$, we write $E \geq c$ to mean that either $s > c$ or $s = c$ and $\not \exists (z,p) \in E$ such that $\mathrm{Re}(z) = s$ and $p > 0$. We write $E > c$ to mean that $s > c$. In particular, if $E$ is the index set for a polyhomogeneous function $f$ at a boundary hypersurface $\mathrm{h}$, $E \geq 0$ means that $f$ is bounded up to $\mathrm{h}$ and $E > 0$ means that $f$ vanishes at $\mathrm{h}$.

Kernels in $\Psi^{\ast,\ast}_{\gas}(M)$ define families of bounded operators between appropriate $0$-Sobolev spaces. We will only need to know these mapping properties for residual operators.

\begin{proposition}
\label{proposition_gas-giant-operator-boundedness}
\allowdisplaybreaks[0]

    Suppose $Q \in \Psi^{-\infty, \mathcal{E}}_\gas(M)$. (That is, $K_Q = \mathcal{K}_Q \sigma$ and $\mathcal{E}$ is the index family of $\mathcal{K}_Q$.) Given a boundary hypersurface $\mathrm{h}$ of $\anysp{M}{012}$, let $E_\mathrm{h}$ be the index set at $\mathrm{h}$. Suppose that
    \begin{align}
        E_{\bf{0}} & \geq 0 \nonumber \\
        E_{\bf{1}}, E_{\bf{2}} &> -\frac12 \nonumber \\
        E_{\tff}, E_{\zff} &\geq 0 \nonumber \\
        E_{\btff}, E_{\bff} &\geq -1 \nonumber \\
        E_{\tf[1]}, E_{\tf[2]} &\geq -\frac12 \label{eqn_index-criteria-for-boundedness}.
    \end{align}
    Moreover, suppose also one of the following:
    \begin{itemize}
        \item $E_{\tf[1]}, E_{\tf[2]} > -1/2$, either $E_{\btff} > -1$ or $E_{\tff} > 0$, and either $E_{\bff} > -1$ or $E_{\zff} > 0$.
        % \item $E_{\tf[1]}, E_{\tf[2]} > -1/2$ and one of the following two hold: $E_{\btff} > -1$ or $E_{\tff} > 0$.
        \item $E_{\btff} > -1$, $E_{\bf{0}} > 0$, $E_{\tff} > 0$, and either $E_{\bff} > -1$ or $E_{\zff} > 0$.
        % \item $E_{\btff} > -1$, $E_{\bf{0}} > 0$, and $E_{\tff} > 0$.
    \end{itemize}
    Then, for any $k, \ell \in \R$, $Q$ defines a family of bounded operators $Q_\epsilon : H^k_0(M, \Omega^{1/2}) \to H^\ell_0(M, \Omega^{1/2})$, with uniformly bounded operator norm, where $H^k_0(M, \Omega^{1/2})$ is a $0$-Sobolev space \cite{mazzeo_differential-edge-operators-I_1991} \cite{hintz_0-calculus-mazzeo-melrose_2023}.

    Moreover, if the first list of inequalities is modified so that all are strict, except possibly for the inequalities $E_{\zff} \geq 0$ and $E_{\bff} \geq -1$, then $\lim_{\epsilon \to 0} \| Q_\epsilon \|_{x^\delta L^2} = 0$ for some $\delta > 0$.
\end{proposition}

\begin{proof}
    By commuting elliptic $0$-operators through $Q$, we reduce to the case $k = \ell = 0$; i.e. we shall prove $Q_\epsilon : L^2(M, \Omega^{1/2}) \to L^2(M, \Omega^{1/2})$ is bounded. It also suffices to consider the case where $Q$ has nonnegative kernel.

    Let $\gamma = (dxdy)^{1/2}$ be a trivializing density for $\Omega^{1/2}(M)$ and let $\mu = (d\epsilon dx dy)^{1/2}$ be a trivializing density for $\Omega^{1/2}([0,1) \times M)$. Let $f(x,y) = x^{-1/2} g(x,y) + (1-g(x,y))$, where $g \in C^\infty(M)$ is a cutoff function equalling $1$ wherever $x \leq 1/2$ and $0$ wherever $x \geq 1$. Of course, this is a polyhomogeneous function on $M$ with index set $(-1/2, 0) + \N_0$. 
    
    Let $f(x,y)$ be a positive polyhomogeneous function on $M$ with index set $(-1/2, 0) + \N_0$.

    Let $Q^*$ be the adjoint of $Q$ on $L^2(M, \Omega^{1/2})$. The index family of $Q^*$ is obtained from $\mathcal{E}$ by exchanging the two factors of $M$ in $\dbsp{M}{\gas}$.
    
    Using the pullback and pushforward theorems \cite{melrose_calculus-conormal-mwc_1992}, and the lifting matrix \eqref{eqn_beta-012-01-lifting-matrix}, formula \eqref{eqn_op-mapping-M-to-single-space}, and the index set criteria \eqref{eqn_index-criteria-for-boundedness}, one finds that $f^{-1} \mu^{-1} Q(f \gamma)$ and $f^{-1} \mu^{-1} Q^*(f \gamma)$ are polyhomogeneous on $\sgsp{M}{\gas}$ with nonnegative index sets at all boundary hypersurfaces, and therefore blow down to bounded functions on $[0,1) \times M$. Thus the families $\{ f^{-1} \gamma^{-1} Q_\epsilon (f \gamma) \}_{\epsilon \in [0,1)}$ and $\{ f^{-1} \gamma^{-1} Q_\epsilon^*(f \gamma) \}_{\epsilon \in [0,1)}$ are uniformly bounded on $M$. By Schur's test, this proves that $Q_\epsilon : L^2(M, \Omega^{1/2}) \to L^2(M, \Omega^{1/2})$ is a uniformly bounded family of operators. Moreover, we get the following estimate: $\| Q_\epsilon \|_{B(L^2)} \leq \sup_M (f^{-1} \gamma^{-1} Q_\epsilon (f \gamma))^{1/2} \sup_M (f^{-1} \gamma^{-1} Q_\epsilon^*(f \gamma))^{1/2}$. Again by the pullback and pushforward theorems, modifying the requisite inequalities to be strict precisely implies that the quantity on the right-hand side is $\mathcal{O}(\epsilon^c)$ for $c > 0$.
\end{proof}

\begin{proposition}
\allowdisplaybreaks[0]
\label{proposition_gas-giant-operator-compactness}

    Suppose $Q \in \Psi^{-\infty, \mathcal{E}}_\gas(M)$. Given a boundary hypersurface $\mathrm{h}$ of $\anysp{M}{012}$, let $E_\mathrm{h}$ be the index set at $\mathrm{h}$. Suppose that
    \begin{align*}
        E_{\bf{0}} &\geq 0 \\
        E_{\bf{1}}, E_{\bf{2}} &> -\frac12 \\
        E_{\tff}, E_{\zff} &> 0 \\
        E_{\btff}, E_{\bff} &> -1 \\
        E_{\tf[1]}, E_{\tf[2]} &> -\frac12.
    \end{align*}
    Then, there exists $\delta ' > 0$ such that, for any $k, \ell \in \R$, $Q$ defines a family of bounded operators $Q_\epsilon : H^k_0(M, \Omega^{1/2}) \to x^{\delta'} H^\ell_0(M, \Omega^{1/2})$ with uniformly bounded operator norm. (In particular, the $Q_\epsilon$ are compact.)

\end{proposition}

\begin{proof}
    Apply Proposition \ref{proposition_gas-giant-operator-boundedness} to the operator $x^{-\delta'} Q$. The index set criteria imply that $x^{-\delta'} Q$ is bounded for some sufficiently small $\delta'$. The operators $Q_\epsilon$ are compact because the inclusion $x^{\delta'} H^\ell_0(M, \Omega^{1/2}) \to H^m_0(M, \Omega^{1/2})$ is compact whenever $\ell > m$ and $\delta' > 0$ \cite{mazzeo_differential-edge-operators-I_1991}.
\end{proof}

We also need the mapping properties for residual \ztrtext{} operators acting on polyhomogeneous functions on $\sgsp{M}{\gas}$. The proof is a straightforward application of the pullback and pushforward theorems.

\begin{proposition}
\label{proposition_gas-giant-phg-mapping-properties}

    Suppose $Q \in \Psi^{-\infty, \mathcal{E}}_\gas(M)$. Let $u \in \mathcal{A}^{\mathcal{F}}_{\mathrm{phg}}(\sgsp{M}{\gas},\Omega^{1/2})$, where $\mathcal{F}$ is an index family on $\sgsp{M}{\gas}$. Then $Qu \in \mathcal{A}_{\mathrm{phg}}^{\mathcal{G}}(\sgsp{M}{\gas})$, where
    \begin{align*}
        G_{\bf{0}} &= (E_{\bf{0}} + F_{\bf{0}}) \overline{\cup} (E_{\tf[2]} + F_{\tf[1]} + 1) \\
        G_{\bf{1}} &= E_{\bf{1}} \overline{\cup} (E_{\zff} + F_{\bf{1}}) \overline{\cup} (E_{\bff} + F_{\bf{1}} + 1) \\
        G_{\tf[1]} &= (E_{\tff} + F_{\tf[1]}) \overline{\cup} (E_{\btff} + F_{\tf[1]} + 1) \overline{\cup} (E_{\tf[1]} + F_{\bf{0}}).
    \end{align*}
\end{proposition}

\section{Normal operators}
\label{section_normal-operators}

The normal operators of a \ztrtext{} differential or pseudodifferential operator are its restrictions to the front faces of $\dbsp{M}{\gas}$. A key step in the parametrix construction for a \ztrtext{} differential operator will be the inversion of its normal operators.

\subsection{Normal operators of differential operators}

\begin{figure}[h]
    \centering
    % https://q.uiver.app/#q=WzAsNCxbMCw0LCJNIl0sWzAsMiwiTV97MDF9Il0sWzAsMCwiTV97MDEyfSJdLFsxLDQsIk1fezAsMn0iXSxbMSwwLCJcXGJldGFeezAxfV8xIl0sWzIsMSwiXFxiZXRhXnswMTJ9X3swMX0iXSxbMiwwLCJcXGJldGFeezAxMn1fMSIsMix7ImN1cnZlIjo1fV0sWzIsMywiXFxiZXRhXnswMTJ9X3swLDJ9IiwwLHsiY3VydmUiOi00fV1d
    \[\begin{tikzcd}
    	{M_{012}} & \\
    	\\
    	{M_{01}} \\
    	\\
    	M & {M_{0,2}}
    	\arrow["{\beta^{012}_{01}}", from=1-1, to=3-1]
    	\arrow["{\beta^{012}_1}"', curve={height=30pt}, from=1-1, to=5-1]
    	\arrow["{\beta^{012}_{0,2}}", curve={height=-24pt}, from=1-1, to=5-2]
    	\arrow["{\beta^{01}_1}", from=3-1, to=5-1]
    \end{tikzcd}\]
    \caption{Commutative diagram showing the maps between the spaces $\anysp{M}{012} = \dbsp{M}{\gas}, \anysp{M}{01} = \sgsp{M}{\gas}$, and $\anysp{M}{0,2} = [0,1) \times M$.}
    \label{fig_commutative-diagram-double-and-triple}
\end{figure}
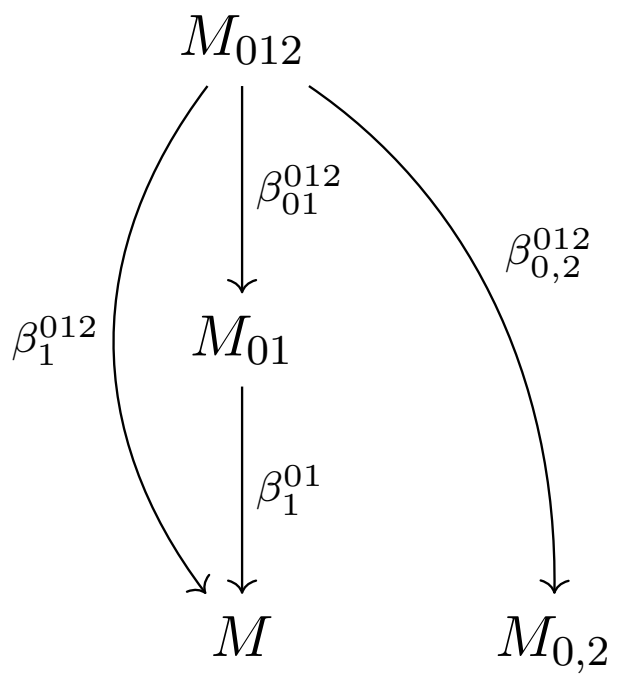

Recall that elements of $\mathrm{Diff}^*_\gas(M)$ are generated by $\mathcal{V}_\gas(\anysp{M}{01})$ with coefficients in $C^\infty(\anysp{M}{01})$. To make sense of the normal operators of an operator $P$, we must lift $\mathcal{V}_\gas(\anysp{M}{01})$ to the double space $\dbsp{M}{\gas} = \anysp{M}{012}$. Recall $\mathcal{V}_\gas(\anysp{M}{01})$ is generated by lifts along $\beta^{01}_1 : \anysp{M}{01} \to M$ of vector fields from $\mathcal{V}_0(M)$. Let $\mathcal{V}_\gas( \anysp{M}{012} )$ be the lift of $\mathcal{V}_\gas(\anysp{M}{01})$ to $\dbsp{M}{\gas}$ along the map $\beta^{012}_{01} : \anysp{M}{012} \to \anysp{M}{01}$ that projects out the second factor of $M$. Equivalently $\mathcal{V}_\gas( \anysp{M}{012} )$ is the set of vector fields generated, over $C^\infty(\anysp{M}{012})$, by lifts of vector fields in $\mathcal{V}_0(M)$ along the map $\beta^{012}_1 : \anysp{M}{012} \to M$ that projects onto the first copy of $M$. In other words, an element of $\mathcal{V}_\gas( \anysp{M}{012} )$ is tangent to the fibers of the projection $\beta^{012}_{0,2} : \anysp{M}{012} \to \anysp{M}{0,2} = [0,1) \times M$ (the map that projects out the first factor of $M$) and, moreover, its restriction to each fiber is $\beta^{012}_1$-related to an element of $\mathcal{V}_0(M)$. (See Figure \ref{fig_commutative-diagram-double-and-triple} for a diagram of the mappings involved.)

Because $0$-vector fields are b-vector fields and b-vector fields lift to b-vector fields, of course $\mathcal{V}_\gas(\anysp{M}{012}) \subset \mathcal{V}_b(\anysp{M}{012})$. Then there is a well-defined restriction map to each of the boundary hypersurfaces of $\anysp{M}{012}$. Letting $\mathrm{h}$ denote a boundary hypersurface, we let $N_\mathrm{h} : \mathcal{V}_\gas(\anysp{M}{01}) \to \mathcal{V}_b(\mathrm{h})$ denote the composition of the restriction map with the map that lifts $\mathcal{V}_\gas(\anysp{M}{01})$ to $\mathcal{V}_\gas( \anysp{M}{012} )$.

Each map extends to a homomorphism on differential operators, called a normal homomorphism. We will characterize the image of each restriction map and describe the action of the normal operators.

We shall describe the normal operators at the three boundary hypersurfaces that intersect the diagonal: $\bf{0}, \zff$, and $\tff$. The normal operators at $\bff$ and $\btff$ will also play a role eventually, but their structure is simpler, so we do not spend time on them here.

Recall that $\bf{0}$ is the proper lift of $\{ \epsilon = 0\}$ in $\anysp{M}{012}$. As $P \in \mathrm{Diff}^*_\gas(M)$ can be interpreted as a family of operators $P_\epsilon$ on $M$ parametrized by $\epsilon$, the normal operator $N_{\bf{0}}(P)$ is identified with $P_0$.

Recalling that $\zff$ is identified with $\mathrm{ff}\left( M^2_0 \right) \times [0,1)$, the normal operator $N_{\uff}(P)$ is identified with the family of operators $N_{\mathrm{ff}\left(M^2_0\right)}(P_\epsilon)$ for $\epsilon \in (0,1]$. Therefore, normal operators at $\zff$ can be understood with the standard theory of $0$-pseudodifferential operators.

The normal homomorphism of most interest to us is $N_{\tff}$, because it captures the transitional behaviour of an operator between $\epsilon > 0$ and $\epsilon = 0$.

For concreteness, we fix coordinate systems on $\dbsp{M}{\gas}$. Recall from Figure \ref{fig_double-space-middle-region-coords} that $(\rho, r, u, x_2, y_2)$ and $(r', s, u, x_2, y_2)$ are coordinate charts for $\dbsp{M}{\gas}$. In each chart, $\tff$ is defined by $x_2 = 0$. The first chart induces a chart $(\rho, r, u, y_2)$ for $\tff$. In this chart, $N_{\tff}(x \partial_x) = r \partial_r$ and $N_{\tff}(x \partial_y) = \rho r \partial_u$. The second chart induces a chart $(r', s, u, y_2)$ for $\tff$, in which $N_{\tff}(x \partial_x) = s \partial_s - r' \partial_{r'}$ and $N_{\tff}(x \partial_y) = s \partial_u$.

Define $\mathcal{V}_\gas(\tff) = N_{\tff} \left( \mathcal{V}_\gas(\anysp{M}{01}) \right) = \mathcal{V}_\gas(\anysp{M}{012})|_{\tff}$. Elements of $\mathcal{V}_\gas(\tff)$ act on the fibers of a fibration of $\tff$ over $\partial M$. We explain this as follows. The map $\beta^{012}_2 : \anysp{M}{012} \to M$, projecting onto the second copy of $M$, restricts to a b-fibration $\tff \to \partial M$. In the aforementioned coordinate charts for $\tff$, the map is given by $y_2$. Because $\mathcal{V}_\gas(\anysp{M}{012})$ is tangent to the fibers of $\beta^{012}_{0,2}$ (and hence also to the fibers of $\beta^{012}_2$), it follows that $\mathcal{V}_\gas(\tff)$ is tangent to the fibers of the projection $\tff \to \partial M$. We can then view $N_{\tff}(P)$ as a family of operators parametrized by $\partial M$. For $y^{(0)} \in \partial M$, write $N_{\tff}(P;y^{(0)})$ for the restriction of $N_{\tff}(P)$ to the fiber over $y^{(0)}$.

The fibers can be projectively identified with the space $[0,\infty)^2_b \times \R^{n-1} = \left[ [0,\infty)^2; (0,0) \right] \times \R^{n-1}$, which is a local model for $\sgsp{M}{\gas}$. To see this, note that, in the chart $(\rho, r, u, y_2)$ for $\tff$, the fibers of $y_2$ are parametrized by $(\rho, r, u)$. In the chart $(r', s, u, y_2)$ for $\tff$, the fibers of $y_2$ are parametrized by $(r', s, u)$. As $r = s/\rho$, $r' = \rho/s$, it follows that these two parametrizations define local diffeomorphisms of each fiber $y_2^{-1}(y^{(0)})$ to $[0,\infty)^2_b \times \R^{n-1}$ that lifts the map $(\rho, s, u) : y_2^{-1}(y^{(0)}) \to [0,\infty)^2 \times \R^{n-1}$.

In fact, the fiber-wise action respects a further fibration. As elements of $\mathcal{V}_\gas( \anysp{M}{012} )$ are lifts of vector fields from the first factor of $M$, they are tangent to the level sets of $\epsilon$ and of $x_2$. Therefore $\mathcal{V}_\gas(\tff)$ acts tangent to the level sets of $\epsilon/x_2 = \rho = r's$, and hence also on the joint level sets of $\rho$ and $y_2$. This can also be verified from the explicit coordinate description of $N_{\tff}$. Each joint level set of $\rho$ and $y_2$ is identified with $\R^n_1 = [0,\infty) \times \R^{n-1}$. Thus, given $P \in \mathrm{Diff}^*_\gas(M), y^{(0)} \in \partial M$, we can view $N_{\tff}(P;y^{(0)})$ as a family of differential operators acting on $\R^n_1$.

Note that the $\epsilon/x_2$ 
level set fibration is coordinate-invariant because the coordinate $\epsilon$ is a canonical coordinate for $[0,1)$ and, although the coordinates $(x_2,y_2)$ for $M$ are not canonical, the level sets of $x_2$ are canonical for any fixed value of $y_2$.

\begin{figure}
    \centering
    \includegraphics[width=0.5\linewidth]{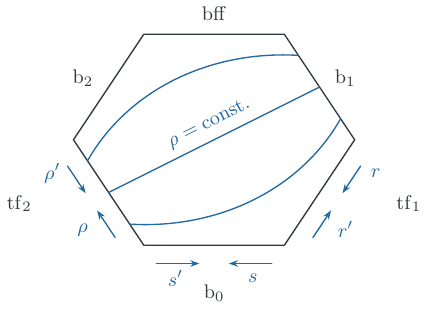}
    \caption{A single fiber of $\tff$ over $\partial M$. A factor of $\R^{n-1}$ parametrized by $u$ (or $v$) is suppressed. Observe that the space itself has a (canonical) fibration given by holding $\rho$ constant.}
    \label{fig_tff-middle-region-fibration}
\end{figure}

Let us compute $N_{\tff}(P; y^{(0)})$ in coordinates. Let $P \in \mathrm{Diff}^m_\gas(M)$ and write $P = \sum_{j + |\gamma| \leq m} a_{j, \gamma}(\epsilon,x,y) (x\partial_x)^j (x\partial_y)^\gamma$. By definition, the functions $a_{j,\gamma}$ lift to be smooth up to the boundary on $\anysp{M}{01} = \sgsp{M}{\gas}$. Away from $\bf{0}$, $(\epsilon, x/\epsilon, y)$ is a smooth chart for this space, so we will write $a_{j,\gamma}(\epsilon, x,y) = \check{a}_{j, \gamma} (\epsilon, x/\epsilon, y)$ for some smooth $\check{a}_{j, \gamma}$. The lift of $P$ to $\anysp{M}{012}$, written in the coordinate chart $(\rho, s, u, x_2, y_2)$, is as follows:
\begin{align*}
    P = \sum_{j+|\gamma|\leq m} \check{a}_{j, \gamma}(x_2 \rho, s/\rho, x_2 u + y_2 ) (s\partial_{s})^j (s \partial_{u})^\gamma.
\end{align*}
Then $N_{\tff}(P;y^{(0)})$ is given by fixing $y_2 = y^{(0)} \in \partial M$ and taking $x_2 \to 0$:
\begin{align}
    N_{\tff}(P;y^{(0)}) = \sum_{j+|\gamma|\leq m} \check{a}_{j, \gamma}(0, s/\rho, y^{(0)} ) (s\partial_{s})^j (s \partial_{u})^\gamma. \label{eqn_normal-operator-tff}
\end{align}
This operator acts on the fiber of $\tff$ over $y^{(0)}$, parametrized by $(\rho, s, u)$. Alternatively, it can be viewed as a family of operators, indexed by $\rho$, acting on a space parametrized by $(s, u)$.

It can be helpful to distinguish the model space $[0,\infty)^2_b \times \R^{n-1}$ from the fibers of $\mathrm{tff}$ themselves. Use tilded coordinates $(\tilde{\epsilon}, \tilde{x}, \tilde{y})$ for the blown-down space $[0,\infty)^2 \times \R^{n-1}$. The mapping from $\tff$ to $[0,\infty)^2_b \times \R^{n-1}$ is given by $\tilde{x} = s$, $\tilde{\epsilon} = \rho$ and $\tilde{y} = u$. Thus, we can write $N_{\tff}(P)$ acting on the model space as follows:
\begin{align*}
    N_{\tff}(P;y^{(0)}) = \sum_{j+|\gamma|\leq m} \check{a}_{j, \gamma}(0,\tilde{x}/\tilde{\epsilon},y^{(0)}) (\tilde{x} \partial_{\tilde{x}})^j (\tilde{x}  \partial_{\tilde{y}})^\gamma.
\end{align*}

The half-density valued kernel of this operator is as follows
\begin{align*}
    K_{N_{\tff}(P;y^{(0)})} &= \mathcal{K}_{N_{\tff}(P;y^{(0)})} \frac{|d\tilde{\epsilon} d \tilde{x}_1 d\tilde{x}_2 d\tilde{y}_1 d\tilde{y}_2|^{1/2}}{\tilde{x}_2^{n}}
\end{align*}
where
\begin{align}
    \mathcal{K}_{N_{\tff}(P;y^{(0)})} &=
    \tilde{x}_2^n \sum_{j+|\gamma|\leq m} \check{a}_{j \alpha}(0,\tilde{x}_1/\tilde{\epsilon},y_2) (\tilde{x}_1 \partial_{\tilde{x}_1})^j (\tilde{x}_1  \partial_{\tilde{y}_1})^\gamma \delta(\tilde{x}_1-\tilde{x}_2) \delta(\tilde{y}_1 - \tilde{y}_2) \nonumber \\
    &= \sum_{j+|\gamma|\leq m} \check{a}_{j \alpha}(0,\tilde{x}_1/\tilde{\epsilon},y_2) (\tilde{x}_1 \partial_{\tilde{x}_1})^j (\tilde{x}  \partial_{\tilde{y}_1})^\gamma \delta\left( \frac{\tilde{x}_1}{\tilde{x}_2} -1 \right) \delta\left( \frac{\tilde{y}_1 - \tilde{y}_2}{\tilde{x}_2} \right) \label{eqn_normal-op-kernel-of-diff-op}
\end{align}

\subsection{Normal operators of pseudodifferential operators}

We can alternatively construct $N_{\tff}(P)$ by restricting the Schwartz kernel of $P$ to $\tff$. The upshot is that this extends the definition of the normal homomorphism to pseudodifferential operators.

Let $P$ be as above. Recall from \eqref{eqn_kernel-gas-diff-op} its half-density-valued Schwartz kernel is
\begin{align*}
    K_P = \sum_{j+|\gamma|\leq m} \check{a}_{j, \gamma}(\epsilon,x_1/\epsilon,y_1) (x_1\partial_{x_1})^j (x_1\partial_{y_1})^\gamma \delta(x_1-x_2) \delta(y_1-y_2) |d\epsilon dx_1 dx_2 dy_1 dy_2|^{1/2}
\end{align*}
Now fix the coordinate system $(\rho, s, u, x_2, y_2)$ on $\anysp{M}{012}$. Changing coordinates, we find that
\begin{align}
   K_P = \mathcal{K}_P(\rho,s,u,x_2,y_2) \frac{|d\rho ds du dx_2 dy_2|^{1/2}}{x_2^{(n-1)/2}} \label{eqn_normal-operator-tff-normalization}
\end{align}
where
\begin{align*}
    \mathcal{K}_P(\rho,s,u,x_2,y_2) = \sum_{j+|\gamma|\leq m} \check{a}_{j, \gamma}(x_2 \rho, s/\rho, x_2 u + y_2 ) (s\partial_{s})^j (s \partial_{u})^\gamma \delta(s-1) \delta(u)
\end{align*}
The density $x_2^{-(n-1)/2} |d\rho ds du dx_2 dy_2|^{1/2}$ is locally equivalent to $\sigma$ as defined before, but they differ in their scaling at infinity.

Fix a point $y_2 = y^{(0)}$ and restrict to $\tff$ by taking $x_2 \to 0$; we find that
\begin{align*}
    \mathcal{K}_P|_{\tff}(\rho,s,u,y^{(0)}) = \sum_{j+|\gamma|\leq m} \check{a}_{j, \gamma}(0,s/\rho,y^{(0)}) (s\partial_{s})^j (s \partial_{u})^\gamma \delta(s-1) \delta(u)
\end{align*}

Observe that this is just $N_{\tff}(P;y^{(0)}) \delta(s-1) \delta(u)$. Also observe, by comparing with \eqref{eqn_normal-op-kernel-of-diff-op}, that
\begin{align*}
    \mathcal{K}_{N_{\tff}(P;y^{(0)})}(\tilde{\epsilon}, \tilde{x}_1, \tilde{x}_2, \tilde{y}_1, \tilde{y}_2) = \mathcal{K}_P|_{\tff}\left( \frac{\tilde{\epsilon}}{\tilde{x}_2}, \frac{\tilde{x}_1}{\tilde{x}_2}, \frac{\tilde{y}_1 - \tilde{y}_2}{\tilde{x}_2}, y^{(0)} \right)
\end{align*}
Using this correspondence between $\mathcal{K}_P|_{\tff}$ and $N_{\tff}(P)$, we can extend $N_{\tff}$ to an arbitrary $P \in \Psi^*_\gas(M)$, and further to an arbitrary element of the large calculus (by normalising by an appropriate power of a boundary-defining function so the leading order at $\tff$ is $0$).

\subsection{Transformed normal operator}

Analogously to the $0$-calculus \cite{mazzeo-melrose_meromorphic-extension-resolvent_1987} \cite{mazzeo_differential-edge-operators-I_1991} \cite{hintz_0-calculus-mazzeo-melrose_2023}, we shall define two additional model operators associated to a \ztrtext{} operator: the transformed normal operator and the Bessel operator. The point is to exploit the symmetries of \ztrtext{} operators.

Let $P \in \mathrm{Diff}^m_\gas(M)$ and write $P = \sum_{j + |\gamma| \leq m} \check{a}_{j, \gamma}(\epsilon,x/\epsilon,y) (x\partial_x)^j (x\partial_y)^\gamma$ as before. The transformed normal operator is given by conjugating $N_{\tff}(P;y^{(0)})$ with the Fourier transform exchanging $u$ for $\eta$.
\begin{align}
    \hat{N}_{\tff}(P;y^{(0)},\eta) = \sum_{j, \gamma} \check{a}_{j, \gamma}(0,s/\rho,y^{(0)}) (s\partial_{s})^j (is \eta)^\gamma. \label{eqn_transformed-normal-operator-tff}
\end{align}
Similarly to the normal operator, this operator can be thought of as acting on the model space $[0,\infty)^2_b$, where the factors of $[0,\infty)$ are parametrized by $\tilde{\epsilon}, \tilde{x}$. In that case, we write $\hat{N}_{\tff}(P;y^{(0)},\eta) = \sum_{j, \gamma} \check{a}_{j, \gamma}\left(0,\tilde{x}/\tilde{\epsilon},y^{(0)}\right) (\tilde{x}\partial_{\tilde{x}})^j (i\tilde{x} \eta)^\gamma$. The kernel of this operator is a distribution on a blow-up of $[0,\infty)^3$, parametrized by $\tilde{\epsilon}, \tilde{x}_1, \tilde{x}_2$. As with the normal operator, we can recover the kernel of $\hat{N}_{\tff}(P;y^{(0)},\eta)$ from the kernel of $P$. Normalize this kernel as $K_{\hat{N}_{\tff}(P;y^{(0)},\eta)} = \mathcal{K}_{\hat{N}_{\tff}(P;y^{(0)},\eta)} x_2^{-1} |d\tilde{\epsilon} d\tilde{x}_1 d\tilde{x}_2|^{1/2}$. Write $\hat{\mathcal{K}}_P|_{\tff}(\rho, s, \eta, y^{(0)})$ for the Fourier transform (exchanging $u$ for $\eta$) of $\mathcal{K}_P|_{\tff}(\rho, s, u, y^{(0)})$. Then,
\begin{align*}
    \mathcal{K}_{\hat{N}_{\tff}(P;y^{(0)},\eta)}(\tilde{\epsilon}, \tilde{x}_1, \tilde{x}_2) = \hat{\mathcal{K}}_P|_{\mathrm{tff}}\left( \frac{\tilde{\epsilon}}{\tilde{x}_2}, \frac{\tilde{x}_1}{\tilde{x}_2},\tilde{x}_2\eta,y_2 \right)
\end{align*}
This formula works equally well to define the transformed normal operator of a general \ztrtext{} pseudodifferential operator.

\subsection{Bessel operator}

The Bessel operator of a differential operator $P$ is obtained from the transformed normal operator by the change of model variables $t = \tilde{x} |\eta|$ and $t_0 = \tilde{\epsilon}|\eta|$. As with the transformed normal operator, it acts on the model space $[0,\infty)^2_b$, where the factors are parametrized by $t_0, t$. However, it can be more appropriate to view it as an operator acting on the space $\left[ [0,\infty]^2 ; (0,0); (0,\infty); (\infty,0); (\infty,\infty) \right]$; i.e., a compactified square with its four corners blown up.
\begin{align}
    B_{\tff}(P;y^{(0)},\hat{\eta})(P) = \sum_{j, \gamma} \check{a}_{j \alpha}(0,t/t_0,y^{(0)}) (t\partial_{t})^j (it \hat{\eta})^\gamma. \label{eqn_Bessel-operator-tff}
\end{align}

Its kernel lives on a blow-up of $[0,\infty]^3$, parametrized by $t_0, t_1, t_2$. We will have more to say about the geometry of this space later. We choose to normalize the half-density-valued kernel of $B_{\tff}(P)$ with respect to a b-half-density, as $K_{B_{\tff}(P;y^{(0)},\hat{\eta})(P)} = \mathcal{K}_{B_{\tff}(P;y^{(0)},\hat{\eta})(P)} (t_0 t_1 t_2)^{-1/2} |dt_0 dt_1 dt_2|^{1/2}$. This is a different kind of normalization to the one used for the normal operator and transformed normal operator. The reason is that it is scale-invariant in $(t_0, t_1, t_2)$, which will be important because we will need to study the asymptotics of Bessel operators at infinity. For now, we remark that the kernel of $B_{\tff}(P)$ can be recovered from the kernel of $P$ by
\begin{align}
    \mathcal{K}_{B_{\tff}(P;y^{(0)},\hat{\eta})}(t_0, t_1, t_2) = (t_0 t_1)^{1/2} t_2^{-1/2} \hat{\mathcal{K}}_P|_{\mathrm{tff}}\left(\frac{t_1}{t_2},\frac{t_0}{t_2},t_2\hat{\eta},y^{(0)}\right) \label{eqn_fourier-transformed-kernel-to-bessel-kernel}
\end{align}
The inverse of this relation is
\begin{align}
    |\eta|^{-1/2} \rho^{-1/2} s^{-1/2} \mathcal{K}_{B_{\tff}(P;y^{(0)},\hat{\eta})}(\rho |\eta|, s|\eta|, |\eta|) =  \hat{\mathcal{K}}_P|_{\mathrm{tff}}\left(\rho,s,\eta,y^{(0)}\right). \label{eqn_bessel-kernel-to-fourier-transformed-kernel}
\end{align}
  \section{The gas giant Laplacian}

We now begin to describe a joint parametrix construction for the family of Laplacians associated to the truncations $M_\epsilon$ of a gas giant $(M, g)$.

Let $M$ be an $n$-dimensional manifold for $n \geq 3$ with a global boundary-defining function $x$. Let $g$ be an $\alpha$-gas giant metric on $M$ for $\alpha \in (0,2)$. That is, $g = \frac{\overline{g}}{x^\alpha}$, where $\overline{g}$ is a smooth Riemannian metric on $M$ up to $\partial M$. There is a normal form for $g$ in local coordinates:

\begin{lemma}
\label{lemma_gas-giant-metric-normal-form}
\cite{dHIKM_gas-giants_2024}

    Let $y$ be a local coordinate system on $\partial M$, so $(x,y)$ is a local coordinate system on $M$ near the boundary. There exists a smooth family of metrics $h_x = h_x(y,dy)$ on $\partial M$ such that
    \begin{align}
        g = \frac{dx^2 + h_x}{x^\alpha} \label{eqn_gas-giant-metric-normal-form}
    \end{align}
    The metric $h_0$ is uniquely determined by $g$.
\end{lemma}

\begin{lemma}
\label{lemma_gas-giant-laplacian}
\cite{dHIKM_gas-giants_2024}

    With $g$ as in \eqref{lemma_gas-giant-metric-normal-form}, the Laplacian of $g$ is expressed as
    \begin{align}
        \Delta_g = x^\alpha \left( \partial_x^2 + q(x,y) \partial_x + \Delta_{h(x)} \right) - x^{\alpha-1} \alpha\left( \frac{n}{2} - 1 \right) \partial_x \label{eqn_gas-giant-laplacian}
    \end{align}
    for a smooth function $q$ depending on $h$.
\end{lemma}

%--Recently promoted from draft
This operator has indicial roots $\gamma_- = 0, \gamma_+ = \alpha(n/2-1) + 1$. In particular, this means that any $f \in L^2(dV_g)$ satisfying $\Delta_g f = 0$ has a polyhomogeneous expansion in $x$ with index set generated by $(\gamma_-, 0), (\gamma_+, 0)$ \cite{dHIKM_gas-giants_2024}, provided that $\gamma_+ - \gamma_-$ is not an integer. (If $\gamma_+ - \gamma_- \in \Z$, then the polyhomogeneous expansion picks up log terms.) For the sake of simplicity, henceforth we assume $\gamma_+ - \gamma_- \notin \Z$.

\begin{proposition}
\label{proposition_gas-giant-laplacian-essentially-self-adjoint}
\cite{dHIKM_gas-giants_2024}

Suppose $\alpha > 2/n$. The operator $-\Delta_g : x^{2-\alpha} H^2_0(M, dV_g) \subset L^2(M, dV_g) \to L^2(M, dV_g)$ is essentially self-adjoint, Fredholm, and has discrete spectrum consisting of eigenvalues $\{ 0 \leq \lambda_0 < \lambda_1 \leq \ldots \}$.
    
\end{proposition}

For the remainder of this paper, we shall impose the restriction that $\alpha > 2/n$. In the non-self-adjoint case, all of our results have analogues for the Dirichlet Laplacian $(\Delta_g, \mathcal{D}_{\mathrm{Dir}})$, where $\mathcal{D}_{\mathrm{Dir}} = \mathcal{D}_\mathrm{max} \cap x^{1-\alpha/2} H^2_0(dV_g)$. Here $\mathcal{D}_{\mathrm{max}}$ is the maximal domain, and the weighted space being chosen to admit $\mathcal{O}(x^{\gamma_+})$ functions and exclude $\mathcal{O}(x^{\gamma_-})$ functions.

The authors of \cite{dHIKM_gas-giants_2024} studied truncations of gas giants. Let $g$ be as in \eqref{lemma_gas-giant-metric-normal-form}. Let $\chi$ be a smooth, nonnegative, nonincreasing function on $[0,\infty)$ satisfying $\chi =1$ on $[0,1]$, $\chi = 0$ on $[4,\infty)$, and $|\chi'| \leq 1/2$ everywhere. Let $M_\epsilon$ be the subset of $M$ where $x \geq \epsilon$. (Rescale $x$ if necessary so that this set is nonempty for $\epsilon \in [0,1)$.) Let $x_\epsilon$ denote the restriction of $x$ to $M_\epsilon$. The metric $g$ restricts to a metric $g_\epsilon$ on $M_\epsilon$. It is a nonsingular metric because $x_\epsilon$ is constrained by $x_\epsilon \geq \epsilon > 0$.
\begin{align*}
    g_\epsilon &= \frac{dx_\epsilon^2 + h_{x_\epsilon}}{x_\epsilon^\alpha}.
\end{align*}
Now let $F(x,y) = (x + \epsilon \chi(x/\epsilon),y)$. This defines a map $M \to M_\epsilon$. Let $g_\epsilon$ be the pullback of $g$ along this map. Then, $g_\epsilon$ pulls back to the following metric (which we also denote by $g_\epsilon$):
\begin{align*}
    g_\epsilon &= \frac{J^2 dx^2 + h_{x_\epsilon}}{(x + \epsilon \chi(x/\epsilon))^\alpha}
\end{align*}
where $J = \frac{\partial x_\epsilon}{\partial x} = 1 + \chi'(x/\epsilon)$. Assume, for simplicity, that $h_x$ is independent of $x$. Allowing $h_x$ to depend on $x$ only produces lower-order terms that make no difference to any of our arguments. Then, the Laplacian of $g_\epsilon$ is
\begin{align*}
     \Delta_{g_\epsilon} &= x_\epsilon^{\alpha n/2} J^{-1} \partial_x \left( x_\epsilon^{-\alpha (n/2-1)} J^{-1} \partial_x \right) + x_\epsilon^\alpha \Delta_h \\
    &= x_\epsilon^\alpha J^{-2} \partial_x^2
    - \alpha(n/2-1) x_\epsilon^{\alpha-1} J^{-1} \partial_x
    - x_\epsilon^\alpha (\partial_x J) J^{-3} \partial_x
    + x_\epsilon^\alpha \Delta_h
\end{align*}
For our analysis, it is preferable to define a family of rescaled operators $P_\epsilon$ that lies in the \ztrtext{} calculus.
\begin{align*}
    P_\epsilon &= \frac{x^2}{x_\epsilon^\alpha} \Delta_{g_\epsilon} \\
    &= x_\epsilon^{\alpha (n/2-1)} J^{-1} x^2 \partial_x \left( x_\epsilon^{-\alpha (n/2-1)} J^{-1} \partial_x \right) + x^2 \Delta_h \\
    &= J^{-2} x^2 \partial_x^2
    - \alpha(n/2-1) \frac{x}{x_\epsilon} J^{-1} x \partial_x
    - (x\partial_x J) J^{-3} x\partial_x
    + x^2 \Delta_h
\end{align*}
The ratio $x/\epsilon$ is smooth on $\sgsp{M}{\gas}$ up to the boundary hypersurfaces $\bf{1}$ and $\tf[1]$ but blows up at the boundary hypersurface $\bf{0}$ (the lift of $\{ \epsilon = 0 \}$). As $\chi'$ is supported on $[1,4]$ and $|\chi'| \leq 1/2$, $J$ lifts to be smooth on $\sgsp{M}{\gas}$ up to the entire boundary, as well as bounded above and below. The ratio $\frac{x}{x_\epsilon} = \frac{1}{1 + (\epsilon/x) \chi(x/\epsilon)}$ also lifts to be smooth on $\sgsp{M}{\gas}$ up to the boundary. Write $P$ to denote $P_\epsilon$ viewed as a differential operator on $\sgsp{M}{\gas}$. The above reasoning shows that $P \in \mathrm{Diff}^2_\gas(M)$.

The rescaled operator $P_0 = x^{2-\alpha} \Delta_g$ maps $x^{2-\alpha} H^2_0(dV_g) \to x^{2-\alpha} L^2(dV_g)$. Then it is natural to view $P_\epsilon$ as mapping $x_\epsilon^{1-\alpha} x H^2_0(dV_{g_\epsilon}) \to x_\epsilon^{1-\alpha} x L^2(dV_{g_\epsilon})$. This reduces to the aforementioned domain when $\epsilon = 0$. When $\epsilon > 0$, $x_\epsilon^{1-\alpha} x L^2(dV_{g_\epsilon}) = x L^2(dV_{\overline{g}})$ is a Dirichlet domain as it excludes functions which are $\mathcal{O}(1)$ up to $\partial M$.

We shall construct parametrices for $P_\epsilon$, i.e. left and right inverses modulo a compact error. Our construction is `uniform' in that the parametrices lie in the \ztrtext{} calculus. The construction follows a standard recipe in the theory of degenerate operators on manifolds with boundary \cite{mazzeo-melrose_meromorphic-extension-resolvent_1987} \cite{mazzeo_differential-edge-operators-I_1991} \cite{melrose_aps-index-theorem_1993}. We construct a first parametrix by inverting the symbol of $P$, modulo a residual error. The error term is noncompact. To improve it to a compact error, we solve model problems for $P_\epsilon$ at the front faces of $\dbsp{M}{\gas}$.

We shall prove the following:
\begin{theorem}
\label{theorem_laplacian-parametrix-with-compact-remainder}

There exists an operator $Q \in \Psi^{-2,*}_{\gas}(M)$, where $Q_\epsilon : x_\epsilon^{1-\alpha} x L^2(dV_{g_\epsilon}) \to x_\epsilon^{1-\alpha} x H^2_0(dV_{g_\epsilon})$ is uniformly bounded for $\epsilon \in [0,1)$ and an operator $R \in \Psi^{-\infty,*}_{\gas}(M)$, where $R_\epsilon: x_\epsilon^{1-\alpha} x L^2(dV_{g_\epsilon}) \to x_\epsilon^{1-\alpha} x L^2(dV_{g_\epsilon})$ is compact and uniformly bounded, such that $PQ = I - R$.

Similarly, there exist $Q' \in \Psi^{-2,*}_{\gas}(M)$, uniformly bounded as a family of operators $x_\epsilon^{1-\alpha} x L^2(dV_{g_\epsilon}) \to x_\epsilon^{1-\alpha} x H^2_0(dV_{g_\epsilon})$, and $R' \in \Psi^{-\infty,*}_{\gas}(M)$, compact and uniformly bounded $x_\epsilon^{1-\alpha} x H^2_0(dV_{g_\epsilon}) \to x_\epsilon^{1-\alpha} x H^2_0(dV_{g_\epsilon})$ such that $Q' P = I - R'$.
\end{theorem}
%--

\section{First parametrix}

Recall from Lemma \ref{lemma_symbol-space} that symbols of \ztrtext{} differential operators lie in $S^{\{m\}}(M) \cong S^{\{m\}} \left( \left( \beta^{01}_1 \right)^*\left({}^0T^*M\right) \right)$. Sections of ${}^0T^*M$ can be expressed as $\xi \frac{dx}{x} + \zeta \cdot \frac{dy}{x}$. Then use $(x,y,\xi,\zeta)$ as  coordinates on ${}^0T^*M$. The operator $P$ has full symbol
\begin{align*}
    p := \sigma_\gas(P) &= - J^{-2} \xi^2 - \left( \beta \frac{x}{x_\epsilon} J^{-1} + J^{-2} - (x \partial_x J) J^{-3} \right) i \xi - |\zeta|_h^2
\end{align*}
where $\beta = \alpha(n/2-1)$. Our previous assumption that the two indicial roots of $P$ do not differ by an integer is equivalent to the assumption that $\beta \notin \Z$. The leading-order part is $-\left( J^{-2} \xi^2 + |\zeta|^2_h \right)$, where $J$ is bounded above and below on $\sgsp{M}{\gas}$. Therefore, $p$ is elliptic and it follows, by quantizing, that there exists $Q_0$ such that $P Q_0 = I - R_0$, where $R_0$ is a residual operator in the small \ztrtext{} calculus.

We must improve $R_0$ to a compact error. Per Proposition \ref{proposition_gas-giant-operator-compactness}, the index family of the remainder must satisfy $E_{\tff}, E_{\zff} > 0, E_{\btff} > -1$, and $E_{\bff} > -1$. However, because $R_0$ lies in the small calculus, it has index set $\emptyset$ at $\btff, \bff$, so we need only concern ourselves with $\tff$ and $\uff$ for the time being.

To produce an error term that is compact, we need to modify our parametrix so that the error is vanishing at $\tff$ and $\uff$. We do this by inverting the normal operators of $P$ at each of these faces.

\section{Transition front face}

\subsection{The model problem}

In Section \ref{section_normal-operators}, we described the transition normal operator $N_{\tff}\left(P;y^{(0)}\right)$ of a \ztrtext{} operator, which can be viewed as a family of operators, parametrized by $y^{(0)} \in \partial M$, acting on the fibers of a fibration of $\tff$ over $\partial M$. The fibers can be locally identified with the model space $[0,\infty)^2_b \times \R^{n-1}$. We also described the transformed normal operator $\hat{N}_{\tff}(P;y^{(0)},\eta)$ and Bessel operator $B_{\tff}(P;y^{(0)},\hat{\eta})$. In this section, we will describe the model problem at $\tff$ and compute these model operators explicitly for the truncated gas giant Laplacian.

Fix a point $y^{(0)} \in \partial M$, a boundary-defining function $x$, and a Riemann normal coordinate system $y$ for $(\partial M, h)$ at $y^{(0)}$ so that $h(y^{(0)},dy) = dy^2$. Use $(x,y)$ as a coordinate system on $M$. This induces a coordinate system $(\rho, s, u, x_2, y_2)$ for $\dbsp{M}{\gas}$, in which $x_2$ is a defining function for $\tff$. This further induces a chart $(\rho, s, u, y_2)$ for $\tff$. In this chart, the fibration of $\tff$ is given by $y_2$. (See Figure \ref{fig_tff-middle-region-fibration}.) Denote a general fiber by $\Sigma$.

By symbolic construction, we obtained $Q_0 \in \Psi^{-2}_\gas(M)$ such that $P Q_0 = I - R_0$ where $R_0 \in \Psi^{-\infty}_\gas(M)$. Passing to transition normal operators:
\begin{align*}
    N_{\tff}\left( P; y^{(0)} \right) N_{\tff}\left( Q_0; y^{(0)} \right) &= N_{\tff}(I) - N_{\tff}\left( R_0; y^{(0)} \right)
\end{align*}

We aim to solve for a correction term $Q_1$ such that
\begin{align}
    N_{\tff}\left( P; y^{(0)} \right) N_{\tff}\left(Q_1; y^{(0)}\right) &= N_{\tff}\left( R_0; y^{(0)} \right) \label{eqn_model-problem-tff-operator-form}
\end{align}
Recall that, for an operator $A$ with kernel $K_A = \mathcal{K}_A \sigma$, $N_{\tff}\left(A; y^{(0)}\right)$ corresponds to the restriction of $\mathcal{K}_A$ to the fiber of $\tff$ over $y^{(0)}$. In the chart $(\rho,s,u,x_2,y_2)$, write $\mathcal{K}_{A} = \mathcal{K}_{A}(\rho,s,u,x_2,y_2)$ and write its restriction as $n\!\left(A;y^{(0)}\right)(\rho,s,u):= \mathcal{K}_{A}|_{\tff}\left(\rho,s,u,y^{(0)}\right) = \mathcal{K}_{A}\left(\rho,s,u,0,y^{(0)}\right)$. Also let $\hat{n}\!\left(A;y^{(0)}\right)(\rho,s,\eta)$ be the Fourier transform (exchanging $u$ for $\eta$) of $n\!\left(A;y^{(0)}\right)(\rho,s,u)$. Then \eqref{eqn_model-problem-tff-operator-form} is equivalent to
\begin{align}
    N_{\tff}\left( P; y^{(0)} \right) n\!\left( Q_1; y^{(0)} \right) = n\!\left( R_0; y^{(0)} \right) \label{eqn_model-problem-tff-function-form}
\end{align}

As $R_0$ belongs to the small calculus, $\mathcal{K}_{R_0}$ is rapidly vanishing at all boundary hypersurfaces of $\dbsp{M}{\gas}$ except for those intersecting the diagonal: $\zff$, $\tff$, and $\bf{0}$. Thus $\mathcal{K}_{R_0}|_{\tff}$ is vanishing at all boundary hypersurfaces of $\tff$ except for $\zff \cap \tff$ and $\bf{0} \cap \tff$. These two boundary hypersurfaces are disjoint (see Figure \ref{fig_tff-corner-structure}): the boundary hypersurface $\zff \cap \tff$ lies in the upper region and $\bf{0} \cap \tff$ lies in the middle and lower regions. 

We shall give a global existence argument for the solution to \eqref{eqn_model-problem-tff-function-form}, then prove that it is polyhomogeneous on $\tff$. To do that, we shall study the model problem separately in the upper region and the middle/lower regions.

\subsection{The normal operator}

First we shall describe the normal operator $N_{\tff}(P;y^{(0)})$. Select the chart $(\rho, s, u, x_2, y_2)$ for $\dbsp{M}{\gas}$ and observe that $x\partial_x = s\partial_s$ and $x\partial_y = s \partial_u$. Re-expressing $P$ in this chart and then taking $x_2 \to 0$ (passing to $\tff$) yields $N_{\tff}(P;y^{(0)})$, expressed in the chart $(\rho, s, u)$ for $\Sigma$.
\begin{align}
    N_{\tff}(P;y^{(0)}) &= (s + \rho\chi)^\beta J^{-1} s^2 \partial_s \left( (s+\rho\chi)^{-\beta} J^{-1}\partial_s \right) + s^2 \Delta_u \nonumber \\
    &= J^{-2} s^2 \partial_s^2 - \beta \frac{s}{s+\rho \chi} J^{-1} s\partial_s - (s\partial_s J) J^{-3} s\partial_s + s^2 \Delta_u \label{eqn_normal-operator-tff-gg-laplacian-general-form}
\end{align}

Simplify this operator by a change of coordinates. Let $(\check{\rho}, \check{s}, \check{u})$ be another copy of the same chart on $\Sigma$. Define a map by $(\check{\rho}, \check{s}, \check{u}) = (\rho, s + \rho(\chi(s/\rho)-1), u)$. This map extends smoothly to the boundaries of the chart. Indeed, this map is a diffeomorphism (in the sense of manifolds with corners) onto its image; in particular, it maps each boundary hypersurface of $\Sigma$ to itself.

Observe that $\frac{\partial \check{s}}{\partial s} = 1 + \chi'(s/\rho) = J(s/\rho)$. Then, changing coordinates in \eqref{eqn_normal-operator-tff-gg-laplacian-general-form},
\begin{align}
    N_{\tff}(P; y^{(0)}) &= s^2 (\check{s} + \check{\rho})^\beta \partial_{\check{s}} \left( (\check{s} + \check{\rho})^{-\beta} \partial_{\check{s}} \right) + s^2 \Delta_u \nonumber \\
    &= \left(\frac{s}{\check{s}} \right)^2 \left( (\check{s})^2 \partial_{\check{s}}^2 - \beta \left( \frac{\check{s}}{\check{s}+\check{\rho}} \right) \check{s} \partial_{\check{s}} + (\check{s})^2 \Delta_u \right). \label{eqn_normal-operator-tff-gg-laplacian-general-form-checked-chart}
\end{align}
The quantity $(s/\check{s})^2$ is smooth and bounded above and below on $\Sigma$. Therefore, nothing is lost by declaring that
\begin{align*}
    N_{\tff}(P; y^{(0)}) &= s^2 \partial_s^2 - \beta \left( \frac{s}{s+\rho} \right) s \partial_s + s^2 \Delta_u.
\end{align*}

where $\chi = \chi(s/\rho)$, $J = 1 + \chi'(s/\rho)$, and $\Delta_u = \sum_i \partial_{u_i}^2$ is the flat Laplacian. 

\subsection{Invertibility of the normal operator}

From now on, we will leave the dependence on the boundary point $y^{(0)}$ implicit. View $N_{\tff}(P)$ as acting on $[0,\infty) \times \R^{n-1}$, parametrized by $(s,u)$. The operator $N_{\tff}(P)$ is symmetric for the volume form $d\tilde{V} = (s+\rho)^{-\beta} s^{-2} \,ds \,du$ on $\R^n_1 = [0,\infty) \times \R^{n-1}$. Define the following quadratic form, for $f, g \in C^\infty_0({\R^n_1})$:
\begin{align*}
    \langle f, g \rangle &= -\int_{\R^n_1} \left( N_{\tff}(P)f \right) g \,d\tilde{V} \\
    &= \int_{\R^n_1} (s \partial_{s} f) (s\partial_{s} g) \,d\tilde{V} + \int_{\R^n_1} (s \nabla_u f) \cdot (s \nabla_u g) \,d\tilde{V}
\end{align*}
The quantity $\left( \|f\|^2_{L^2(d\tilde{V})} + \langle f, f\rangle \right)^{1/2}$ is equivalent to the norm on the $0$-Sobolev space $H^1_0( d\tilde{V})$. In fact, the pairing is coercive for this norm. It suffices to show that $\int f^2 \,d\tilde{V} \lesssim \int (s \partial_{s} f)^2 \,d\tilde{V}$ for $f \in C^\infty_0(\R^n_1)$. This is a Hardy-type inequality. Start by integrating by parts:
\begin{align*}
	\int_{\R^n_1} f(s, u)^2 \,d\tilde{V} &= \int_{\R^{n-1}} \int_0^\infty \frac{f(s, u)^2}{s^2 (s+1)^{\beta}} \,ds du \\
	&= \int_{\R^{n-1}} \int_0^\infty \partial_{s} \left( \frac{f(s, u)^2}{(s+1)^\beta} \right) \frac{ds du}{s}
\end{align*}
Expand the derivative with the product rule, discarding the term where $\partial_{s}$ hits on $1/(s+1)^\beta$ as it is nonpositive:
\begin{align*}
	\int_{\R^n_1} f(s, u)^2 \,d\tilde{V} &\lesssim \int_{\R^{n-1}} \int_0^\infty \frac{f(s,u) \partial_{s} f(s, u)}{s (s+1)^{\beta}} \,ds du \\
    &= \int_{\R^{n-1}} \int_0^\infty f(s,u) s \partial_{s} f(s, u) \,d\tilde{V}
\end{align*}
Applying the Cauchy-Schwarz inequality with respect to $d\tilde{V},$
\begin{align*}
    \int_{\R^n_1} f(s, u)^2 \,d\tilde{V} &\lesssim \left( \int \int_0^\infty f(s, u)^2 \,d\tilde{V} \right)^{1/2} \left( \int \int_0^\infty \left(s \partial_{s} f(s, u) \right)^2 \, d\tilde{V} \right)^{1/2}
\end{align*}
Dividing by $\left( \int f^2 \,d\tilde{V} \right)^{1/2}$ yields the inequality. Thus $N_{\tff}(P) : H^1_0(\R^n_1, d\tilde{V}) \to H^{-1}(\R^n_1, d\tilde{V})$ is invertible.

In the model problem \eqref{eqn_model-problem-tff-function-form}, view $n(R_0; y^{(0)}) = n(R_0; y^{(0)})(\rho,s,u)$ as a function on $\Sigma$. Recall that $n(R_0; y^{(0)})$ is smooth up to $\bf{0}, \tff$, and $\zff$ and rapidly vanishing at all other boundary hypersurfaces. For fixed $\rho \in (0,\infty)$, $s = 0$ corresponds to $\bf{1}$, $s = \infty$ corresponds to $\bf{2}$, and $|u| = \infty$ corresponds to $\btff$. Therefore, for fixed $\rho$, $n(R_0; y^{(0)})(\rho, s, u)$ is rapidly vanishing as $s \to 0$ and rapidly decaying as $|(s,u)| \to \infty$, and therefore lies in $L^2(\R^n_1, d\tilde{V}) \subset H^{-1}(\R^n_1, d\tilde{V})$. Then, by invertibility of $N_{\tff}(P)$, \eqref{eqn_model-problem-tff-function-form} has a unique solution $n(Q_1)$ satisfying the property that $n(Q_1)(\rho,s,u) \in H^1_0(\R^n_1, d\tilde{V})$ for each fixed $\rho$.

We shall now establish the polyhomogeneity of the solution on $\Sigma$, using separate arguments in the upper, middle, and lower regions.

\subsection{The model problem in the upper region}

In the upper region, express $N_{\tff}(P)$ in the chart $(\rho', s, u)$ for $\Sigma$, where $\rho' = 1/\rho$.
\begin{align}
    N_{\tff}(P;y^{(0)}) = s^2 \partial_s^2 - \beta \frac{s\rho'}{s\rho' + 1}  s\partial_s + s^2 \Delta_u \label{eqn_normal-operator-tff-gg-laplacian-general-form-s-rhop-chart}
\end{align}
Recall that the upper region of $\tff$ can be locally identified with $\dbsp{\tilde{M}}{0'}$, where $\tilde{M} = \R^n_1$, parametrized by $(\tilde{x}, \tilde{y})$. Let $\bf{i}[][\widetilde]$ denote the lift of $\{\tilde{x}_i = 0\}$ from $\tilde{M}^2$ to $\dbsp{\tilde{M}}{0'}$, $\ff[][][\widetilde]$ denote the front face of $\dbsp{\tilde{M}}{0'}$, and $\bff[][][\widetilde]$ denote the `b-front face' \cite{hintz_0-calculus-mazzeo-melrose_2023}.

Utilize the coordinate chart $(\epsilon, s, u, \rho', y_2)$ for the upper region of $\dbsp{M}{\gas}$, which induces the chart $(s, u, \rho', y_2)$ for $\tff$. Using the coordinate chart $(\tilde{s}, \tilde{u}, \tilde{x}_2, \tilde{y}_2)$ on $\dbsp{\tilde{M}}{0'}$, the identification is given by the map $\tilde{s} = s, \tilde{u} = u, \tilde{x}_2 = \rho', \tilde{y}_2 = y_2$. By pulling back \eqref{eqn_normal-operator-tff-gg-laplacian-general-form-s-rhop-chart} along this map, and then changing coordinates according to $\tilde{x}_1 = \tilde{x}_2 \tilde{s}$, $\tilde{y}_1 = \tilde{x}_2 \tilde{u} + \tilde{y}_2$, we find that $N_{\tff}(P)$ is identifiable with the lift of the following operator from $\tilde{M}$ to $\dbsp{\tilde{M}}{0'}$ on the left:

\begin{align}
    \tilde{P} = {\tilde{x}}^2 \partial_{\tilde{x}}^2 - \beta \frac{{\tilde{x}}}{{\tilde{x}} + 1} {\tilde{x}}\partial_{\tilde{x}} + {\tilde{x}}^2 \Delta_{\tilde{y}} \label{eqn_normal-operator-tff-gg-laplacian-upper-modified}
\end{align}
This is an elliptic $0$-operator on $\tilde{M}$ with indicial roots $0$ and $1$.

Similarly, identify $n(R_0)$ with a kernel $\mathcal{K}_{\tilde{R}_0}$ on $\dbsp{\tilde{M}}{0'}$. In the upper region, $n(R_0)$ is rapidly vanishing at $\bf{1}$ and $\bf{2}$ and smooth up to $\zff$, so $\mathcal{K}_{\tilde{R}_0}$ is rapidly vanishing at $\tilde{\mathrm{b}}_1$ and $\tilde{\mathrm{b}}_2$ and smooth up to $\ff[][][\widetilde]$, hence an element of the small $0$-calculus.

The model equation \eqref{eqn_model-problem-tff-function-form} becomes the following:
\begin{align}
    \tilde{P} \tilde{Q}_1 &= \tilde{R}_0 \label{eqn_model-problem-tff-upper-modified}
\end{align}

By adapting the argument for the invertibility of $N_{\tff}(P)$, the operator $\tilde{P} : H^1_0( d\tilde{V}_2) \to H^{-1}(d\tilde{V}_2)$ is invertible, where $d \tilde{V}_2 = (\tilde{x}+1)^{-\beta} \tilde{x}^{-2} \,d\tilde{x} d\tilde{y}$. Equivalently, the operator $\tilde{L} = (\tilde{x}+1)^{-\beta/2} \tilde{P} (\tilde{x}+1)^{\beta/2}$ is invertible from $\tilde{x}H^1_0(d\tilde{x}d\tilde{y}) \to \tilde{x} H^{-1}(d\tilde{x}d\tilde{y})$. By the parametrix construction for the $0$-calculus \cite{mazzeo_differential-edge-operators-I_1991} \cite{hintz_0-calculus-mazzeo-melrose_2023}, the inverse $\tilde{L}^{-1}$ lies in the large $0$-calculus on $\tilde{M}$. Its index sets at $\bf{1}[][\widetilde], \bf{2}[][\widetilde], \ff[][][\widetilde], \bff[][][\widetilde]$ are (with respect to the normalization chosen in \cite{mazzeo_differential-edge-operators-I_1991}) generated by $(1,0)$, $(-1,0)$, $(0,0)$, and $(1/2,0)$, respectively. Because the factor $(\tilde{x} +1)$ is smooth and nonvanishing, the same is true of $\tilde{P}^{-1} = (\tilde{x}+1)^{\beta/2} \tilde{L}^{-1} (\tilde{x}+1)^{-\beta/2}$. Using the composition formula for the $0$-calculus, it follows that the solution $\tilde{Q} = \tilde{P}^{-1} \tilde{R}_0$ to \eqref{eqn_model-problem-tff-upper-modified} has the same index family as $\tilde{P}^{-1}$.

\begin{proposition}
    The solution $n(Q_1)$ to \eqref{eqn_model-problem-tff-function-form} satisfies polyhomogeneous estimates on $\Sigma$ up to $\bf{1}, \bf{2}$, $\zff$, and $\bff$, with index sets generated by $(1,0)$, $(n-1,0)$, $(0,0)$, and $(n,0)$, respectively.
\end{proposition}
\begin{proof}
    Undo the local identification of $\tff$ with $\dbsp{\tilde{M}}{0'}$. The change in index sets between $\tilde{Q}_1$ and $n(Q_1)$ is owed to the difference between the choice of normalization in \cite{mazzeo_differential-edge-operators-I_1991} and the choice we made in \eqref{eqn_normal-operator-tff-normalization}.
\end{proof}

\subsection{The model problem in the middle and lower regions}

The remainder of this section will be dedicated to proving that the solution to the model problem \eqref{eqn_model-problem-tff-function-form} is polyhomogeneous in the middle and lower regions of $\tff$. We shall study the problem by reducing it to the level of the Bessel operator.

Modify \eqref{eqn_model-problem-tff-function-form} by taking a Fourier transform exchanging $u$ for $\eta \in \R^{n-1}$:
\begin{align}
    \hat{N}_{\tff}(P;\eta) \hat{n}(Q_1)(\rho,s,\eta) = \hat{n}(R_0)(\rho,s,\eta) \label{eqn_transformed-model-problem-tff-lower}
\end{align}
Here $\hat{N}_{\tff}(P;\eta)$ is the transformed normal operator of $P$, defined as in \eqref{eqn_transformed-normal-operator-tff}:
\begin{align*}
    \hat{N}_{\tff}(P;\eta) &= s^2 \partial_s^2 - \beta \left( \frac{s}{s+\rho} \right) s \partial_s - s^2 |\eta|^2.
\end{align*}
Now, pass to Bessel operators using the relation \eqref{eqn_fourier-transformed-kernel-to-bessel-kernel}:
\begin{align}
    B_{\tff}(P) B_{\tff}(Q_1;\hat{\eta}) &= B_{\tff}(R_0;\hat{\eta}) \label{eqn_Bessel-model-problem-tff-lower}
\end{align}
where $B_{\tff}(P)$ is defined as in \eqref{eqn_Bessel-operator-tff} by
\begin{align*}
    B_{\tff}(P) &= t^2 \partial_{t}^2 - \beta \frac{t}{t+t_0} t\partial_{t} - t^2
\end{align*}
Note that, in general, the Bessel operator of an arbitrary element of $\mathrm{Diff}^*_\gas(M)$ depends on $\hat{\eta} = \eta/|\eta|$. However, in this particular case, there is no such dependence.

Solve for $B_{\tff}(Q_1;\hat{\eta})$ by setting
\begin{align}
    B_{\tff}(Q_1;\hat{\eta}) &= B_{\tff}(G) B_{\tff}(R_0;\hat{\eta}) \label{eqn_Bessel-model-problem-solution}
\end{align}
where $B_{\tff}(G)$ is a Green's function of $B_{\tff}(P)$, satisfying $B_{\tff}(P) B_{\tff}(G) = \delta(t_1-t_2) |dt_0 dt_1 dt_2|^{1/2}$.

\subsection{The Green's function for $B_{\tff}(P)$}

We solve now for $B_{\tff}(G)$, the Dirichlet Green's function for $B_{\tff}(P)$. On $[0,\infty)_t$, for a fixed $t_0$, $B_{\tff}(P)$ has solutions $u_0$, $u_1$ which vanish at $t = 0$ and decay exponentially as $t \to \infty$, respectively.
\begin{align*}
    u_0(t;t_0,\beta) &= (t+t_0)^{\frac{\beta+1}{2}} \left( I_{\frac{\beta+1}{2}}(t+t_0) - \frac{I_{\frac{\beta+1}{2}}(t_0)}{K_{\frac{\beta+1}{2}}(t_0)} K_{\frac{\beta+1}{2}}(t+t_0) \right) \\
    u_1(t;t_0,\beta) &= (t+t_0)^{\frac{\beta+1}{2}} K_{\frac{\beta+1}{2}}(t+t_0)
\end{align*}
where $I_\nu, K_\nu$ are modified Bessel functions with parameter $\nu$, given by (for $\nu$ not an integer)
\begin{align*}
    I_\nu(x) &= \frac{1}{2^\nu} x^\nu\sum_{m=0}^\infty \frac{1}{m! \Gamma(m+\nu+1) 2^{2m}} x^{2m} \\
    K_\nu(x) &= \frac{\pi}{2\sin \nu \pi} \left( I_{-\nu}(x) - I_\nu(x) \right) \\
    &= \frac{2^{\nu-1}\pi}{\sin \nu \pi} x^{-\nu} \sum_{m=0}^\infty \frac{1}{m! \Gamma(m-\nu+1) 2^{2m}} x^{2m} \\
    &\quad - \frac{2^{-\nu-1}\pi}{\sin \nu \pi} x^{\nu} \sum_{m=0}^\infty \frac{1}{m! \Gamma(m+\nu+1) 2^{2m}} x^{2m}
\end{align*}
The Wronskian factor $W = u_0 u_1' - u_0' u_1$ is given by $W(t;t_0,\beta) = -(t+t_0)^\beta$. It follows then that $B_{\tff}(G)$ has the following half-density-valued kernel:
\begin{align}
    K_{B_{\tff}(G)}(t_0,t_1,t_2) &= -\frac{1}{t_2^2 (t_2+t_0)^\beta} \Big( u_0(t_1;t_0,\beta) u_1(t_2,t_0,\beta) H(t_2-t_1) + \nonumber \\
    \quad &+ u_0(t_2;t_0,\beta) u_1(t_1;t_0,\beta) H(t_1-t_2) \Big) |dt_0 dt_1 dt_2|^{1/2} \label{eqn_Bessel-Greens-function}
\end{align}
where $H$ is the Heaviside step function.

The next step is to construct an appropriate pseudodifferential calculus in which the Green's function and the residual operator $B_{\tff}(R_0;\hat{\eta})$ lie. We shall study composition of operators in this calculus, and use it to determine the properties of the solution \eqref{eqn_Bessel-model-problem-solution} to the model problem \eqref{eqn_Bessel-model-problem-tff-lower}.

\section{The double b-scattering-transition calculus}

The Green's function lies in a calculus that is essentially a concatenation of four calculi: two copies each of a b-transition calculus and a scattering-transition calculus. Here we define each of these two calculi.

\subsection{The b-transition calculus}

The b-transition calculus is a slight modification of the b-surgery calculus, as described in \cite{mcdonald_laplacian-cone-singularities_1990} \cite{mazzeo-melrose_analytic-surgery_1995}. McDonald considered the Laplacian of a one-parameter family of metrics that develops a conic singularity at an interior hypersurface as the parameter goes to zero. In our situation, the relevant singularity occurs at the boundary, so additional blow-ups are required to resolve the intersection of the diagonal with the boundary.

Let $\tilde{T} = [0,\infty)$. We will construct the double space $\tilde{T}^2_{\btr}$ as a blow-up of $[0,\infty) \times \tilde{T}^2$, where the first factor of $[0,\infty)$ is not denoted by $\tilde{T}$, so as to emphasize that it corresponds to a transition parameter rather than a spatial coordinate. Parametrize $[0,\infty) \times \tilde{T}^2$ by $\left(\tilde{t}_0, \tilde{t}_1, \tilde{t}_2\right)$, where $\tilde{t}_1, \tilde{t}_2$ are viewed as spatial variables. Let $\bf{i}[][\widetilde] = \left\{ \tilde{t}_i = 0 \right\} \subset [0,\infty) \times \tilde{T}^2$. Let $\anydiag{12}[][\widetilde] = \{ \tilde{t}_1 = \tilde{t}_2 \}$ be the diagonal. As we have done previously, we will consistently use the same notation for a submanifold as its proper lift to a blow-up.

Define the following submanifolds of $[0,\infty) \times \tilde{T}^2$.
\begin{align}
    B_{12} &:= \bf{1}[][\widetilde] \cap \bf{2}[][\widetilde] \nonumber \\
    C_{12} &:= \bf{0}[][\widetilde] \cap \bf{1}[][\widetilde] \cap \bf{2}[][\widetilde] \nonumber \\
    C_i &:= \bf{0}[][\widetilde] \cap \bf{1}[][\widetilde] \text{ for $i = 1,2$.}
    \label{eqn_b-tr-double-space-submanifolds-to-blow-up}
\end{align}
Define $\dbsp{\tilde{T}}{\btr} = \big[ [0,\infty) \times \tilde{T}^2; C_{12}; C_1; C_2; B_{12} \big]$. This is identical to the double space of \cite{mcdonald_laplacian-cone-singularities_1990} \cite{mazzeo-melrose_analytic-surgery_1995} except it involves an additional blow-up, namely that of $B_{12}$. This double space was also used in \cite{guillarmou-hassell_low-energy-resolvent-asymptotically-conic_2007}. Denote the front face of the blow-up of $B_{12}$ by $\bff[][][\widetilde]$, that of $C_{12}$ by $\btff[][][\widetilde]$, and that of the blow-up of $C_i$ by $\tf[i][][\widetilde]$.

\begin{figure}
    \centering
    \includegraphics[width=0.5\linewidth]{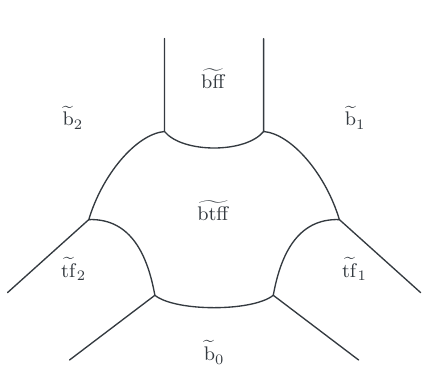}
    \caption{The b-transition double space $\dbsp{\tilde{T}}{\btr}$.}
    \label{fig_b-transition-double-space}
\end{figure}

The small calculus $\Psi^*_{\btr}(\tilde{T})$ consists of distributional half-density sections $K_A = \mathcal{K}_A (\tilde{t}_0 \tilde{t}_1 \tilde{t}_2)^{-1/2} \left|d\tilde{t}_0 \,d\tilde{t}_1 \, d\tilde{t}_2\right|^{1/2}$ on $\dbsp{\tilde{T}}{\btr}$, where $\mathcal{K}_A$ is conormal at $\anydiag{12}[][\widetilde]$, smooth up to the boundary hypersurfaces $\bf{0}[][\widetilde], \btff[][][\widetilde]$, and $\bff[][][\widetilde]$, rapidly vanishing at all other boundary hypersurfaces, and rapidly decaying as $|(\tilde{t}_1, \tilde{t}_2)| \to \infty$. The large calculus $\Psi^{*,*}_{\btr}(\tilde{T})$ consists of sections where $\mathcal{K}_A$ is conormal at $\anydiag{12}[][\widetilde]$, polyhomogeneous up to all boundary hypersurfaces, and rapidly decaying as $|(\tilde{t}_1, \tilde{t}_2)| \to \infty$. As with the 0-transition calculus, the calculus is filtered by spaces $\Psi^{\ell,\mathcal{E}}_{\btr}(\tilde{T})$, where $\mathcal{E}$ is the index family of $\mathcal{K}_A$ and $\ell$ is the order of its conormal singularity. (See \cite{mcdonald_laplacian-cone-singularities_1990} \cite{mazzeo-melrose_analytic-surgery_1995} for more details.)

This calculus obeys the following composition law:

\begin{proposition}
\label{proposition_b-transition-calculus-composition-formula}

    Suppose $A \in \Psi^{\ell, \mathcal{E}}_\btr(\tilde{T})$ and $B \in \Psi^{m, \mathcal{F}}_\btr(\tilde{T})$. Suppose moreover that $E_{\bf{2}} + F_{\bf{1}} > 0$. Then the composition $AB$ is well-defined and $AB \in \Psi^{\ell+m, \mathcal{G}}_\btr(\tilde{T})$, where $\mathcal{G}$ is given by
    \begin{align*}
         G_{\bf{0}}&=(E_{\bf{0}}+F_{\bf{0}}-1/2)
         \mathbin{\overline\cup}(E_{\tf[2]}+F_{\tf[1]}-1/2)\\
         G_{\bf{1}}&=E_{\bf{1}}
         \mathbin{\overline\cup}(E_{\bff}+F_{\bf{1}})\\
         G_{\bf{2}}&=F_{\bf{2}}
         \mathbin{\overline\cup}(E_{\bf{2}}+F_{\bff})\\
         G_{\btff}&=(E_{\btff}+F_{\btff}-1/2)
         \mathbin{\overline\cup}(E_{\tf[1]}+F_{\tf[2]}-1/2)\\
         G_{\tf[1]}&=(E_{\tf[1]}+F_{\bf{0}}-1/2)
         \mathbin{\overline\cup}(E_{\btff}+F_{\tf[1]}-1/2)\\
         G_{\tf[2]}&=(E_{\tf[2]}+F_{\btff}-1/2)
         \mathbin{\overline\cup}(E_{\bf{0}}+F_{\tf[2]}-1/2)\\
         G_{\bff}&=(E_{\bf{1}}+F_{\bf{2}})
         \mathbin{\overline\cup}(E_{\bff}+F_{\bff}).
    \end{align*}
\end{proposition}
\begin{proof}
    The proof is very similar to the proof of the composition formula in \cite{mcdonald_laplacian-cone-singularities_1990}, with only the minor added complication of the face $\bff[][][\widetilde]$. One constructs a triple space $\trsp{\tilde{T}}{\btr}$ as a blow-up of $[0,\infty) \times \tilde{T}^3$ and shows that it admits a triplet of b-fibrations $\trsp{\tilde{T}}{\btr} \to \dbsp{\tilde{T}}{\btr}$ that extend the projections $[0,\infty) \times \tilde{T}^3 \to [0,\infty) \times \tilde{T}^2$ which project out the first, second, and third copies of $\tilde{T}$. Then, the pushforward theorem yields the composition formula.
\end{proof}

\subsection{The scattering-transition calculus}

The scattering-transition calculus is a slight modification of the cusp-surgery calculus described in \cite{albin-rochon-sher_fibered-cusps_2021} \cite{anghel_dirac-spectrum-degenerating-riemann-surfaces_2024}. As with the b-transition calculus, the modification is the presence of additional scattering-type blow-ups to resolve the intersection of the diagonal with the corner.

The scattering-transition double space $\dbsp{\tilde{T}}{\sctr}$ is constructed from $\dbsp{\tilde{T}}{\btr}$. We define the following submanifolds of $\dbsp{\tilde{T}}{\btr}$.
\begin{align*}
    E_{12} &:= \anydiag{12}[][\widetilde] \cap \bff[12][][\widetilde] \\
    F_{12} &:= \anydiag{12}[][\widetilde] \cap \btff[12][][\widetilde]
\end{align*}
Define $\dbsp{\tilde{T}}{\sctr} = \left[ \dbsp{\tilde{T}}{\btr}; F_{12}; E_{12} \right]$. This differs from the double space of \cite{mcdonald_laplacian-cone-singularities_1990} \cite{mazzeo-melrose_analytic-surgery_1995} by the blow-up of $E_{12}$ (and the blow-up of $B_{12}$ involved in $\dbsp{\tilde{T}}{\btr}$). A very similar space appears in \cite{guillarmou-hassell_low-energy-resolvent-asymptotically-conic_2007}, but $\dbsp{\tilde{T}}{\sctr}$ differs from that by the blow-up of $F_{12}$. Denote the front face of the blow-up of $E_{12}$ by $\scff[][][\widetilde]$ and that of the blow-up of $F_{12}$ by $\sctff[][][\widetilde]$.

\begin{figure}
    \centering
    \includegraphics[width=0.5\linewidth]{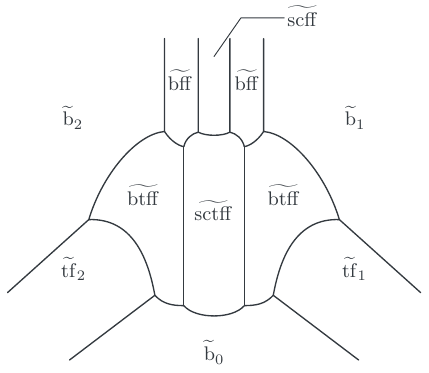}
    \caption{The scattering-transition double space $\dbsp{\tilde{T}}{\sctr}$.}
    \label{fig_sc-transition-double-space}
\end{figure}

The small scattering-transition calculus $\Psi^*_{\sctr}(\tilde{T})$ consists of those sections $K_A = \mathcal{K}_A (\tilde{t}_0 \tilde{t}_1 \tilde{t}_2)^{-1/2} \left|d\tilde{t}_0 \,d\tilde{t}_1 \, d\tilde{t}_2\right|^{1/2}$ on $\dbsp{\tilde{T}}{\sctr}$ where $\mathcal{K}_A$ is conormal at the diagonal, smooth up to $\bf{0}[][\widetilde], \scff[][][\widetilde]$, and $ \sctff[][][\widetilde]$, and rapidly vanishing elsewhere. The large calculus $\Psi^*_{\sctr}(\tilde{T})$ is defined identically to $\Psi^*_{\btr}(\tilde{T})$ except that the kernels live on $\dbsp{\tilde{T}}{\sctr}$. (See \cite{albin-rochon-sher_fibered-cusps_2021} \cite{anghel_dirac-spectrum-degenerating-riemann-surfaces_2024} for more details.)

\begin{proposition}
\label{proposition_sc-transition-calculus-composition-formula}

    Suppose $A \in \Psi^{\ell, \mathcal{E}}_\sctr(\tilde{T})$ and $B \in \Psi^{m, \mathcal{F}}_\sctr(\tilde{T})$. Suppose moreover that $E_{\bf{2}} + F_{\bf{1}} > 0$. Then the composition $AB$ is well-defined and $AB \in \Psi^{\ell+m, \mathcal{G}}_\sctr(\tilde{T})$, where $\mathcal{G}$ is given by
    \begin{align*}
         G_{\bf{0}}&=(E_{\bf{0}}+F_{\bf{0}}-1/2)
         \mathbin{\overline\cup}(E_{\tf[2]}+F_{\tf[1]}-1/2)\\
         G_{\bf{1}}&=E_{\bf{1}}
         \mathbin{\overline\cup}(E_{\bff}+F_{\bf{1}})
         \mathbin{\overline\cup}(E_{\scff}+F_{\bf{1}}+1/2)\\
         G_{\bf{2}}&=F_{\bf{2}}
         \mathbin{\overline\cup}(E_{\bf{2}}+F_{\bff})
         \mathbin{\overline\cup}(E_{\bf{2}}+F_{\scff}+1/2)\\
         G_{\btff}&=(E_{\btff}+F_{\btff}-1/2)
         \mathbin{\overline\cup}(E_{\sctff}+F_{\btff})
         \mathbin{\overline\cup}(E_{\tf[1]}+F_{\tf[2]}-1/2)
         \mathbin{\overline\cup}(E_{\btff}+F_{\sctff})\\
         G_{\tf[1]}&=(E_{\tf[1]}+F_{\bf{0}}-1/2)
         \mathbin{\overline\cup}(E_{\btff}+F_{\tf[1]}-1/2)
         \mathbin{\overline\cup}(E_{\sctff}+F_{\tf[1]})\\
         G_{\tf[2]}&=(E_{\tf[2]}+F_{\btff}-1/2)
         \mathbin{\overline\cup}(E_{\bf{0}}+F_{\tf[2]}-1/2)
         \mathbin{\overline\cup}(E_{\tf[2]}+F_{\sctff})\\
         G_{\bff}&=(E_{\bf{1}}+F_{\bf{2}})
         \mathbin{\overline\cup}(E_{\bff}+F_{\bff})
         \mathbin{\overline\cup}(E_{\scff}+F_{\bff}+1/2)
         \mathbin{\overline\cup}(E_{\bff}+F_{\scff}+1/2)\\
         G_{\sctff}&=(E_{\btff}+F_{\btff})
         \mathbin{\overline\cup}(E_{\tf[1]}+F_{\tf[2]})
         \mathbin{\overline\cup}(E_{\sctff}+F_{\sctff})\\
         G_{\scff}&=(E_{\bf{1}}+F_{\bf{2}}+1/2)
         \mathbin{\overline\cup}(E_{\bff}+F_{\bff}+1/2)
         \mathbin{\overline\cup}(E_{\scff}+F_{\scff}+1/2)
    \end{align*}
\end{proposition}
\begin{proof}
    The proof is very similar to the proof of the composition formula in \cite{albin-rochon-sher_fibered-cusps_2021}, with only the minor added complication of the faces $\bff[][][\widetilde]$ and $\scff[][][\widetilde]$. As with the composition formula for the scattering-transition calculus, one constructs a triple space $\trsp{\tilde{T}}{\sctr}$, shows that it admits a triplet of b-fibrations to the $\dbsp{\tilde{T}}{\sctr}$, and invokes the pushforward theorem.
    % \note{Essentially identical to the \cite{albin-rochon-sher_fibered-cusps_2021} proof, only difference is you have to check that a few extra blow-ups commute.}
\end{proof}

\subsection{The double b-scattering-transition calculus}

Now we shall build the double b-scattering transition calculus from the calculi constructed above. The calculus lives on a blow-up of $[0,\infty]^3$, parametrized by $t_0, t_1, t_2$. Near $t_0 = 0, (t_1, t_2) = 0$ and near $t_0 = 0, (t_1, t_2) = \infty$, the space resembles the b-transition double space. Near $t_0 = \infty, (t_1, t_2) = 0$ and near $t_0 = \infty, (t_1, t_2) = \infty$, the space resembles the scattering-transition double space.

(Note that a calculus called a `b-scattering transition calculus' appears in \cite{guillarmou-hassell_low-energy-resolvent-asymptotically-conic_2007} \cite{kottke_callias-type_2015}. This is purely a coincidence; that calculus is adapted to a different kind of transition behavior.)

\begin{figure}
    \centering
    \includegraphics[width=0.5\linewidth]{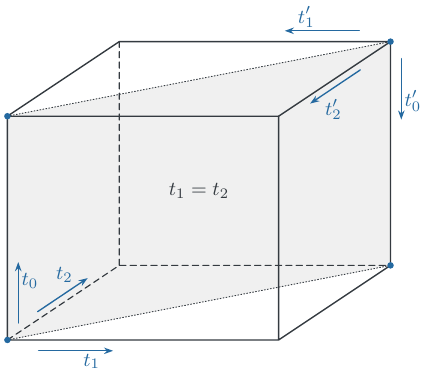}
    \caption{The space $[0,\infty] \times T^2$, with the four key coordinate charts indicated.}
    \label{fig_Bessel-triple-space-blowndown}
\end{figure}

Let $T = [0,\infty]$ be the compactification of the half-open interval $[0,\infty)$. If the interval is parametrized by $t$, then $t' := 1/t$ is a boundary-defining function of $\infty$. Let $[0,\infty] \times T^2$ be parametrized by $t_0, t_1, t_2$, where the first factor of $[0,\infty]$ is not denoted by $T$ to emphasize that it is non-spatial. This space is a cube with six boundary hypersurfaces, $\bf{i}[0][^{T}] = \{ t_i = 0 \}$ and $\bf{i}[\infty][^{T}] = \{ t_i' = 0 \}$ for $i = 0,1,2$. Let $\anydiag{12}[][^{T}] = \{ t_1 = t_2 \}$ denote the diagonal. (See Figure \ref{fig_Bessel-triple-space-blowndown}.)

Consider four coordinate charts whose origins lie at the four vertices of $[0,\infty] \times T^2$ touching the diagonal: $(t_0, t_1, t_2)$, $(t_0', t_1, t_2)$,$(t_0, t_1', t_2')$, and $(t_0', t_1', t_2')$. We would like to perform blow-ups in each of these charts. However, to ensure that the blow-ups commute, we must first perform blow-ups at off-diagonal edges and vertices. Define
\begin{align*}
    L_{12}^{jk} &= \bf{1}[j][^{T}] \cap \bf{2}[k][^{T}] \\
    L_{012}^{ijk} &= \bf{0}[i][^{T}] \cap \bf{1}[j][^{T}] \cap \bf{2}[k][^{T}]
\end{align*}
where $i,j,k \in \{0,\infty\}$.
Then construct the space
\begin{align*}
    \mathring{T}^2 = \left[[0,\infty] \times T^2 ; L_{012}^{00\infty}; L_{012}^{0\infty0}; L_{012}^{\infty0\infty}; L_{012}^{\infty\infty0}; L_{12}^{0\infty}; L_{12}^{\infty0} \right].
\end{align*}
Collectively, refer to the front faces of the above blow-ups as `far-from-diagonal faces'. Now we can perform the b-transition and scattering-transition calculi blow-ups independently on $\mathring{T}^2$.

\begin{figure}
    \centering
    \includegraphics[width=0.5\linewidth]{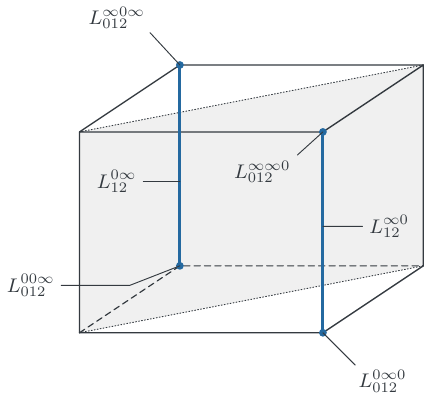}
    \caption{The space $[0,\infty] \times T^2$, with the far-from-diagonal submanifolds that are blown up to produce $\mathring{T}^2$ highlighted.}
    \label{fig_Bessel-triple-space-far-from-diagonal-submanifolds}
\end{figure}

Identify each of the aforementioned charts with the model coordinates $\left(\tilde{t}_0, \tilde{t}_1, \tilde{t}_2\right)$. Then, the subset $\left\{ \tilde{t}_0, \tilde{t}_1, \tilde{t}_2 < \infty \right\} \subset \mathring{T}^2$ can be identified with the space $[0,\infty) \times \tilde{T}^2 = [0,\infty)^3$. Construct $T^2_{\dbsctr}$ by performing the blow-ups that turn $[0,\infty) \times \tilde{T}^2$ into $\dbsp{\tilde{T}}{\btr}$ in the $(t_0, t_1, t_2)$ and $(t_0', t_1, t_2)$ charts, and performing the blow-ups that turn $[0,\infty) \times \tilde{T}^2$ into $\dbsp{\tilde{T}}{\sctr}$ in the $(t_0, t_1', t_2')$ and $(t_0', t_1', t_2')$ charts. Thanks to the far-from-diagonal blow-ups, the double b-transition construction commutes with the double scattering-transition construction. The two b-transition constructions are not entirely independent from one another as there is one submanifold ($B_{12}$ in the notation of \eqref{eqn_b-tr-double-space-submanifolds-to-blow-up}) whose blow-up is required by both constructions. It is only necessary to blow this submanifold up once. The rest of the two b-transition constructions commute with each other. Similarly, there are two submanifolds involved in both of the scattering-transition constructions, namely $B_{12}$ and $E_{12}$. Again, it suffices to blow these up once; the rest of the two constructions commute.

Define $\beta^i_{\btr} : \dbsp{T}{\dbsctr} \to \dbsp{\tilde{T}}{\btr}$ by
\begin{align*}
    \tilde{t}_0 &= \begin{cases}
        t_0 \text{ if i = 0,} \\
        t_0' \text{ if $i = \infty$}
    \end{cases} \\
    \tilde{t}_1 &= t_1 \\
    \tilde{t}_2 &= t_2
\end{align*}
Equivalently, $\beta^0_{\btr}$ (resp. $\beta^\infty_{\btr}$) is a blow-down map. It undoes all blow-ups not involved in the construction of $\dbsp{\tilde{T}}{\btr}$ in the $(t_0, t_1, t_2)$ chart (resp. the $(t_0', t_1, t_2)$ chart), as well as undoing the far-from-diagonal blow-ups.

Analogously define $\beta^i_{\sctr} : \dbsp{T}{\dbsctr} \to \dbsp{\tilde{T}}{\sctr}$ by
\begin{align*}
    \tilde{t}_0 &= \begin{cases}
        t_0 \text{ if i = 0,} \\
        t_0' \text{ if $i = \infty$}
    \end{cases} \\
    \tilde{t}_1 &= t_1' \\
    \tilde{t}_2 &= t_2'
\end{align*}
Equivalently, $\beta^0_{\sctr}$ (resp. $\beta^\infty_{\sctr}$) is a blow-down map undoing all blow-ups not involved in the construction of $\dbsp{\tilde{T}}{\sctr}$ in the $(t_0, t_1', t_2')$ chart (resp. the $(t_0', t_1', t_2')$ chart), as well as undoing the far-from-diagonal blow-ups.

The targets of the maps $\beta^i_{\btr}$ and $\beta^i_{\sctr}$ are the `constituent double spaces' of $\dbsp{T}{\dbsctr}$. We refer to a constituent space with $i = 0$ (respectively $i = \infty$) as a `lower' (respectively `upper') space. Each of the boundary hypersurfaces of $\dbsp{T}{\dbsctr}$, excluding the far-from-diagonal faces, is the inverse image of a boundary hypersurface of some constituent space. Some faces belong to multiple constituent spaces. For example, the two constituent b-transition double spaces share a b-front face, i.e. $\big(\beta^0_{\btr}\big)^{-1}\big( \bff[][][\widetilde] \big) = \big(\beta^\infty_{\btr}\big)^{-1}\big( \bff[][][\widetilde] \big)$. As another example, the two lower constituent spaces share a copy of $\bf{0}[][\widetilde]$, i.e. $\big(\beta^0_{\btr}\big)^{-1}\big( \bf{0}[][\widetilde] \big) = \big(\beta^0_{\sctr}\big)^{-1}\big( \bf{0}[][\widetilde] \big)$.

\begin{definition}

    The small double b-scattering transition calculus $\Psi^*_{\dbsctr} \left( T\right)$ consists of half-density-valued distributions $K_Q = \mathcal{K}_Q (t_0 t_1 t_2)^{-1/2} |dt_0 \,dt_1 \, dt_2|^{1/2}$ on $T^2_{\dbsctr}$, where $\mathcal{K}_Q$ is conormal with respect to $\anydiag{12}[][^{T}]$, smooth down to all boundary hypersurfaces intersecting $\anydiag{12}[][^{T}]$, and vanishing to infinite order at all other boundary hypersurfaces.
\end{definition}

\begin{definition}
    Suppose $\mathcal{E}$ is an index family on $\dbsp{T}{\dbsctr}$ whose index set at all far-from-diagonal faces is $\emptyset$. An element of $\Psi^{m,\mathcal{E}}_{\dbsctr} \left( T\right)$ is a half-density-valued distribution $K_Q = \mathcal{K}_Q (t_0 t_1 t_2)^{-1/2} |dt_0 \,dt_1 \, dt_2|^{1/2}$, where $\mathcal{K}_Q$ is conormal of order $m$ at $\anydiag{12}[][^{T}]$ and polyhomogeneous on $T^2_{\dbsctr}$ with index family $\mathcal{E}$.
\end{definition}

Owing to the assumption of rapid vanishing at far-from-diagonal faces, operators in $\Psi^{*,*}_{\dbsctr}$ are essentially diagonalized with respect to the constituent calculi. Composition of operators reduces to composition in the constituent calculi. Suppose $\mathcal{E}$ is an index family on $\dbsp{T}{\dbsctr}$ whose index set at all far-from-diagonal faces is $\emptyset$. By identifying faces of $\dbsp{T}{\dbsctr}$ with faces of its constituent spaces using the maps $\beta^i_{\btr}$ and $\beta^i_{\sctr}$, this index family induces two index families $\mathcal{E}_{\btr}^i$ on $\dbsp{\tilde{T}}{\btr}$ and two index families $\mathcal{E}_{\sctr}^i$ on $\dbsp{\tilde{T}}{\sctr}$, for $i = 0, \infty$. To compose two pseudodifferential operators $K_A \in \Psi^{m, \mathcal{E}}_{\dbsctr}(T)$, $K_B \in \Psi^{m', \mathcal{F}}_{\dbsctr}(T)$, it suffices to apply our previous composition formulae (Propositions \ref{proposition_b-transition-calculus-composition-formula} and \ref{proposition_sc-transition-calculus-composition-formula}) to the pairs $(\mathcal{E}_{\btr}^i, \mathcal{F}_{\btr}^i)$ and $(\mathcal{E}_{\sctr}^i, \mathcal{F}_{\sctr}^i)$, for $i = 0,\infty$. This suffices to prove the following composition formula.

\begin{proposition}
\label{proposition_dbsctr-calculus-composition-formula} 
    Suppose $A \in \Psi^{\ell, \mathcal{E}}_{\dbsctr}(T)$, $B \in \Psi^{m, \mathcal{F}}_{\dbsctr}(T)$, where $E_{\bf{2}[0][^{T}]} + F_{\bf{1}[0][^{T}]} > 0$ and $E_{\bf{2}[\infty][^{T}]} + F_{\bf{1}[\infty][^{T}]} > 0$. Then $\mathcal{E}_{\btr}^i$ and $\mathcal{F}_{\btr}^i$ satisfy the conditions for composition to be well-defined in the b-transition calculus; let $\mathcal{G}_{\btr}^i$ be the index family of the composition. The same is true for $\mathcal{E}_{\sctr}^i$ and $\mathcal{F}_{\sctr}^i$ in the scattering-transition calculus; let $\mathcal{G}_{\sctr}^i$ be the index family of that composition.
    
    Define an index family $\mathcal{G} = \{ G(\mathrm{h}) \}$ on $\dbsp{T}{\dbsctr}$ as follows. Let $\mathrm{h}$ be any boundary hypersurface of $\dbsp{T}{\dbsctr}$. If it is a far-from-diagonal face, then let $G(\mathrm{h}) = \emptyset$. Otherwise, let $\beta(\mathrm{h})$ be the set of faces in the constituent double spaces that are identified with $\mathrm{h}$. For each $\mathrm{h}' \in \beta(\mathrm{h})$, let $G(\mathrm{h}')$ be the index set of $\mathrm{h}'$ in whichever of $\{\mathcal{G}_{\btr}^0, \mathcal{G}_{\btr}^\infty, \mathcal{G}_{\sctr}^0, \mathcal{G}_{\sctr}^\infty\}$ corresponds to the constituent space of $\mathrm{h}'$. Then, let $G(\mathrm{h}) = \bigcup_{\mathrm{h}' \in \beta(\mathrm{h})} G(\mathrm{h}')$.

    The operator $AB$ is well-defined and has kernel $K_{AB} \in \Psi^{\ell + m, \mathcal{G}}_{\dbsctr}(T)$.
\end{proposition}

Armed with the calculus $\Psi^{*,*}_{\dbsctr}(T)$, we can return to the model problem at $\tff$.

\section{Solving the lower model problem at $\tff$}

\subsection{The solution lies in $\Psi^{*,*}_{\dbsctr}(T)$}

Recall that we wish to solve the model problem $N_{\tff}(P) n(Q_1) = n(R_0)$ for $n(Q_1)$ and that we reduced this to a model problem for Bessel operators, $B_{\tff}(P) B_{\tff}(Q_1;\hat{\eta}) = B_{\tff}(R_0;\hat{\eta})$, whose solution is $B_{\tff}(Q_1;\hat{\eta}) = B_{\tff}(G) B_{\tff}(R_0;\hat{\eta})$. To understand the solution, we shall show that $B_{\tff}(G)$ and $B_{\tff}(R_0;\hat{\eta})$ lie in $\Psi^{*,*}_{\dbsctr}(T)$.

\begin{lemma}
    The Green's function $B_{\tff}(G)$ lies in $\Psi^{-2,*}_{\dbsctr}(T)$.
\end{lemma}

\begin{proof}
    Recall the expression for $K_{B_{\tff}(G)}$ in $\eqref{eqn_Bessel-Greens-function}$. We must verify that this is polyhomogeneous on $\dbsp{T}{\dbsctr}$ with a conormal singularity at $\diaglift$.

    On $\dbsp{T}{\dbsctr}$, the functions $t_0$, $t_1$, $t_2$, $t_1+t_0$, $t_2+t_0$, and $t_1-t_2$ lift, in any neighborhood of the boundary, to a product of defining functions, or a reciprocal of such a product. It follows that polyhomogeneous functions of $t_1 + t_0$ and $t_2 + t_0$ are polyhomogeneous wherever $t_1 + t_0$ and $t_2 + t_0$ are small.

    This handles the Bessel functions appearing in \eqref{eqn_Bessel-Greens-function} whenever their arguments are small. To handle the situation where their arguments are large, consider their asymptotic forms:
    \begin{align*}
        I_\nu(z) &\sim \frac{1}{\sqrt{2\pi}} z^{-1/2} e^z \sum_{i=0}^\infty (-1)^i\frac{a_{\nu,i}}{z^i} \\
        K_\nu(z) &\sim \sqrt{\frac{\pi}{2}} z^{-1/2} e^{-z} \sum_{i=0}^\infty \frac{a_{\nu,i}}{z^i}
    \end{align*}
    as $z \to \infty$. The asymptotic expansion is a polyhomogeneous expansion in $1/z$, meaning $\frac{I_\nu(z)}{\frac{1}{\sqrt{2\pi}} z^{-1/2} e^z} - \sum_{i=0}^k (-1)^i\frac{a_{\nu,i}}{z^i} = O(z^{-(k+1)})$ and this estimate is stable under conormal derivatives. (The analogous statement is true of $K_\nu$.)
    
    It follows that, if $z = t_0$, $t_1+t_0$, or $t_2+t_0$ is large, Bessel functions of $z$ are equal to an exponential factor multiplied by a polyhomogeneous expansion in $1/z$, which is polyhomogeneous on $\dbsp{T}{\dbsctr}$ because $1/z$ is a product of boundary-defining functions.

    It must be checked that the exponential factors are smooth up to all boundaries. It suffices to check this in the case that $t_1 < t_2$, as the Green's function is, up to a polyhomogeneous prefactor, symmetric under reflection across the diagonal. In the regime $t_1 < t_2$, the Green's function, normalized with respect to the b-half-density $(t_0 t_1 t_2)^{-1/2} |dt_0 dt_1 dt_2|^{1/2}$ is given as follows:
    \begin{align}
        \mathcal{K}_{B_{\tff}(G)}(t_0, t_1, t_2) &= -t_0^{1/2} t_1^{1/2} t_2^{-3/2}(t_1+t_0)^{\frac{\beta+1}{2}} (t_2+t_0)^{\frac{-\beta+1}{2}} K_{\frac{\beta+1}{2}}(t_0)^{-1} K_{\frac{\beta+1}{2}}(t_2+t_0) \nonumber \\
        &\quad \cdot \left( I_{\frac{\beta+1}{2}}(t_1+t_0) K_{\frac{\beta+1}{2}}(t_0) - I_{\frac{\beta+1}{2}}(t_0) K_{\frac{\beta+1}{2}}(t_1+t_0) \right) \label{eqn_Bessel-greens-function-t1<t2}
    \end{align}
    The second term in the parenthetical is more rapidly decaying than the first, as $K_{\frac{\beta+1}{2}}$ is exponentially decaying while $I_{\frac{\beta+1}{2}}$ is growing, so the second term can be ignored for the following analysis.

    In the $t_1 < t_2$ regime, there are these sub-regimes to consider:
    \begin{itemize}
        \item $t_0$ large,
        \item $t_0$ small, $t_2$ small,
        \item $t_0$ small, $t_1$ large,
        \item $t_0$ small, $t_1$ small, $t_2$ large.
    \end{itemize}
    In the first regime, \eqref{eqn_Bessel-greens-function-t1<t2} is asymptotic (up to a polyhomogeneous prefactor) to $e^{t_1 - t_2}$. In the second, the arguments of the Bessel functions are small, and in this case polyhomogeneity has already been established. In the third, \eqref{eqn_Bessel-greens-function-t1<t2} is asymptotic to $e^{t_1 - t_2}$. In the fourth, \eqref{eqn_Bessel-greens-function-t1<t2} is asymptotic to $e^{-t_2}$.

    The functions $e^z$ and $e^{-1/z}$ are smooth down to $z = 0$. As $\dbsp{T}{\dbsctr}$ resolves $t_2$ and $t_1 - t_2$, it follows that the exponential factors involved in the Green's function are smooth.

    The Green's function is conormal at $\anydiag{12}[][^{T}]$ because $\dbsp{T}{\dbsctr}$ resolves $t_1 - t_2$. The conormal order of $-2$ is owed to the Heaviside function. Rapid decay at the required faces can be checked manually. This proves that the $B_{\tff}(G)$ lies in $\Psi^{-2,*}_{\dbsctr}(T)$.
\end{proof}

Using the series expansions of Bessel functions, the index sets of $\mathcal{K}_{B(G)}$ can be reliably computed. We only record the index sets corresponding to the lower constituent b-calculus. These sets are generated by 
\begin{itemize}
        \item $(1/2 + k(\beta+1),0)$ for $k \in \N_0$ at $\left(\beta^0_{\btr}\right)^{-1}\big( \bf{0}[][\widetilde] \big)$
        \item $(3/2,0)$ at $\left(\beta^0_{\btr}\right)^{-1}\big( \bf{1}[][\widetilde] \big)$
        \item $(-1/2,0)$ at $\left(\beta^0_{\btr}\right)^{-1}\big( \bf{2}[][\widetilde] \big)$
        \item $(1/2 + k(\beta+1),0)$ at $\left(\beta^0_{\btr}\right)^{-1}\big( \btff[][][\widetilde] \big)$
        \item $(1+(k+1)(\beta+1),0)$ for $k \in \N_0$ at $\left(\beta^0_{\btr}\right)^{-1}\big( \tf[1][][\widetilde] \big)$
        \item $(k(\beta+1),0)$ for $k \in \N_0$ at $\left(\beta^0_{\btr}\right)^{-1}\big( \tf[2][][\widetilde] \big)$
        \item $(0,0)$ at $\left(\beta^0_{\btr}\right)^{-1}\big( \bff[][][\widetilde] \big)$
\end{itemize}
If $\beta$ is an odd positive integer, then $(\beta+1)/2$ is a positive integer and the series expansions of the Bessel functions take a different form. The consequence is that the above index sets must be modified by increasing some logarithmic exponents.

\begin{lemma}

    For any $\hat{\eta} \in S^{n-2}$, $B_{\tff}(R_0;\hat{\eta}) \in \Psi^{-\infty}_{\dbsctr} \left( T\right)$. It has index sets generated by the following indices:
    \begin{itemize}
        \item $(1/2,0)$ at $\left(\beta^0_{\btr}\right)^{-1}\big( \bf{0}[][\widetilde] \big)$
        \item $(1/2,0)$ at $\left(\beta^0_{\btr}\right)^{-1}\big( \btff[][][\widetilde] \big)$
        \item $(0,0)$ at $\left(\beta^0_{\btr}\right)^{-1}\big( \bff[][][\widetilde] \big)$
        \item $(-1/2,0)$ at $\left(\beta^\infty_{\btr}\right)^{-1}\big( \bf{0}[][\widetilde] \big)$
        \item $(-1/2,0)$ at $\left(\beta^\infty_{\btr}\right)^{-1}\big( \btff[][][\widetilde] \big)$
    \end{itemize}
    Its index sets at all other boundary hypersurfaces are $\emptyset$. All polyhomogeneous estimates are uniform in $\hat{\eta}$.
\end{lemma}
\begin{proof}
    The scalar kernel $b(R_0;\hat{\eta}) := \mathcal{K}_{B_{\tff}(R_0;\hat{\eta})}$ is obtained from $n(R_0)$ using \eqref{eqn_fourier-transformed-kernel-to-bessel-kernel}. Let us determine the asymptotics of $n(R_0)(\rho,s,u)$. Recall that this function is rapidly vanishing at all boundary hypersurfaces of $\tff$ except for $\zff \cap \tff$ and $\bf{0} \cap \tff$, where it is smooth.
   
    Assuming $\rho, s$ are fixed, taking $|u| \to \infty$ corresponds to approaching the boundary hypersurface $\mathrm{btff} \cap \tff$, so $n(R_0)(\rho,s,u)$ is rapidly decaying in $u$. Also note that $\tff \cap \{\rho = 0 \} = \tff \cap (\bf{0} \cup \tf[1])$ and $\tff \cap \{s = 0\} = \tff \cap (\tf[1] \cup \bf{1})$. It follows then that $n(R_0)$ is smooth up to $\rho = 0$ and rapidly vanishing as $s \to 0$. Also note that, assuming $\rho, u$ are fixed, $s \to \infty$ corresponds to approaching the boundary hypersurface $\bf{2} \cap \tff$. Thus $n(R_0)$ is rapidly decaying as $s \to 0$. Finally, note that $\tff \cap \{\rho = \infty\} = \tff \cap \left( \bf{2} \cup \zff \cup \bff \right)$, so $n(R_0)$ is bounded as $\rho \to \infty$. All estimates are stable under b-derivatives in $\rho, s$.

    Due to the rapid decay in $u$, the function $\hat{n}(R_0)(\rho,s,\eta)$ has the same asymptotics as $n(R_0)$, except with $u$ exchanged for $\eta$.

    Therefore, $b(R_0;\hat{\eta})$ is rapidly vanishing at any boundary hypersurface of $\dbsp{T}{\dbsctr}$ where $t_1/t_2 = 0, \infty$ or $t_0/t_2 = \infty$, or $t_2 = \infty$. The index sets at other boundary hypersurfaces are similarly read off from \eqref{eqn_fourier-transformed-kernel-to-bessel-kernel}.
\end{proof}

As $B_{\tff}(G)$ and $B_{\tff}(R_0,\hat{\eta})$ both lie in $\Psi^{*,*}_{\dbsctr} \left( T\right)$, and moreover $B_{\tff}(R_0,\hat{\eta})$ is residual, we can invoke the composition law (Proposition \ref{proposition_dbsctr-calculus-composition-formula}) to conclude that $B_{\tff}(Q_1,\hat{\eta}) \in \Psi^{-\infty,*}_{\dbsctr} \left( T\right)$ and moreover compute the index family of its kernel.

\begin{proposition}
\label{eqn_Bessel-model-problem-solution-is-phg}
    For any $\hat{\eta} \in S^{n-2}$, the operator $B_{\tff}(Q_1,\hat{\eta})$ lies in $\Psi^{-\infty,*}_{\dbsctr} \left( T\right)$ and has index sets generated by the following indices:
    \begin{itemize}
        \item $(1/2+k(\beta+1),0)$ for $k \in \N_0$ at $\left(\beta^0_{\btr}\right)^{-1}\big( \bf{0}[][\widetilde] \big)$
        \item $(3/2,0)$ at $\left(\beta^0_{\btr}\right)^{-1}\big( \bf{1}[][\widetilde] \big)$
        \item $(-1/2,0)$ at $\left(\beta^0_{\btr}\right)^{-1}\big( \bf{2}[][\widetilde] \big)$
        \item $(1/2 + k(\beta+1),0)$ for $k \in \N_0$ at $\left(\beta^0_{\btr}\right)^{-1}\big( \btff[][][\widetilde] \big)$
        \item $(1 + (k+1)(\beta+1),0)$ for $k \in \N_0$ at $\left(\beta^0_{\btr}\right)^{-1}\big( \tf[1][][\widetilde] \big)$
        \item $(k(\beta+1),0)$ for $k \in \N_0$ at $\left(\beta^0_{\btr}\right)^{-1}\big( \tf[2][][\widetilde] \big)$
        \item $(0,0)$ at $\left(\beta^0_{\btr}\right)^{-1}\big( \bff[][][\widetilde] \big)$
    \end{itemize}
    If $\beta$ is an odd positive integer, then these index sets are modified by increasing some logarithmic exponents.
    
    Moreover, the operator has nonempty index sets at $\left(\beta^\infty_{\btr}\right)^{-1}\big( \bf{0}[][\widetilde] \big)$, $\left(\beta^\infty_{\btr}\right)^{-1}\big( \btff[][][\widetilde] \big)$, $\left(\beta^\infty_{\btr}\right)^{-1}\big( \tf[1][][\widetilde] \big)$, and $\left(\beta^\infty_{\btr}\right)^{-1}\big( \tf[2][][\widetilde] \big)$ although we do not record these here. The index sets at all other boundary hypersurfaces are $\emptyset$. In particular, the kernel is rapidly decaying as $(t_1, t_2) \to \infty$.
\end{proposition}
Only the index sets of $B_{\tff}(Q_1,\hat{\eta})$ at the seven faces indicated above will be relevant in what follows. These are the faces of the lower constituent \btrtext{} calculus. Henceforth, we designate $\big(\beta^0_{\btr}\big)^{-1}\big( \bf{i}[][\widetilde] \big)$ by $\bf{i}[][^{T}]$, $\big(\beta^0_{\btr}\big)^{-1}\big( \btff[][][\widetilde] \big)$ by $\btff[][][^{T}]$, and so on.

\subsection{Lifting the solution to the Bessel problem}

The half-density-valued kernel $K_{B_{\tff}(Q_1,\hat{\eta})} = b(Q_1;\hat{\eta})(t_0, t_1, t_2) (t_0 t_1 t_2)^{-1/2} |dt_0 dt_1 dt_2|^{1/2}$ solves the model problem \eqref{eqn_Bessel-model-problem-tff-lower}. Using \eqref{eqn_bessel-kernel-to-fourier-transformed-kernel}, lift this solution to a solution of \eqref{eqn_transformed-model-problem-tff-lower} by setting $\hat{n}(Q_1)(\rho,s,\eta) = |\eta|^{-1/2} \rho^{-1/2} s^{-1/2} b(Q_1;\hat{\eta})\left(\rho|\eta|, s|\eta|, |\eta|\right)$. Further lift this to a solution of \eqref{eqn_model-problem-tff-function-form} by inverting the Fourier transform; set $n(Q_1)(\rho,s,u) = \int e^{i u \cdot \eta} \hat{n}(Q_1)(\rho,s,\eta) \,d\eta$.

We shall use the polyhomogeneity of $b(Q_1;\hat{\eta})$ on $\dbsp{T}{\dbsctr}$ to prove conormal-type and polyhomogeneous-type estimates for $\hat{n}(Q_1)$, and go on to prove estimates for $n(Q_1)$ in the lower and middle regions of $\tff$. More generally, we replace $b(Q_1;\hat{\eta})$ with an arbitrary distribution $b(a)$ on $\dbsp{T}{\dbsctr}$ satisfying conormal or polyhomogeneous estimates.

We introduce suitable coordinate systems. Recall that $(\rho, s, u, y_2)$ is a coordinate system on $\tff$ and that $\tff$ fibers over $\partial M$, where the basepoint is given by $y_2$. Each fiber $\Sigma$ is dual, under the Fourier transform exchanging $u$ for $\eta$, to a space parametrized by $(\rho,s,\eta)$, which we denote by $\hat{\Sigma}$. These coordinates are related to coordinates on $\dbsp{T}{\dbsctr}$ by $t_0 = \rho |\eta|$, $t_1 = s |\eta|$, $t_2 = |\eta|$.

Actually, we must use a family of coordinate systems on $\Sigma$ that collectively cover all corners of the middle region. In Figure \ref{fig_double-space-middle-region-coords}, we provided six such coordinate systems. In fact, three suffice. We define one such system below. We will prove conormal and polyhomogeneous estimates in that coordinate system. The arguments in the remaining charts are similar.

Define:
\begin{align*}
    \nu &= \frac{\epsilon (x_1+x_2)}{x_1 x_2} = \frac{t_0(t_1+t_2)}{t_1 t_2} \\
    q &= t_1 + t_2 \\
    \tau &= \frac{x_1 - x_2}{x_1 + x_2} = \frac{t_1 - t_2}{t_1 + t_2} \\
    U &= \frac{y_1 - y_2}{x_1 + x_2} = \frac{x_2}{x_1 + x_2} u = \frac{1-\tau}{2} u \\
    H &= \frac{x_1 + x_2}{x_2} \eta = \left( \frac{1-\tau}{2} \right)^{-1} \eta
\end{align*}

On $\Sigma$, use the chart $(\nu, \tau, U)$, where $U$ takes values in $\R^{n-1}$. This chart is smooth up to $\Sigma \cap \bf{0}, \Sigma \cap \tf[1]$, and $\Sigma \cap \tf[2]$. On $\hat{\Sigma}$, use the corresponding chart $(\nu, \tau, H)$, where $H$ takes values in $\R^{n-1}$. On $\dbsp{T}{\dbsctr}$, use the coordinate chart $(\nu, \tau, q) \in [0,\infty) \times [-1,1] \times [0,\infty)$. This chart is smooth up to the boundary hypersurfaces $\bf{0}[][^{T}], \tf[1][][^{T}], \tf[2][][^{T}]$, and $\btff[][][^{T}]$.   

In what follows, let $b(a) = b(a)(\nu, \tau, q)$ be a conormal distribution on $\dbsp{T}{\dbsctr}$. Define $\hat{a}(\nu, \tau, H) = \nu^{-1/2} ((1+\tau)/2)^{-1} ((1-\tau)/2)^{n-1} |H|^{-1/2} b(a)(\nu, \tau, |H|)$. Then define $a(\nu, \tau, U) = \int e^{i U \cdot H} \hat{a}(\nu, \tau, H) \,dH$. The choice of normalization of $\hat{a}$ is made so that $a$ is coordinate-invariant. The following proofs are similar to the proof of Proposition 2.20 in \cite{hintz_0-calculus-mazzeo-melrose_2023}.

\begin{figure}
    \centering
    \includegraphics[width=0.5\linewidth]{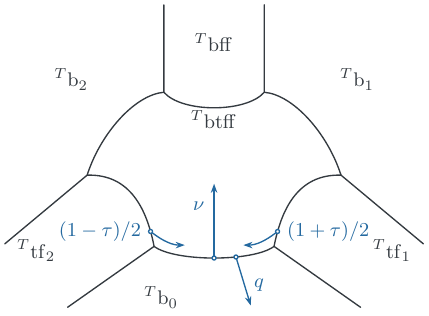}
    \caption{Coordinates on the lower constituent \btrtext{} calculus of $\dbsp{T}{\dbsctr}$.}
    \label{fig_b-transition-double-space-symmetric-coords}
\end{figure}

\begin{figure}
    \centering
    \includegraphics[width=0.5\linewidth]{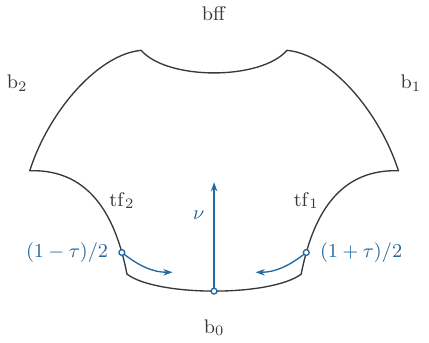}
    \caption{Coordinates on $\tff$. The coordinate $U \in \R^{n-1}$ is suppressed.}
    \label{fig_tff-middle-region-symmetric-coords}
\end{figure}

\begin{lemma}
\label{lemma_conormal-Bessel-operator-to-transformed-normal-operator}

    Suppose that $b(a) = b(a)(\nu, \tau, q)$ is conormal of order $e_0, f_1, f_2, f_{12}$ at $\bf{0}[][^{T}], \tf[1][][^{T}], \tf[2][][^{T}]$, and $\btff[][][^{T}]$, respectively. Suppose also that $b(a)$ is rapidly decaying as $q = t_1+t_2 \to \infty$. Then, for any $N > 0$,
    \begin{align}
        |\hat{a}(\nu, \tau, H)| \lesssim \nu^{e_0-1/2} \left( \frac{1+\tau}{2} \right)^{f_1-1} \left( \frac{1-\tau}{2} \right)^{f_2+n-1} |H|^{f_{12}-1/2} \left( 1+ |H| \right)^{-N} \label{eqn_conormal-transformed-normal-operator-estimate}
    \end{align}
    Moreover, this estimate is stable under the application of $(\nu \partial_{\nu})^j$, $\left((1-\tau^2) \partial_\tau \right)^{k}$, and $H^{\delta} \partial_H^{\gamma}$ for $j, k \in \N_0$, $\delta, \gamma \in \N_0^{n-1}$ satisfying $|\delta| \geq |\gamma|$.

    Moreover, suppose that $b(a)$ depends smoothly on parameters $\hat{H} \in S^{n-2}$ and $y^{(0)} \in \R^{n-1}$ and that its conormal estimates at all boundary hypersurfaces depend smoothly on $\hat{H}$ and $y^{(0)}$. Then the estimate \eqref{eqn_conormal-transformed-normal-operator-estimate} also depends smoothly on $\hat{H}$ and $y^{(0)}$.
\end{lemma}

\begin{proof}
    The assumptions imply that
    \begin{align}
        |b(a)(\nu, \tau, q)| &\lesssim \nu^{e_0} \left( \frac{1+\tau}{2} \right)^{f_1} \left( \frac{1-\tau}{2} \right)^{f_2} q^{f_{12}} \left(1 + q \right)^{-N} \label{eqn_conormal-bessel-operator-estimate}
    \end{align}
    for all $N \in \N_0$. Moreover, this estimate is stable under the application of $(\nu \partial_{\nu})^j$, $\left((1-\tau^2) \partial_\tau \right)^{k}$, and $(q \partial_q)^\ell$ for $j, k, \ell \in \N_0$. Then, the basic estimate \eqref{eqn_conormal-transformed-normal-operator-estimate} is immediate from \eqref{eqn_conormal-bessel-operator-estimate}, via the identification $|H| = q$.

    Stability under derivatives is immediate. We briefly remark on stability under $H^\delta \partial_{H}^{\gamma}$. The coordinates $H$ and $q$ are related by $q = |H|$, so $\partial_{H_i} = \hat{\eta}_i \partial_{q} + |H|^{-1} \partial_{\hat{H}_i}$. Then, because $|\delta| \geq |\gamma|$, $|H^\delta \partial_{\eta}^{\gamma} \hat{a}|$ can be controlled by $(q \partial_q)^\ell (\partial_{\hat{\eta}})^\zeta b(a)$ for bounded $\ell, \zeta$.
\end{proof}

\begin{lemma}
\label{lemma_conormal-transformed-normal-operator-to-normal-operator}
    Suppose that $\hat{a}$ satisfies the conclusion of Lemma \ref{lemma_conormal-Bessel-operator-to-transformed-normal-operator}. Additionally suppose that $f_{12} > -n + 3/2$. Then,
    \begin{align}
        |a(\nu, \tau, U)| \lesssim \nu^{e_0-1/2} \left( \frac{1+\tau}{2} \right)^{f_1-1} \left( \frac{1-\tau}{2} \right)^{f_2+n-1} |U|^{-(f_{12}+n-3/2)} \label{eqn_conormal-normal-operator-estimate}
    \end{align}
    Moreover, this estimate is stable under the application of $(\nu \partial_{\nu})^j$, $\left((1-\tau^2) \partial_\tau \right)^{k}$, and $U^{\gamma} \partial_U^{\delta}$, for $j, k \in \N_0$, $\delta, \gamma \in \N_0^{n-1}$ satisfying $|\delta| \geq |\gamma|$. That is, $a$ is conormal on $\Sigma$ up to $\Sigma \cap \bf{0}$, $\Sigma \cap \tf[1]$, and $\Sigma \cap \tf[2]$ of order $e_0 - 1/2, f_1 - 1$, and $f_2 + n - 1$, respectively.
\end{lemma}

\begin{proof}
    The Fourier transform is well-defined because $\hat{a}$ is rapidly decaying in $H$. We split the Fourier transform into two terms:
    \begin{align}
        a(\nu, \tau, U) &\lesssim \int_{|H| < |U|^{-1}} |\hat{a}(\nu, \tau, H)| \,dH + \left| \int_{|H| > |U|^{-1}} e^{i U \cdot H} \hat{a}(\nu, \tau, H) \, dH \right| \label{eqn_normal-operator-from-transformed-normal-operator-split-integral}
    \end{align}
    Handle the $|H| < |U|^{-1}$ integral by invoking Lemma \ref{lemma_conormal-Bessel-operator-to-transformed-normal-operator}. Define the quantity $P = \nu^{e_0-1/2} \left( \frac{1+\tau}{2} \right)^{f_1-1} \left( \frac{1-\tau}{2} \right)^{f_2+n-1}$. Then,
    %-- Deprecated 5 Sep 26 - old normalisation
    % Handle the $|H| < |U|^{-1}$ integral by invoking Lemma \ref{lemma_conormal-Bessel-operator-to-transformed-normal-operator}. Let $P = \nu^{e_0} \left( \frac{1+\tau}{2} \right)^{f_1} \left( \frac{1-\tau}{2} \right)^{n+f_2}$. Then,
    %--
    \begin{align*}
        \int_{|H| < |U|^{-1}} |\hat{a}(\nu, \tau, H)| \,dH &\lesssim P \int_{|H| < |U|^{-1}} |H|^{f_{12} -1/2} \left( 1 + |H| \right)^{-N}  \,dH \\
        &\lesssim P \int_0^{|U|^{-1}} R^{f_{12} - 1/2 + n - 2} \,dR \\
        &\approx P |U|^{-(f_{12} + n - 3/2)}
    \end{align*}
    
    Handle the $|H| > |U|^{-1}$ integral by integrating by parts. Write $e^{i U \cdot H} = \left(\frac{U \cdot D_H}{|H|^2} \right)^{L} e^{i U \cdot H}$ for large $L \in \N$ to be determined. Then integrate by parts $L$ times. The boundary term produced by the $\ell$th integration by parts is controlled by 
    \begin{align*}
        & \int_{|H| = |U|^{-1}} \left| \left(\frac{U \cdot D_H}{|H|^2} \right)^{L - \ell} e^{i U \cdot H} \frac{|U|}{|H|^2} \left(\frac{U \cdot D_H}{|H|^2} \right)^{\ell-1} \hat{a} \right| \,dH \\
        &= \int_{|H| = |U|^{-1}} \left| \frac{1}{|U|} \left(\frac{U \cdot D_H}{|U|^2} \right)^{\ell-1} \hat{a} \right| \,dH
    \end{align*}
    By Lemma \ref{lemma_conormal-Bessel-operator-to-transformed-normal-operator},
    \begin{align*}
        \left| \left(\frac{U \cdot D_H}{|U|^2} \right)^{\ell-1} \hat{a} \right| \lesssim \frac{1}{|U|^{\ell-1} |H|^{\ell-1}} |H|^{f_{12}-1/2} \left( 1 + |H| \right)^{-N} P
    \end{align*}
    Restricting to $|H| = |U|^{-1}$, it follows then that the $\ell$th boundary term is controlled by
    \begin{align*}
        P |U|^{-(f_{12}+1/2)} \mathrm{vol}\left\{ |H| = |U|^{-1} \right\} \approx P |U|^{-(f_{12} + n - 3/2)}.
    \end{align*}
    Now consider the bulk part of the $|H| > |U|^{-1}$ integral in \eqref{eqn_normal-operator-from-transformed-normal-operator-split-integral}. Again invoke Lemma \ref{lemma_conormal-Bessel-operator-to-transformed-normal-operator} to obtain
    \begin{align*}
        &\int_{|H| > |U|^{-1}} \left| \left(\frac{U \cdot D_H}{|H|^2} \right)^{L} \hat{a} \right| \,dH \\
        &\lesssim P |U|^{-L} \int_{|H| > |U|^{-1}} \frac{1}{|H|^L} |H|^{f_{12}-1/2} \left( 1 + |H| \right)^{-N} \,dH \\
        &\lesssim P |U|^{-L} \int_{|U|^{-1}}^\infty R^{f_{12} - 1/2 - L + n - 2} \,dR \\
        &\approx P |U|^{-(f_{12}+n-3/2)}
    \end{align*}
    assuming that $L$ is chosen to be sufficiently large.

    We find that every term in the split Fourier transform obeys the estimate \eqref{eqn_conormal-normal-operator-estimate}. We prove that conormal derivatives of $a$ satisfy the same estimate by differentiating under the integral sign.
\end{proof}

Under conormality assumptions on $b(a)$ up to $\bf{0}[][^{T}], \tf[1][][^{T}], \tf[2][][^{T}]$, and $\btff[][][^{T}]$, we have proven estimates for $a$ on $\Sigma$ up to $\bf{0}, \tf[1]$, and $\tf[2]$, but not yet up to $\btff$.

Fix this with a change of coordinates. Introduce $Z = \frac{x_1 + x_2}{|Y|}$ and $\hat{Y} = \frac{Y}{|Y|}$, where $Y = y_1 - y_2$. Then $(\nu, \tau, Z, \hat{Y},|Y|,y_2)$ is a coordinate chart on $\dbsp{M}{\gas}$. In this chart, $|Y|$ is a defining function for $\tff$. Restricting to $\Sigma$ yields the chart $(\nu, \tau, Z, \hat{Y})$. This chart is related to the old chart by the relation $U = \hat{Y}/Z$. The chart is smooth on $\Sigma$ up to $\bf{0}, \tf[1], \tf[2]$, and $\btff$.

\begin{figure}
    \centering
    \includegraphics[width=0.5\linewidth]{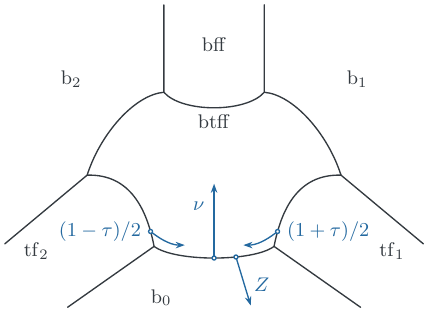}
    \caption{Coordinates on $\tff$ smooth up to $\btff$.}
    \label{fig_tff-lower-region-symmetric-coords}
\end{figure}

\begin{lemma}
\label{lemma_conormal-transformed-normal-operator-to-normal-operator-at-btff}
    Suppose $a = a(\nu, \tau, U)$ is a distribution on $\Sigma$ satisfying the conclusion of Lemma \ref{lemma_conormal-transformed-normal-operator-to-normal-operator}. Define $a^{\mathrm{II}}(\nu, \tau, Z, \hat{Y}) = a(\nu, \tau, \hat{Y}/Z)$. Then $a^{\mathrm{II}}$ satisfies
    \begin{align}
        |a^{\mathrm{II}}(\nu, \tau, Z, \hat{Y})| \lesssim \nu^{e_0-1/2} \left( \frac{1+\tau}{2} \right)^{f_1-1} \left( \frac{1-\tau}{2} \right)^{f_2+n-1} Z^{f_{12}+n-3/2} \label{eqn_conormal-normal-operator-estimate-at-btff}
    \end{align}
    Moreover, this estimate is stable under the application of $(\nu \partial_\nu)^j, \left((1-\tau^2) \partial_\tau \right)^k$, $(Z \partial_Z)^\ell$, $\partial_{\hat{Y}}^\delta$ for $j, k, \ell \in \N_0$, $\delta \in \N_0^{n-1}$. That is, $a^\mathrm{II}$ is conormal on $\Sigma$ up to $\Sigma \cap \bf{0}$, $\Sigma \cap \tf[1]$, $\Sigma \cap \tf[2]$, and $\Sigma \cap \btff$ of order $e_0 - 1/2, f_1 - 1$, $f_2 + n - 1$, and $f_{12} + n - 3/2$, respectively.
\end{lemma}
\begin{proof}
    This is an immediate consequence of Lemma \ref{lemma_conormal-transformed-normal-operator-to-normal-operator}. Stability under $\partial_{\hat{Y}}$ follows from stability under $U$-derivatives by writing $\partial_{\hat{Y}} = \frac{1}{Z} \partial_U = \frac{1}{|U|} \partial_U$. Stability under $Z \partial_Z$ follows from the observation that $Z = 1/|U|$.
\end{proof}

Repeating the above arguments yields conormal estimates up to all boundary hypersurfaces of $\Sigma$ in the middle region.

\begin{lemma} \label{lemma_conormal-Bessel-operator-to-normal-operator-all-charts}
Suppose that $b(a) = b(a)(t_0,t_1,t_2)$ is conormal on $\dbsp{T}{\dbsctr}$ of order $\mathfrak{e}$ at the boundary hypersurfaces $\bf{0}[][^{T}], \bf{1}[][^{T}], \bf{2}[][^{T}], \btff[][][^{T}], \tf[1][][^{T}], \tf[2][][^{T}]$, and $\bff[][][^{T}]$, where $\mathfrak{e} = (e_0, e_1, e_2, f_{12}, f_1, f_2, g_{12})$. Suppose also that $b(a)$ is rapidly decaying as $q = t_1+t_2 \to \infty$. Define, as in \eqref{eqn_bessel-kernel-to-fourier-transformed-kernel}, $\hat{a}(\rho, s, \eta) = |\eta|^{-1/2} \rho^{-1/2} s^{-1/2} b(a)(\rho|\eta|, s|\eta|, |\eta|)$ and $a(\rho,s,u) = \int e^{i u \cdot \eta} \hat{a}(\rho,s,\eta) \,d\eta$. Then, $a$ is conormal of order $\mathfrak{e} + (-1/2, -1/2, n-1/2, n-3/2, -1, n-1, n-1)$ on $\Sigma$ at the boundary hypersurfaces $\bf{0}[], \bf{1}[], \bf{2}[], \btff[][], \tf[1][], \tf[2][]$, and $\bff[][][]$, respectively.
\end{lemma}
\begin{proof}
Define a family of charts on $\Sigma$, similar to the $(\nu, \tau, U)$ chart, that cover the middle region. In each chart, repeat the arguments of Lemmas \ref{lemma_conormal-Bessel-operator-to-transformed-normal-operator}, \ref{lemma_conormal-transformed-normal-operator-to-normal-operator}, and \ref{lemma_conormal-transformed-normal-operator-to-normal-operator-at-btff}.
\end{proof}

Assuming $b(a)$ is improved to polyhomogeneous, we can also prove $a$ is polyhomogeneous.

\begin{proposition} \label{proposition_phg-Bessel-operator-to-normal-operator-all-charts}
Suppose that $b(a) = b(a)(t_0,t_1,t_2)$ is polyhomogeneous on $\dbsp{T}{\dbsctr}$. Let $E_0, E_1, E_2, F_{12}, F_1, F_2$, and $G_{12}$ denote its index sets at the boundary hypersurfaces $\bf{0}[][^{T}], \bf{1}[][^{T}], \bf{2}[][^{T}], \btff[][][^{T}], \tf[1][][^{T}], \tf[2][][^{T}]$, and $\bff[][][^{T}]$, respectively. Suppose also that $b(a)$ is rapidly vanishing at all boundary hypersurfaces where $q = t_1+t_2 = \infty$. Define $a$ as in Lemma \ref{lemma_conormal-Bessel-operator-to-normal-operator-all-charts}. Then $a$ satisfies polyhomogeneous estimates on $\Sigma$ at the boundary hypersurfaces $\bf{0}[], \bf{1}[], \bf{2}[], \btff[][], \tf[1][], \tf[2][]$, and $\bff[][][]$. Its index sets are $E_0 - 1/2, E_1 - 1/2, E_2 + n - 1/2, F_{12} + n - 3/2, F_1 - 1, F_2 + n - 1$, and $G_{12} + n - 1$, respectively.
\end{proposition}
\begin{proof}
As in the proof of Lemma \ref{lemma_conormal-Bessel-operator-to-normal-operator-all-charts}, one defines a family of charts that covers the middle region of $\Sigma$. We sketch a proof of polyhomogeneity in the $(\nu, \tau, U)$ chart, which is smooth up to $\bf{0}, \btff, \tf[1]$, and $\tf[2]$. As before, use the corresponding chart $(\nu, \tau, q)$ on $\dbsp{T}{\dbsctr}$.

Write $b(a) = b(a)(\nu, \tau, q)$. Because $b(a)$ is polyhomogeneous, the estimate \eqref{eqn_conormal-bessel-operator-estimate} improves to
\begin{align}
    \left| B_{E_0, C}(\nu \partial_\nu) b(a)(\nu, \tau, q)\right| &\lesssim \nu^{C} \left( \frac{1+\tau}{2} \right)^{f_1} \left( \frac{1-\tau}{2} \right)^{f_2} q^{f_{12}} \left(1 + q \right)^{-N} \label{eqn_phg-bessel-operator-estimate}
\end{align}
where
\begin{align*}
    B_{E_0, C}(\nu \partial_\nu) = \prod_{\substack{(z,p) \in E_0 \\ \mathrm{Re} z \leq C}} (\nu \partial_\nu - z)
\end{align*}
in addition to corresponding estimates that capture the expansions of $b(a)$ in $(1+\tau), (1-\tau)$, and $q$. These estimates hold for all $C \in \R$ and remain stable under the application of $(\nu \partial_{\nu})^j$, $\left((1-\tau^2) \partial_\tau \right)^{k}$, and $(q \partial_q)^\ell$ for $j, k, \ell \in \N_0$.

To prove that $a$ is polyhomogeneous, we must improve \eqref{eqn_conormal-normal-operator-estimate-at-btff} to a polyhomogeneous estimate. For example, the desired polyhomogeneous estimate in $\nu$ is
\begin{align}
    \left| B_{E_0, C}(\nu\partial_\nu) a^{\mathrm{II}}(\nu, \tau, Z, \hat{Y})\right| \lesssim \nu^{C-1/2} \left( \frac{1+\tau}{2} \right)^{f_1-1} \left( \frac{1-\tau}{2} \right)^{f_2+n-1} Z^{f_{12}+n-3/2}. \label{eqn_phg-normal-operator-estimate-at-btff}
\end{align}
Follow essentially the same steps as in Lemmas \eqref{lemma_conormal-Bessel-operator-to-transformed-normal-operator}, \eqref{lemma_conormal-transformed-normal-operator-to-normal-operator}, and \eqref{lemma_conormal-transformed-normal-operator-to-normal-operator-at-btff}, except every time we would invoke \eqref{eqn_conormal-bessel-operator-estimate}, we invoke \eqref{eqn_phg-bessel-operator-estimate} instead.
\end{proof}

\begin{corollary}
\label{corollary_model-problem-tff-lower-solution-is-phg}
The solution $n(Q_1)$ to the model problem \eqref{eqn_model-problem-tff-function-form} satisfies polyhomogeneous estimates on $\Sigma$ at the boundary hypersurfaces $\bf{0}[], \bf{1}[], \bf{2}[], \btff[][], \tf[1][], \tf[2][]$, and $\bff[][][]$. Its index sets are generated by
\begin{itemize}
    \item $(k(\beta+1),0)$ for $k \in \N_0$ at $\bf{0}[][]$
    \item $(1,0)$ at $\bf{1}[][]$
    \item $(n-1,0)$ at $\bf{2}[][]$
    \item $(n-1 + k(\beta+1),0)$ for $k \in \N_0$ at $\btff[][][]$
    \item $((k+1)(\beta+1),0)$ for $k \in \N_0$ at $\tf[1][][]$
    \item $(n-1+k(\beta+1),0)$ for $k \in \N_0$ at $\tf[2][][]$
    \item $(n-1,0)$ at $\bff[][][]$.
\end{itemize}
If $\beta$ is an odd positive integer, then these index sets are modified by increasing some logarithmic exponents.
\end{corollary}
\begin{proof}
Combine Proposition \eqref{eqn_Bessel-model-problem-solution-is-phg} and Proposition \eqref{proposition_phg-Bessel-operator-to-normal-operator-all-charts}.
\end{proof}

Combined with the previous results for the upper model problem, we finally resolve the model problem at $\tff$. Recall that, although we have notationally suppressed it, $n(Q_1) = n(Q_1;y^{(0)})$ is a distribution on a fiber $\Sigma$ of $\tff$ with basepoint $y^{(0)}$. All of our polyhomogeneous estimates vary smoothly in $y^{(0)}$, and hence $n(Q_1)$ defines a polyhomogeneous distribution on $\tff$ by allowing $y^{(0)}$ to vary.

\begin{corollary}
\label{corollary_laplacian-parametrix-with-vanishing-remainder-at-tff}
    There exists an index family $\mathcal{E}$ on $\dbsp{M}{\gas}$ and an operator $Q_1 \in \Psi^{-\infty,\mathcal{E}}_{\gas}(M)$, such that
    \begin{align*}
        P(Q_0 + Q_1) = I - R_1
    \end{align*}
    where $R_1 \in \Psi^{-\infty,\mathcal{R}_1}_{\gas}(M)$.

    The index family $\mathcal{E}$ is generated by
    \begin{itemize}
        \item $(k(\beta+1),0)$ for $k \in \N_0$ at $\bf{0}$
        \item $(1,0)$ at $\bf{1}$
        \item $(-1,0)$ at $\bf{2}$
        \item $(0,0)$ at $\tff$
        \item $(-1 + k(\beta+1),0)$ for
        $k \in \N_0$ at $\btff[]$
        \item $(0,0)$ at $\zff$
        \item $(-1,0)$ at $\bff$
        \item $((k+1)(\beta+1),0)$ for $k \in \N_0$ at $\tf[1]$
        \item $(-1+k(\beta+1),0)$ for $k \in \N_0$ at $\tf[2]$.
    \end{itemize}
    If $\beta$ is an odd positive integer, then these index sets are modified by increasing some logarithmic exponents.
    
    The index family $\mathcal{R}_1$ is identical to $\mathcal{E}$ except that it is generated at $\tff$ by $(1,0)$.
\end{corollary}
\begin{proof}
Define $Q_1$ such that its restriction to $\tff$ is $n(Q_1)$, e.g. by choosing a collar neighborhood of $\tff$ and extending $n(Q_1)$ to be constant in the normal variable. The difference in index sets between $n(Q_1)$ and $Q_1$ is owed to the difference between the half-density factor in \eqref{eqn_normal-operator-tff-normalization} and the half-density $\sigma$.

Recall that $PQ_0 = I - R_0$. By construction, $N_{\tff}(P; y^{(0)}) N_{\tff}(Q_1; y^{(0)}) = N_{\tff}(R_0; y^{(0)})$, so $R_1 := R_0 - PQ_1$ vanishes to first order at $\tff$. As $R_0$ and $P$ are in the small calculus, the composition law for $\Psi^{\ast,\ast}_{\gas}$ implies that $R_1$ has the same index set as $Q_1$ at other faces. 
\end{proof}

\section{Upper front face}

Having solved the model problem at the transition front face, we next solve the model problem at the upper front face $\zff$. We shall solve for $Q_2$ such that
\begin{align*}
    P(Q_0 + Q_1 + Q_2) = I - R_2
\end{align*}
where $R_2$ has index set $(1,0) + \N_0$ at both $\tff$ and $\zff$.

The upper front face $\uff$ is naturally diffeomorphic to $[0,1) \times \ff(\dbsp{M}{0'})$, where $\ff(\dbsp{M}{0'})$ is the front face of the overblown $0$-double space $\dbsp{M}{0'}$ \cite{hintz_0-calculus-mazzeo-melrose_2023}. The diffeomorphism is given by restricting the natural map $\beta^{012}_{0,12} : \anysp{M}{012} = \dbsp{M}{\gas} \to \anysp{M}{0,12} := [0,1) \times \dbsp{M}{0'}$. This map is given by blowing down all blown-up faces of $\dbsp{M}{\gas}$ contained in $\{ \epsilon = 0 \}$.

We choose a coordinate system near $\zff$. As we did for the model problem at $\tff$, fix $y^{(0)} \in \partial M$, a boundary-defining function $x$, and a Riemann normal coordinate system $y$ for $(\partial M, h)$ at $y^{(0)}$ so that $h(y^{(0)},dy) = dy^2$. Use the coordinate system $(\epsilon, s, u, \rho', y_2)$ for $\dbsp{M}{\gas}$, where $s = x_1/x_2, u = (y_1 - y_2)/x_2, \rho' = x_2/\epsilon$. This chart is smooth near the intersection of $\zff, \tff$, and $\bf{1}$. (See Figure \ref{fig_double-space-upper-region-coords}.) In this chart, $\rho'$ is a defining function for $\zff$. The chart induces a chart $(\epsilon, s, u)$ for $\zff$, in which $\epsilon$ is a boundary-defining function for $\tff \cap \zff$.

In these coordinates, $P$ lifts to
\begin{align*}
    P = (s \partial_s)^2 - \left( \beta\frac{s}{s+1/\rho'} \right) s\partial_s + s^2 \Delta_u
\end{align*}
Taking $\rho' \to 0$, we find that
\begin{align}
    N_{\zff}\left(P; y^{(0)}\right) &= s^2 (\partial_s^2 + \Delta_u). \label{eqn_normal-operator-zff-gg-laplacian}
\end{align}

In particular, observe that this normal operator is independent of $\epsilon$. It can be viewed as acting on the model space $\tilde{M} = \R^n_1$, parametrized by $(s,u)$.

For an operator $A$ on $M$, let $\mathcal{K}_A$ denote the normalization of its kernel with respect to the density $\epsilon^{-n/2} (\rho')^{-(n+1)/2} |d\epsilon ds du d\rho' dy_2|^{1/2}$, which is locally equivalent to the density $\sigma$. Write $\mathcal{K}_{A} = \mathcal{K}_{A}(\epsilon, s, u, \rho', y_2)$ and write its restriction to $\zff$ as $n_{\zff}\!\left(A;y^{(0)}\right)(\epsilon, s, u) := \mathcal{K}_{A}|_{\zff}(\epsilon, s, u, y^{(0)}) = \mathcal{K}_{A}\left(\epsilon, s, u, 0, y^{(0)}\right)$. Then, the model problem at $\zff$ (suppressing the $y^{(0)}$-dependence) can be written as
\begin{align}
    N_{\zff}\left( P \right) n_{\zff}\!\left(Q_1 \right)(\epsilon, s,u) &= n_{\zff}\!\left( R_1 \right)(\epsilon, s,u) \label{eqn_model-problem-zff-function-form}
\end{align}
By the standard parametrix construction for the 0-calculus \cite{mazzeo_differential-edge-operators-I_1991} \cite{hintz_0-calculus-mazzeo-melrose_2023}, $N_{\zff}\left(P\right) : sH^2_0 \to sL^2$ has an inverse $N_{\zff}\left(P\right)^{-1}$ that lies in $\Psi^{-2,\ast}_0(\tilde{M})$, the $0$-calculus on $\tilde{M}$. Then set $n_{\zff}\!\left(Q_2 \right) = N_{\zff}\left(P\right)^{-1} n_{\zff}\!\left( R_1 \right)$.

Because $N_{\zff}\left(P\right)^{-1}$ lies in a $0$-calculus, we can compute the index family of $n_{\zff}\!\left(Q_2 \right)$ from that of $n_{\zff}\!\left( R_1 \right)$. In particular, because $R_1$ has index set $(1,0) + \N_0$ at $\tff$ by construction, $n_{\zff}\!\left( R_1 \right)$ has index set $(1,0) + \N_0$ at $\tff \cap \zff$. Because $\epsilon$ is a boundary-defining function for $\tff \cap \zff$ and $N_{\zff}\left(P\right)^{-1}$ is $\epsilon$-independent, it follows that the same is true of $n_{\zff}\!\left( Q_2 \right)$.

Then extend $n_{\zff}\!\left( Q_2 \right)$ off $\zff$ to obtain $Q_2$, which has index set $\N_0$ at $\zff$ (by construction) and index set $(1,0) + \N_0$ at $\tff$. Setting $R_2 = R_1 - PQ_2$, we write
\begin{align*}
    P(Q_0 + Q_1 + Q_2) = I - R_2
\end{align*}
By \eqref{eqn_model-problem-zff-function-form}, $R_2$ has index set $(1,0) + \N_0$ at $\zff$. Because both $R_1$ and $Q_2$ have index set $(1,0) + \N_0$ at $\tff$, the same is true of $R_2$. We have therefore obtained a right parametrix $Q^{(2)} = Q_0 + Q_1 + Q_2$ for $P$ with remainder term vanishing at both $\tff$ and $\zff$.

We summarize our result as follows:
\begin{proposition}
\label{proposition_laplacian-parametrix-with-vanishing-remainder-at-zff}

There exist operators $Q_0 \in \Psi^{-2}_{\gas}(M)$ and $Q^{(2)} \in \Psi^{-\infty,\mathcal{E}}_{\gas}(M)$, such that
\begin{align*}
    P(Q_0 + Q^{(2)}) = I - R_2
\end{align*}
where $\mathcal{E}$ is as in Corollary \ref{corollary_laplacian-parametrix-with-vanishing-remainder-at-tff} and $R_2 \in \Psi^{-\infty,\mathcal{R}_2}_{\gas}(M)$, where $\mathcal{R}_2$ is identical to $\mathcal{E}$ except that it is generated by $(1,0)$ at $\tff$ and $\zff$.
\end{proposition}

\section{b-front faces}

The small-calculus parametrix $Q_0$ satisfies $PQ_0 = I - R_0$, where $R_0 \in \Psi^{-\infty}_{\gas}(M)$. In particular, $R_0$ is rapidly vanishing at $\btff$ and $\bff$. However, this property is not retained by the error term $R_2$, which has leading order $-1$ at both of these boundary hypersurfaces. By Proposition \ref{proposition_gas-giant-operator-compactness}, the leading order at both faces must be strictly greater than $-1$ for the remainder term to be compact.

We resolve this by solving model problems at $\btff$ and $\bff$. Work in the chart $(\rho, s, z_2, Y, y_2)$ in the lower region of $\dbsp{M}{\gas}$, where $Y = y_1 - y_2$ and $z_2 = x_2/|Y|$. (See Figure \ref{fig_double-space-lower-region-coords}.) This chart is smooth up to $\bf{0}, \bf{1}, \tff$, and $\btff$; these faces have the boundary-defining functions $\rho, s, |Y|$, and $z_2$, respectively. Lifting $P$ to $\dbsp{M}{\gas}$ in this chart, fixing $(Y, y_2) = (Y^{(0)}, y^{(0)})$, and taking $z_2 \to 0$, we obtain
\begin{align*}
    N_{\btff}(P; Y^{(0)}, y^{(0)}) &= s^2 \partial_s^2 - \beta \frac{s}{s+\rho} s \partial_s
\end{align*}
which is to be viewed as an operator acting on the fiber $(Y, y_2) = (Y^{(0)}, y^{(0)})$ of $\btff$. The fiber admits a further fibration by holding $\rho$ constant. The normal operator acts tangent to these $\rho$-level sets.

\begin{figure}
    \centering
    \includegraphics[width=0.5\linewidth]{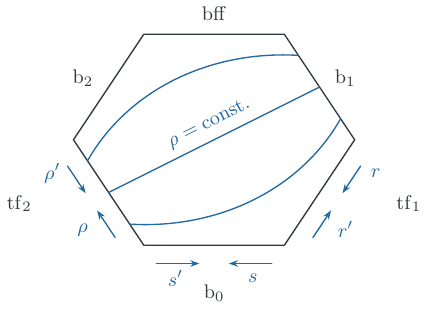}
    \caption{A fiber of $\btff$ over $(Y^{(0)}, y^{(0)})$ and its $\rho$-level sets.}
    \label{fig_btff-fibration}
\end{figure}

Our model problem is
\begin{align}
    N_{\btff}(P; Y^{(0)}, y^{(0)}) N_{\btff}(Q_3; Y^{(0)}, y^{(0)}) = N_{\btff}(R_2; Y^{(0)}, y^{(0)}) \label{eqn_model-problem-btff-operator-form}
\end{align}
Identify the fiber over $(Y^{(0)}, y^{(0)})$ with $[0,\infty)^2_b$, parametrized by $(\rho,s)$. We write $n_{\btff}(A; Y^{(0)}, y^{(0)})(\rho,s)$ to denote the restriction of the scalar-valued kernel of $A$, normalized as in \eqref{eqn_normal-operator-tff-normalization}, to this fiber. Rewrite the problem as
\begin{align*}
    N_{\btff}(P; Y^{(0)}, y^{(0)}) n_{\btff}(Q_3; Y^{(0)}, y^{(0)}) = n_{\btff}(R_2; Y^{(0)}, y^{(0)})
\end{align*}
On each $\rho$-level set, the problem has a unique solution $n_{\btff}(Q_3; Y^{(0)}, y^{(0)})$ satisfying a Dirichlet boundary condition at $s = 0$ and decaying as $s \to \infty$. We must prove that this solution is polyhomogeneous on $\btff$.

By identifying the $(Y^{(0)}, y^{(0)})$-fiber with $[0,\infty)^2_b$, we view $N_{\btff}(P; Y^{(0)}, y^{(0)})$ as an operator in the b-transition calculus on $[0,\infty)$. Writing an explicit Green's function with Dirichlet boundary conditions and invoking mapping properties of the b-transition calculus (cf. \cite{mcdonald_laplacian-cone-singularities_1990} \cite{mazzeo-melrose_analytic-surgery_1995}) proves that $n_{\btff}(Q_3; Y^{(0)}, y^{(0)})$ is polyhomogeneous on $[0,\infty)^2_b$. On $\btff$, this corresponds to polyhomogeneous estimates up to $\bf{1}, \bf{0}$, and $\tf[1]$.

A similar technique yields polyhomogeneous estimates up to $\bf{1}, \bf{2}$, and $\bff$. Do this by choosing the coordinate chart $(\rho', r, z_0, Y, y_2)$ for $\dbsp{M}{\gas}$, which restricts to the chart $(\rho', r)$ for the fiber of $\btff$ over $(Y^{(0)}, y^{(0)})$. In this chart, the normal operator is $N_{\btff}(P; Y^{(0)}, y^{(0)}) = r^2 \partial_r^2 - \beta \frac{r}{r+1} r \partial_r$. It acts on $[0,\infty)^2_b$, parametrized by $(\rho', r)$. Again, using an explicit Green's function and invoking a mapping formula yields polyhomogeneity of the solution to the model problem on $[0,\infty)^2_b$, which corresponds to estimates up to $\bf{1}, \bf{2}$, and $\bff$.

It remains to obtain polyhomogeneous estimates up to $\bf{0}, \bf{2}$, and $\tf[2]$. Here, the problem is especially simple. Near $\bf{2} \cap \tf[2]$, use the chart $(\rho', r')$ and observe that $N_{\btff}(P; Y^{(0)}, y^{(0)}) = (r')^2 \partial_{r'} + \big( \beta \frac{1}{1+r'} + 2 \big) r' \partial_{r'}$, a $\rho'$-independent operator with indicial roots $-(\beta+1), 0$. Near $\bf{0} \cap \tf[2]$, use the chart $(\rho, s')$ and observe that $N_{\btff}(P; Y^{(0)}, y^{(0)}) = (s')^2 \partial_{s'} + \big( \beta \frac{1}{1+\rho s'} + 2 \big) s' \partial_{s'}$. The operator is dependent on $\rho$, but its indicial roots are not; they are identically $-(\beta+1)$ and $0$. The decay condition $s \to \infty$ imposed on the solution $n_{\btff}(Q_3; Y^{(0)}, y^{(0)})$ is equivalent to selecting a solution corresponding to the upper indicial root $0$.

Putting all of this together, the model problem has a polyhomogeneous solution $N_{\btff}(Q_3; Y^{(0)}, y^{(0)})$. By allowing $Y^{(0)}, y^{(0)}$ to vary again, write $N_{\btff}(Q_3)$ for the solution viewed as a distribution on $\btff$. As $R_2$, normalized by $\sigma$, vanishes to first order at $\tff$ by construction, the right-hand side of \eqref{eqn_model-problem-btff-operator-form} vanishes to first order in $|Y|$, which is a boundary-defining function for $\tff$. As the normal operator has no dependence on the basepoint $(Y^{(0)}, y^{(0)}$, it follows that $N_{\btff}(Q_3)$ is also vanishing to first order at $\tff$.

Extend the solution off $\btff$ to obtain a kernel $Q_3$. As $\zff$ is disjoint from $\btff$, we can safely assume that $Q_3$ is rapidly vanishing at $\zff$. Setting $R_3 = R_2 - PQ_3$, we obtain
\begin{align*}
    P(Q_0 + Q_1 + Q_2 + Q_3) &= I - R_3
\end{align*}
The kernel $Q_3$ is (normalized by $\sigma$) vanishing to first order at $\tff$, so $PQ_3$ is vanishing to first order at $\tff$. It follows from this and the rapid vanishing of $Q_3$ at $\zff$ that $R_3$ retains the property of $R_2$ of having index set $(1,0) + \N_0$ at $\tff$ and $\zff$. Finally, because $N_{\btff}(P) N_{\btff}(Q_3) = N_{\btff}(R_2)$, the leading term of $R_2$ is killed off.

It remains to solve away the error at $\bff$, but this is even simpler because the normal operator is the indicial operator of $P$, which is constant in $\epsilon$. Namely, $N_{\bff}(P; Y^{(0)}, y^{(0)}) = s^2 \partial_s^2$.

We summarize our results as:

\begin{theorem}
\label{theorem_laplacian-parametrix-with-compact-remainder-phg-index-family}

There exist operators $Q_0 \in \Psi^{-2}_{\gas}(M), Q^{(4)} \in \Psi^{-\infty,\mathcal{E}'}_{\gas}(M)$, such that
\begin{align*}
    P\left(Q_0 + Q^{(4)}\right) = I - R_4
\end{align*}
where $\mathcal{E}'$ agrees with the index family $\mathcal{E}$ of Corollary \ref{corollary_laplacian-parametrix-with-vanishing-remainder-at-tff}, except some logarithmic exponents are increased, and $R_4 \in \Psi^{-\infty,\mathcal{R}_4'}_{\gas}(M)$, where $\mathcal{R}_4'$ is identical to the index family $\mathcal{R}_2$ of Proposition \ref{proposition_laplacian-parametrix-with-vanishing-remainder-at-zff}, except that it is generated by $(0,0)$ at $\bff$ and by $\{ (0, 0) \} \cup \{ ((k+1)(\beta+1),0) : k \in \N_0 \}$ at $\btff$, and that possibly some logarithmic exponents are increased. In particular, $(R_4)_\epsilon : x_\epsilon^{1-\alpha} x L^2(dV_{g_\epsilon}) \to x_\epsilon^{1-\alpha} x L^2(dV_{g_\epsilon})$ is a family of compact operators for $\epsilon \in [0,1)$.
\end{theorem}

This proves the existence of a right parametrix as described in Theorem \ref{theorem_laplacian-parametrix-with-compact-remainder}. The left parametrix construction can be done by applying the same steps to the adjoint of $P$ with respect to $x_\epsilon^{1-\alpha} x L^2(dV_{g_\epsilon})$.

In fact, under the assumption that $P_0 = x^{2-\alpha} \Delta_g : x^{2-\alpha} H^2_0(dV_g) \to x^{2-\alpha} L^2(dV_g)$ is invertible, we can ensure that the remainder term is also vanishing to first order at $\bf{0}$.

\begin{theorem}
\label{theorem_laplacian-parametrix-with-compact-remainder-vanishing-at-b0-phg-index-family}

Suppose that $P_0 = x^{2-\alpha} \Delta_g : x^{2-\alpha} H^2_0(dV_g) \to x^{2-\alpha} L^2(dV_g)$ is invertible. Then there exist operators $Q_0 \in \Psi^{-2}_{\gas}(M), \check{Q}^{(4)} \in \Psi^{-\infty,\mathcal{E}'}_{\gas}(M)$, such that
\begin{align*}
    P\left(Q_0 + \check{Q}^{(4)}\right) = I - \check{R}_4
\end{align*}
where $\mathcal{E}'$ agrees with the index family $\mathcal{E}$ of Corollary \ref{corollary_laplacian-parametrix-with-vanishing-remainder-at-tff}, except some logarithmic exponents are increased, and $\check{R}_4 \in \Psi^{-\infty,\check{\mathcal{R}}_4'}_{\gas}(M)$, where $\check{\mathcal{R}}_4'$ is identical to the index family $\mathcal{R}_4'$ of Theorem \ref{theorem_laplacian-parametrix-with-compact-remainder-phg-index-family}, except that it is generated at $\bf{0}$ by $\{ (1, 0) \} \cup \{ ((k+1)(\beta+1), 0) : k \in \N_0 \}$, and that possibly some logarithmic exponents are increased. In particular, $(\check{R}_4)_\epsilon : x_\epsilon^{1-\alpha} x L^2(dV_{g_\epsilon}) \to x_\epsilon^{1-\alpha} x L^2(dV_{g_\epsilon})$ is a family of compact operators for $\epsilon \in [0,1)$.
\end{theorem}
\begin{proof}
    Let $Q_0, Q^{(4)}, R_4$ be as in Theorem \ref{theorem_laplacian-parametrix-with-compact-remainder-phg-index-family}. We must solve away the leading term of $R_4$ at $\bf{0}$. This is possible as long as the normal operator of $P$ at $\bf{0}$ is invertible. But $\bf{0}$ is identified with the 0-double space $\dbsp{M}{0'}$ and the normal operator is precisely $P_0$ itself. Then $P_0^{-1}$ lies in the 0-calculus on $M$ \cite{mazzeo_differential-edge-operators-I_1991} \cite{hintz_0-calculus-mazzeo-melrose_2023} and the solution to the model problem lifts off $\bf{0}$ to a kernel in $\Psi^{-\infty}_{\gas}(M)$.
\end{proof}

\section{Solving away the error at side faces}

Having solved for a parametrix with error term vanishing at all faces with leading order 1, our aim is to iterate the error away to infinite order. We shall do so with a Neumann series argument. In preparation, we must solve away the error at the side faces $\tf[1]$ and $\bf{1}$ to infinite order, to ensure that the Neumann series converges in the \ztrtext{} calculus. We give a technique that assumes the hypotheses of Theorem \ref{theorem_laplacian-parametrix-with-compact-remainder-vanishing-at-b0-phg-index-family}

Begin with the face $\tf[1]$. Work in one of two coordinate charts for $\dbsp{M}{\gas}$: $(h, r, u, x_2, y_2)$ or $(h,r', u, x_2, y_2)$, where $h = s + \rho$. These are slight modifications to charts depicted in Figure \ref{fig_double-space-middle-region-coords}; the reason is that $h$ is a global boundary-defining function for $\tf[1]$.

We shall solve $P Q_5 = R_4$ modulo polyhomogeneous distributions vanishing to infinite order at $\tf[1]$. The action of $P$ on polyhomogeneous expansions at $\tf[1]$ is, to leading order, governed by its indicial operator \cite{mazzeo_differential-edge-operators-I_1991}. After the same change of coordinates as in \eqref{eqn_normal-operator-tff-gg-laplacian-general-form-checked-chart}, this operator can be identified with
\begin{align*}
    I_{\tf[1]}(P) &= \left( \frac{r}{r+1} h\partial_h + r\partial_r \right)^2 - \left( \frac{\beta r}{r+1} + 1 \right) \left( \frac{r}{r+1} h\partial_h + r\partial_r \right) \\
    &= \left( \frac{1}{1+r'} h\partial_h - r' \partial_{r'} \right)^2 - \left( \frac{\beta}{1+r'} + 1 \right) \left( \frac{1}{1+r'} h\partial_h - r' \partial_{r'} \right)
\end{align*}
Given a polyhomogeneous expansion $f = \sum_{(z,p) \in E} a_{z,p} h^z (\log h)^p = a_{z_0,p_0} h^{z_0} (\log h)^{p_0} + \mathrm{h.o.t.}$, $P$ acts on it by
\begin{align*}
    Pf &= I_{\tf[1]}(P)(a_{z_0,p_0} h^{z_0} (\log h)^{p_0}) + \mathrm{h.o.t.} \\
    &= \left( I_{\tf[1]}(P,z_0)a_{z_0,p_0} \right) s^{h_0} (\log h)^{p_0} + \mathrm{h.o.t.}
\end{align*}
where $I_{\tf[1]}(P,z)$ is the indicial family, given by
\begin{align*}
    I_{\tf[1]}(P,z) &= \left( \frac{r}{r+1} z + r\partial_r \right)^2 - \left( \frac{\beta r}{r+1} + 1 \right) \left( \frac{r}{r+1} z + r\partial_r \right) \\
    &= \left( \frac{1}{1+r'} z - r' \partial_{r'} \right)^2 - \left( \frac{\beta}{1+r'} + 1 \right) \left( \frac{1}{1+r'} z - r' \partial_{r'} \right)
\end{align*}
As a differential operator in $r$, $I_{\tf[1]}(P,z)$ has indicial roots $0, 1$. As a differential operator in $r'$, it has indicial roots $z, z - (\beta+1)$. The expansion of $R_4$ contains the orders $(k+1)(\beta+1)$ for $k \in \N_0$, so the indicial roots of this operator are of the form $(k+1)(\beta+1)$ and $k (\beta+1)$. 

Expand $R_4$ at $\tf[1]$ and solve $P Q_5 = R_4$ termwise using the indicial family of $P$. Each termwise equation takes the form $I_{\tf[1]}(P,z) q = f$ and has a general solution that is polyhomogeneous at $r' = 0$ (i.e., at $\bf{0}$) and $r = 0$ (i.e., at $\bf{1}$). From Theorem \ref{theorem_laplacian-parametrix-with-compact-remainder-vanishing-at-b0-phg-index-family}, $f$ has leading order $1$ at $\bf{0}$ and $\bf{1}$. As $I_{\tf[1]}(P,z)$ has indicial root $0$ at $r = 0$ and indicial root $k(\beta+1)$ (for $k \in \N_0$) at $r' = 0$, a general solution $q$ has a $0$th-order term in its expansion at $r = 0$ and possibly also at $r' = 0$. We impose boundary conditions at $r = 0$ and $r' = 0$ to eliminate this possibility. All other indicial roots of $I_{\tf[1]}(P,z)$ are at least $1$, so it follows that the solution has leading order 1 as $r \to 0$ and $r' \to 0$. This is extremely important, as it means we have not undone our hard work from before. Sum the termwise solutions using Borel's lemma to produce an operator $Q_5 \in \Psi_{\gas}^{-\infty, \check{\mathcal{R}}_4'}$, where $\check{\mathcal{R}}_4'$ is as in Theorem \ref{theorem_laplacian-parametrix-with-compact-remainder-vanishing-at-b0-phg-index-family}, although possibly with increased logarithmic exponents. By construction, $PQ_5 - R_4$ is vanishing to infinite order at $\tf[1]$. For more details on this technique, see \cite{mazzeo_differential-edge-operators-I_1991}.

Now solve away the error at $\bf{1}$. Work in the chart $(\epsilon, s, u, \rho', y_2)$, as depicted in Figure \ref{fig_double-space-upper-region-coords}. After a change of coordinates, the indicial operator and indicial family at this face are
\begin{align*}
    I_{\bf{1}}(P) &= (s \partial_s)^2 - s \partial_s \\
    I_{\bf{1}}(P;z) &= z^2 - z
\end{align*}
Acting termwise on polyhomogeneous expansions at $\bf{1}$, the indicial operator is invertible except at terms corresponding to the indicial roots of $I_{\bf{1}}(P)$, i.e. terms of index $(z,p)$ for $z = 0, 1$. Nonetheless, given a polyhomogeneous expansion $g$ with index set $E$ at $\bf{1}$, there exists polyhomogeneous $f$ with index set $E'$, which differs from $E$ by increasing the logarithmic exponents corresponding to the indicial roots $0, 1$, such that $Pf = g$. Solve the equation $PQ_6 = R_4 - PQ_5$ for $Q_6$, modulo polyhomogeneous functions vanishing to infinite order at $\bf{1}$, once again by summing termwise solutions using Borel's lemma. We produce $Q_6$ that has index family identical to $\check{\mathcal{R}}_4'$ except that it, like $R_4 - PQ_5$, has index set $\emptyset$ at $\tf[1]$. This is because the indicial operator is independent of $\epsilon$, the boundary-defining function of $\tf[1]$. Then $P(Q_5 + Q_6) - R_4$ vanishes to infinite order at both $\bf{1}$ and $\tf[1]$.

We summarize our results as follows:

\begin{proposition}
\label{proposition_laplacian-parametrix-with-compact-remainder-vanishing-at-side-faces-phg-index-family}

Suppose the hypotheses of Theorem \ref{theorem_laplacian-parametrix-with-compact-remainder-vanishing-at-b0-phg-index-family} hold. Then there exist operators $Q_0 \in \Psi^{-2}_{\gas}(M), Q^{(6)} \in \Psi^{-\infty,\mathcal{E}'}_{\gas}(M)$, such that
\begin{align*}
    P\left(Q_0 + Q^{(6)}\right) = I - R_6
\end{align*}
where $\mathcal{E}'$ (possibly different from the $\mathcal{E}'$ from previous results) agrees with the index family $\mathcal{E}$ of Corollary \ref{corollary_laplacian-parametrix-with-vanishing-remainder-at-tff}, except some logarithmic exponents are increased, and $R_6 \in \Psi^{-\infty,\mathcal{R}'_6}_{\gas}(M)$, where $\mathcal{R}'_6$ is identical to the index family $\check{\mathcal{R}}_4'$ of Theorem \ref{theorem_laplacian-parametrix-with-compact-remainder-vanishing-at-b0-phg-index-family}, except that its index sets at $\bf{1}$ and $\tf[1]$ are $\emptyset$ and that possibly some logarithmic exponents are increased.
\end{proposition}
  \section{A composition formula}

Also as a prelude to a Neumann series argument, we must prove a composition formula for $\Psi^{*,*}_{\gas}$. To do this, we construct a triple space $\trsp{M}{\gas}$ by blowing up $[0,1) \times M^3$. Denote the former by $\anysp{M}{0123}$ and the latter by $\anysp{M}{0,1,2,3}$, where commas indicate a (un-blown-up) Cartesian product. Introduce coordinates $(\epsilon, x_1, x_2, x_3, y_1, y_2, y_3)$ on $[0,1) \times M^3$, where $\epsilon$ parametrizes $[0,1)$ and $(x_i, y_i)$ is a coordinate system for the $i$th copy of $M$. We define a family of submanifolds of $[0,1) \times M^3$ using Table \ref{table_triple-space-submanifolds}. Each column is a linear function of the aforementioned coordinates and each row represents a submanifold. Each submanifold is defined as the vanishing set of the functions in whose corresponding column a $0$ is marked.

\begin{table}[ht]
    \centering
    \begin{tabular}{|l|l|l|l|l|l|l|l|l|}
    \hline
              & $x_1$ & $x_2$ & $x_3$ & $y_1-y_2$ & $y_1-y_3$ & $y_2-y_3$ & $\epsilon$ & Front face \\ \hline
    $W^{123}$ & 0     & 0     & 0     & 0         & 0         & 0         & 0          & $\mathrm{tff}_{123}$  \\ \hline
    $V^{12}$  & 0     & 0     & 0     & 0         &           &           & 0          & $\mathrm{v}_{12}$             \\ \hline
    $V^{13}$  & 0     & 0     & 0     &           & 0         &           & 0          & $\mathrm{v}_{13}$             \\ \hline
    $V^{23}$  & 0     & 0     & 0     &           &           & 0         & 0          & $\mathrm{v}_{23}$             \\ \hline
    $W^{12}$  & 0     & 0     &       & 0         &           &           & 0          & $\mathrm{tff}_{12}$   \\ \hline
    $W^{13}$  & 0     &       & 0     &           & 0         &           & 0          & $\mathrm{tff}_{13}$   \\ \hline
    $W^{23}$  &       & 0     & 0     &           &           & 0         & 0          & $\mathrm{tff}_{23}$   \\ \hline
    $C^{123}$ & 0     & 0     & 0     &           &           &           & 0          & $\mathrm{btff}_{123}$ \\ \hline
    $C^{12}$  & 0     & 0     &       &           &           &           & 0          & $\mathrm{btff}_{12}$  \\ \hline
    $C^{13}$  & 0     &       & 0     &           &           &           & 0          & $\mathrm{btff}_{13}$  \\ \hline
    $C^{23}$  &       & 0     & 0     &           &           &           & 0          & $\mathrm{btff}_{23}$  \\ \hline
    $C^1$     & 0     &       &       &           &           &           & 0          & $\mathrm{tf}_1$      \\ \hline
    $C^2$     &       & 0     &       &           &           &           & 0          & $\mathrm{tf}_2$      \\ \hline
    $C^3$     &       &       & 0     &           &           &           & 0          & $\mathrm{tf}_3$      \\ \hline
    $S^{123}$ & 0     & 0     & 0     & 0         & 0         & 0         &            & $\mathrm{zff}_{123}$  \\ \hline
    $U^{12}$  & 0     & 0     & 0     & 0         &           &           &            & $\mathrm{u}_{12}$             \\ \hline
    $U^{13}$  & 0     & 0     & 0     &           & 0         &           &            & $\mathrm{u}_{13}$             \\ \hline
    $U^{23}$  & 0     & 0     & 0     &           &           & 0         &            & $\mathrm{u}_{23}$             \\ \hline
    $S^{12}$  & 0     & 0     &       & 0         &           &           &            & $\mathrm{zff}_{12}$   \\ \hline
    $S^{13}$  & 0     &       & 0     &           & 0         &           &            & $\mathrm{zff}_{13}$   \\ \hline
    $S^{23}$  &       & 0     & 0     &           &           & 0         &            & $\mathrm{zff}_{23}$   \\ \hline
    $B^{123}$ & 0     & 0     & 0     &           &           &           &            & $\mathrm{bff}_{123}$  \\ \hline
    $B^{12}$  & 0     & 0     &       &           &           &           &            & $\mathrm{bff}_{12}$   \\ \hline
    $B^{13}$  & 0     &       & 0     &           &           &           &            & $\mathrm{bff}_{13}$   \\ \hline
    $B^{23}$  &       & 0     & 0     &           &           &           &            & $\mathrm{bff}_{23}$   \\ \hline
    \end{tabular}
    \caption{Submanifolds of $[0,1) \times M^3$ that are blown up in the construction of $\trsp{M}{\gas}$. The last column gives the name of the front face corresponding to the given submanifold.}
    \label{table_triple-space-submanifolds}
\end{table}

Note that the $W$-, $V$-, and $C$-type submanifolds lie in the boundary $\{\epsilon = 0 \}$, which we call the `parameter boundary'. Note also each of these submanifolds is a p-submanifold of $M^3$ and that any chain of inclusions among them is a chain of p-submanifolds \cite{melrose_differential-analysis-manifold-corners_1996}.

We shall construct $\trsp{M}{\gas}$ in nine stages. Each stage entails blowing up one or more submanifolds whose blow-ups mutually commute. First construct $T_1 = [[0,1) \times M^3; W^{123}]$. The submanifolds $V^{ij}$ intersect only in $W^{123}$ and, moreover, their lifts to $T_1$ are disjoint. Now construct $T_2 = \left[T_1; V^{12}; V^{13}; V^{23} \right]$, where we abuse notation by writing $V^{ij}$ to refer to its own lift. Lifted to $T_2$, the submanifolds $W^{ij}$ and $C^{123}$ are pairwise disjoint. Therefore, we can construct the next space: $T_3 = \left[ T_2; W^{12}; W^{13}; W^{23}; C^{123} \right]$. We can similarly form $T_4 = \left[ T_3; C^{12}; C^{13}; C^{23} \right]$ and $T_5 = \left[ T_4; C^{1}; C^{2}; C^{3} \right]$. This concludes the blow-ups performed within the parameter boundary. The remaining blow-ups are as follows: $T_6 = \left[ T_5; S^{123} \right]$, $T_7 = \left[ T_6; U^{12}; U^{13}; U^{23} \right]$, $T_8 = \left[ T_7; S^{12}; S^{13}; S^{23}; B^{123} \right]$, and $T_9 = \left[ T_8; B^{12}; B^{13}; B^{23} \right]$, which we denote, as previously mentioned, by $\trsp{M}{\gas}$ or $\anysp{M}{0123}$.

\begin{figure}
    \centering
    \begin{adjustbox}{max width=0.95\linewidth}
    \begin{tikzcd}
	&&&&&&& {W^{123}} &&&& \\
	\\
	&&&& {S^{123}} && {V^{12}} & {V^{13}} & {V^{23}} \\
	\\
	&&& {U^{12}} & {U^{13}} & {U^{23}} & {\textcolor{red}{W^{12}}} & {W^{13}} & {W^{23}} && {C^{123}} \\
	\\
	& {B^{123}} && {\textcolor{red}{S^{12}}} & {S^{13}} & {S^{23}} &&&& {\textcolor{red}{C^{12}}} & {C^{13}} & {C^{23}} \\
	\\
	{\textcolor{red}{B^{12}}} & {B^{13}} & {B^{23}} &&&&&&& {C^1} & {C^2} & {C^3}
	\arrow[no head, from=1-8, to=3-5]
	\arrow[no head, from=1-8, to=3-7]
	\arrow[no head, from=1-8, to=3-8]
	\arrow[no head, from=1-8, to=3-9]
	\arrow[no head, from=3-5, to=5-4]
	\arrow[no head, from=3-5, to=5-5]
	\arrow[no head, from=3-5, to=5-6]
	\arrow[no head, from=3-7, to=5-4]
	\arrow[no head, from=3-7, to=5-7]
	\arrow[no head, from=3-7, to=5-11]
	\arrow[no head, from=3-8, to=5-5]
	\arrow[no head, from=3-8, to=5-8]
	\arrow[no head, from=3-8, to=5-11]
	\arrow[no head, from=3-9, to=5-6]
	\arrow[no head, from=3-9, to=5-9]
	\arrow[no head, from=3-9, to=5-11]
	\arrow[no head, from=5-4, to=7-2]
	\arrow[no head, from=5-4, to=7-4]
	\arrow[no head, from=5-5, to=7-2]
	\arrow[no head, from=5-5, to=7-5]
	\arrow[no head, from=5-6, to=7-2]
	\arrow[no head, from=5-6, to=7-6]
	\arrow[no head, from=5-7, to=7-4]
	\arrow[no head, from=5-7, to=7-10]
	\arrow[no head, from=5-8, to=7-5]
	\arrow[no head, from=5-8, to=7-11]
	\arrow[no head, from=5-9, to=7-6]
	\arrow[no head, from=5-9, to=7-12]
	\arrow[no head, from=5-11, to=7-10]
	\arrow[no head, from=5-11, to=7-11]
	\arrow[no head, from=5-11, to=7-12]
	\arrow[no head, from=7-2, to=5-11]
	\arrow[no head, from=7-2, to=9-1]
	\arrow[no head, from=7-2, to=9-2]
	\arrow[no head, from=7-2, to=9-3]
	\arrow[no head, from=7-4, to=9-1]
	\arrow[no head, from=7-5, to=9-2]
	\arrow[no head, from=7-6, to=9-3]
	\arrow[no head, from=7-10, to=9-1]
	\arrow[no head, from=7-10, to=9-10]
	\arrow[no head, from=7-10, to=9-11]
	\arrow[no head, from=7-11, to=9-2]
	\arrow[no head, from=7-11, to=9-10]
	\arrow[no head, from=7-11, to=9-12]
	\arrow[no head, from=7-12, to=9-3]
	\arrow[no head, from=7-12, to=9-11]
	\arrow[no head, from=7-12, to=9-12]
    \end{tikzcd}
    \end{adjustbox}
    \caption{The inclusion poset, ascending from top to bottom, of the submanifolds defined in Table \ref{table_triple-space-submanifolds}.}
    \label{fig_triple-space-submanifolds-poset}
\end{figure}
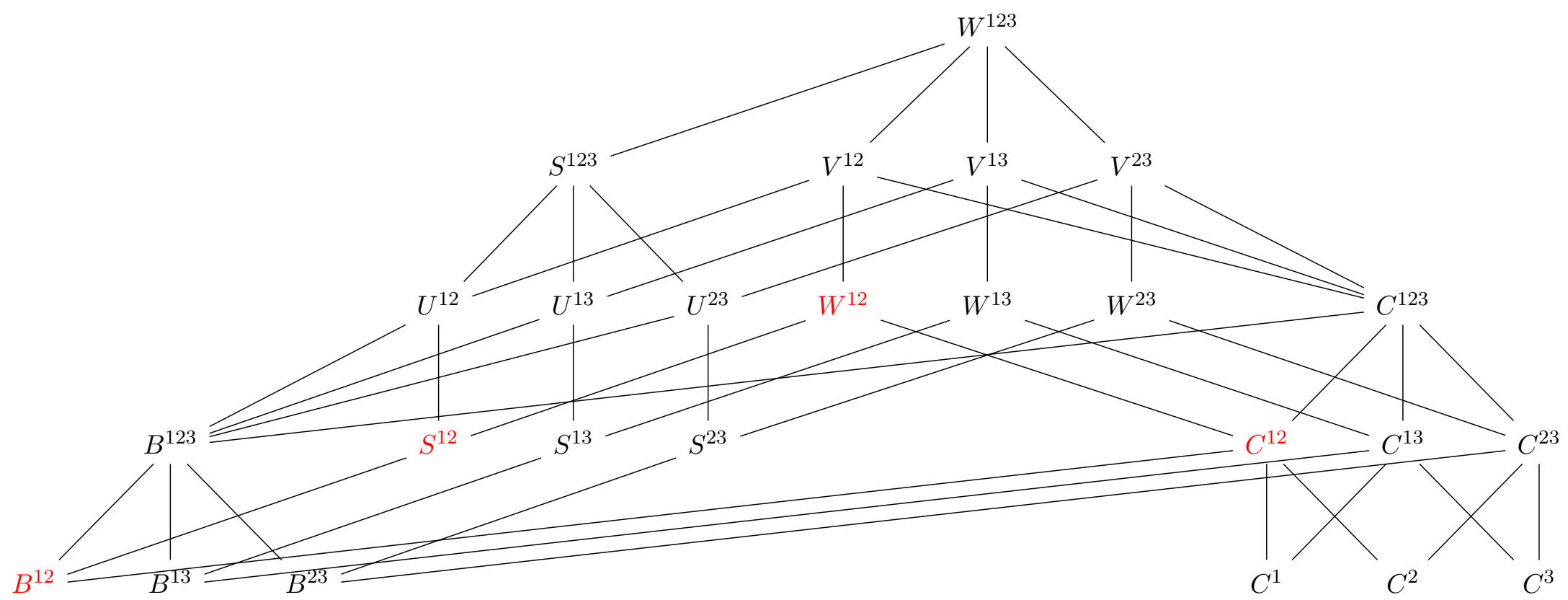

In Figure \ref{fig_triple-space-submanifolds-poset}, the red labels correspond to submanifolds involved in the construction of the double space $\dbsp{M}{\gas}$. To be precise, let $\anysp{M}{0,1,2} := [0,1) \times M^2$, where we think of $M^2$ as consisting of the first two copies of $M$ and let $\beta^{0,1,2,3}_{0,1,2} : \anysp{M}{0,1,2,3} \to \anysp{M}{0,1,2}$ be the map that projects out the third copy of $M$. Then, the red-labeled submanifolds $W^{12}, S^{12}, B^{12}, C^{12}$ are the lifts, along $\beta^{0,1,2,3}_{0,1,2}$, of the submanifolds of $M^2$ that are blown up in the construction of $\dbsp{M}{\gas}$. (Previously, these submanifolds were themselves denoted by $W^{12}, S^{12}, B^{12}, C^{12}$.)

Analogously to the notation $\anysp{M}{0,1,2}$, let $\anysp{M}{0,1,3} = \anysp{M}{0,2,3} = [0,1) \times M^2$, interpreted as consisting of the first and third or second and third copies of $M$, respectively. Similarly define $\anysp{M}{0ij}$ to be $\dbsp{M}{\gas}$, formed from the $i$th and $j$th copy of $M$. In order for composition to be well-defined in the \ztrtext{} calculus, we require the following lemma.

\begin{lemma}
    For $(i,j) = (1,2), (1,3), (2,3)$, the projection $\beta^{0,1,2,3}_{0,i,j} : \anysp{M}{0,1,2,3} \to \anysp{M}{0,i,j}$ extends to a b-fibration $\beta^{0123}_{0ij} : \anysp{M}{0123} \to \anysp{M}{0ij}$.
\end{lemma}
\begin{proof}
    Existence can be shown elegantly by reordering blow-ups so that $\anysp{M}{0123} = [\anysp{M}{012} \times [0,1) ; \mathcal{Y}]$, where $\mathcal{Y}$ is a list of submanifolds of $\anysp{M}{012}$.

    This can be done algorithmically. Suppose $\mathcal{X} = (X_1, \ldots, X_k), \mathcal{X}' = (X_1', \ldots, X_k')$ are enumerations of the submanifolds in Table \ref{table_triple-space-submanifolds}. Write that $\mathcal{X}, \mathcal{X}'$ are equivalent if the spaces $[[0,1) \times M^3; \mathcal{X}], [[0,1) \times M^3; \mathcal{X}']$ are canonically diffeomorphic. Now suppose $\mathcal{X}, \mathcal{X}'$ are related by transposing two adjacent submanifolds; i.e. there exists $j \in \{1,\ldots,k-1\}$ such that $X_i = X_i'$ for all $i \notin \{j,j+1\}$, $X_j = X'_{j+1}$, and $X_{j+1} = X'_j$. Then $\mathcal{X}, \mathcal{X}'$ are equivalent if $X_j \cap X_{j+1} = X_\ell$ for some $\ell \leq j+1$. 
    
    To show this, consider the cases $\ell < j$ and $\ell \in \{j,j+1\}$. In the first case, blowing up $X_\ell$ causes $X_j, X_{j+1}$ to lift to disjoint submanifolds, and therefore the blow-ups of $X_j, X_{j+1}$ commute. In the case where $\ell = j$ or $\ell = j+1$ (without loss of generality, assume the former), the pair $X_j \subset X_{j+1}$ is a chain of p-submanifolds in the sense of \cite{melrose_differential-analysis-manifold-corners_1996}, and therefore again their blow-ups commute.

    Let $\mathcal{X}_0$ be the enumeration of the submanifolds in Table \ref{table_triple-space-submanifolds} in the order in which they are blown up. Using the equivalence rule, we can show that $\mathcal{X}_0$ is equivalent to an enumeration of the form $\mathcal{X}_* = (W^{12}, C^{12}, S^{12}, B^{12}, C^1, C^2, \ldots)$. The first six submanifolds are identified with the submanifolds blown up in the construction of $\anysp{M}{012} = \dbsp{M}{\gas}$, so $\anysp{M}{0123} = [\anysp{M}{0,1,2,3}; \mathcal{X}_0] \cong [\anysp{M}{0,1,2,3}; \mathcal{X}_*] = [ \anysp{M}{012,3} ; \mathcal{Y}]$, where $\anysp{M}{012,3} := \anysp{M}{012} \times M$ and $\mathcal{Y}$ is the sublist of $\mathcal{X}_*$ that follows $C^2$.
    
    Let $\beta^{0123}_{012,3} : \anysp{M}{0123} \to \anysp{M}{012,3}$ be the map that blows down $\mathcal{Y}$. Let $\beta^{012,3}_{012} : \anysp{M}{012,3} \to \anysp{M}{012}$ be the map that projects out the third copy of $M$. Both maps are b-submersions, and hence so is $\beta^{0123}_{012} = \beta^{012,3}_{012} \circ \beta^{0123}_{012,3}$. In fact, this is a b-fibration, which can be shown by explicitly computing its lifting matrix.
    % \begin{align*}
    % \setcounter{MaxMatrixCols}{29}
    %     e^{0123}_{012} = \begin{pmatrix}
    %          1 & 0 & 0 & 0 & 0 & 0 & 0 & 0 & 0 & 0 & 0 & 0 & 0 & 0 & 0 & 0 & 0 & \
    %         1 & 0 & 0 & 0 & 0 & 0 & 0 & 0 & 0 & 0 & 0 & 0 \\
    %          0 & 1 & 0 & 0 & 0 & 0 & 0 & 0 & 0 & 0 & 0 & 0 & 0 & 0 & 0 & 0 & 0 & \
    %         0 & 0 & 0 & 0 & 0 & 0 & 1 & 0 & 0 & 0 & 1 & 0 \\
    %          0 & 0 & 1 & 0 & 0 & 0 & 0 & 0 & 0 & 0 & 0 & 0 & 0 & 0 & 0 & 0 & 0 & \
    %         0 & 0 & 0 & 0 & 0 & 0 & 0 & 1 & 0 & 0 & 0 & 1 \\
    %          0 & 0 & 0 & 0 & 1 & 1 & 0 & 0 & 1 & 0 & 0 & 0 & 0 & 0 & 0 & 0 & 0 & \
    %         0 & 0 & 0 & 0 & 0 & 0 & 0 & 0 & 0 & 0 & 0 & 0 \\
    %          0 & 0 & 0 & 0 & 0 & 0 & 1 & 1 & 0 & 0 & 0 & 1 & 1 & 0 & 0 & 0 & 0 & \
    %         0 & 0 & 0 & 0 & 0 & 0 & 0 & 0 & 0 & 0 & 0 & 0 \\
    %          0 & 0 & 0 & 0 & 0 & 0 & 0 & 0 & 0 & 0 & 0 & 0 & 0 & 0 & 0 & 0 & 0 & \
    %         0 & 1 & 1 & 0 & 0 & 1 & 0 & 0 & 0 & 0 & 0 & 0 \\
    %          0 & 0 & 0 & 0 & 0 & 0 & 0 & 0 & 0 & 0 & 0 & 0 & 0 & 0 & 0 & 0 & 0 & \
    %         0 & 0 & 0 & 1 & 1 & 0 & 0 & 0 & 1 & 1 & 0 & 0 \\
    %          0 & 0 & 0 & 0 & 0 & 0 & 0 & 0 & 0 & 1 & 0 & 0 & 0 & 1 & 0 & 1 & 0 & \
    %         0 & 0 & 0 & 0 & 0 & 0 & 0 & 0 & 0 & 0 & 0 & 0 \\
    %          0 & 0 & 0 & 0 & 0 & 0 & 0 & 0 & 0 & 0 & 1 & 0 & 0 & 0 & 1 & 0 & 1 & \
    %         0 & 0 & 0 & 0 & 0 & 0 & 0 & 0 & 0 & 0 & 0 & 0 \\
    %         \end{pmatrix}
    % \end{align*}
    By symmetry, an identical argument constructs b-fibrations $\beta^{0123}_{0ij}$ for $(i,j) = (1,3),(2,3)$. 
\end{proof}

This allows us to write the composition of \ztrtext{} pseudodifferential operators as follows. Suppose $A \in \Psi^{\ell, \mathcal{E}}_{\gas}(M)$ and $B \in \Psi^{m, \mathcal{F}}_{\gas}(M)$, with half-density-valued kernels $K_A, K_B$ on $\dbsp{M}{\gas}$. Then, assuming it is well-defined, the composition $C = AB$ has kernel
\begin{align*}
    K_C = \left(\beta^{0123}_{013}\right)_* \left( \left(\beta^{0123}_{012}\right)^* (K_A) \left(\beta^{0123}_{023}\right)^* (K_B) \right)
\end{align*}

\begin{proposition}
\label{proposition_0-transition-calculus-composition-formula}
    Suppose $A \in \Psi^{\ell,\mathcal{E}}_\gas(M)$ and $B \in \Psi^{m,\mathcal{F}}_\gas(M)$. Suppose moreover that $E_{\bf{2}} + F_{\bf{1}} > -1$.
    
    Then the composition $AB$ is well-defined and $AB \in \Psi^{\ell+m,\mathcal{G}}_\gas(M)$, where $\mathcal{G}$ is given by
    \begin{align*}
        G_{\bf{0}}
        &=(E_{\bf{0}}+F_{\bf{0}})
        \mathbin{\overline\cup}(E_{\tf[2]}+F_{\tf[1]}+1)
        \\
        G_{\bf{1}}
        &=E_{\bf{1}}
        \mathbin{\overline\cup}(E_{\zff}+F_{\bf{1}})
        \mathbin{\overline\cup}(E_{\bff}+F_{\bf{1}}+1)
        \\
        G_{\bf{2}}
        &=F_{\bf{2}}
        \mathbin{\overline\cup}(E_{\bf{2}}+F_{\zff})
        \mathbin{\overline\cup}(E_{\bf{2}}+F_{\bff}+1)
        \\
        G_{\tff}
        &=(E_{\tff}+F_{\tff})
        \mathbin{\overline\cup}(E_{\btff}+F_{\btff}+n+1)
        \mathbin{\overline\cup}(E_{\tf[1]}+F_{\tf[2]}+n)
        \\
        G_{\btff}
        &=(E_{\btff}+F_{\btff}+1)
        \mathbin{\overline\cup}(E_{\tff}+F_{\btff})
        \mathbin{\overline\cup}(E_{\btff}+F_{\tff})
        \mathbin{\overline\cup}(E_{\tf[1]}+F_{\tf[2]})
        \\
        G_{\zff}
        &=(E_{\zff}+F_{\zff})
        \mathbin{\overline\cup}(E_{\bf{1}}+F_{\bf{2}}+n)
        \mathbin{\overline\cup}(E_{\bff}+F_{\bff}+n+1)
        \\
        G_{\bff}
        &=(E_{\bf{1}}+F_{\bf{2}})
        \mathbin{\overline\cup}(E_{\bff}+F_{\bff}+1)
        \mathbin{\overline\cup}(E_{\zff}+F_{\bff})
        \mathbin{\overline\cup}(E_{\bff}+F_{\zff})
        \\
        G_{\tf[1]}
        &=(E_{\tf[1]}+F_{\bf{0}})
        \mathbin{\overline\cup}(E_{\tff}+F_{\tf[1]})
        \mathbin{\overline\cup}(E_{\btff}+F_{\tf[1]}+1)
        \\
        G_{\tf[2]}
        &=(E_{\bf{0}}+F_{\tf[2]})
        \mathbin{\overline\cup}(E_{\tf[2]}+F_{\tff})
        \mathbin{\overline\cup}(E_{\tf[2]}+F_{\btff}+1).
    \end{align*}
\end{proposition}

Equipped with the composition formula, we can once again improve the parametrix.

\begin{theorem}
\label{theorem_laplacian-parametrix-with-infinite-order-remainder}

Suppose the hypotheses of Theorem \ref{theorem_laplacian-parametrix-with-compact-remainder-vanishing-at-b0-phg-index-family} hold. Then there exist operators $Q_0 \in \Psi^{-2}_{\gas}(M), Q^{(7)} \in \Psi^{-\infty,\mathcal{E}'}_{\gas}(M)$, such that
\begin{align*}
    P\left(Q_0 + Q^{(7)}\right) = I - R_7
\end{align*}
where $\mathcal{E}'$ (possibly different from previous $\mathcal{E}'$) agrees with the index family $\mathcal{E}$ of Corollary \ref{corollary_laplacian-parametrix-with-vanishing-remainder-at-tff}, except some logarithmic exponents are increased, and $R_7 \in \Psi^{-\infty,\mathcal{R}_7'}_{\gas}(M)$, where $\mathcal{R}'_7$ is an index family whose index set at all boundary hypersurfaces is $\emptyset$, except at $\bf{2}$, where it is $(-1,0) + \N_0$ except possibly with some increased logarithmic indices.

\end{theorem}
\begin{proof}
    Let $Q_0, Q^{(6)}, R_6$ be as in Proposition \ref{proposition_laplacian-parametrix-with-compact-remainder-vanishing-at-side-faces-phg-index-family} so that
    \begin{align*}
        P\left(Q_0 + Q^{(6)}\right) &= I - R_6
    \end{align*}
    The leading conormal orders of $R_6$ are $(1,\infty,-1,1,0,1,0,\infty,-1)$ at the boundary hypersurfaces $\bf{0}, \bf{1}, \bf{2}, \tff, \btff, \zff, \bff, \tf[1], \tf[2]$, respectively. Using the composition formula (Proposition \ref{proposition_0-transition-calculus-composition-formula}), we find that $R_6^j$ has leading orders $(j,\infty,-1,j,j-1,j,j-1,\infty,j-2)$ (relative to the ordering of boundary hypersurfaces given in the statement of Proposition \ref{proposition_0-transition-calculus-composition-formula}). Using Borel's lemma, construct $S \approx \sum_{j=1}^\infty R_6^j$, i.e. $S \in \Psi^{-\infty,\ast}_{\gas}(M)$ such that $S - \sum_{j=1}^k R_6^j$ has leading orders $(k+1,\infty,-1,k+1,k,k+1,k,\infty,k-1)$. Collecting residual terms, write $(Q_0 + Q^{(6)})(I+S) = Q_0 + Q^{(7)}$, where $Q^{(7)} \in \Psi^{-\infty,\ast}_{\gas}(M)$. Then,
    \begin{align*}
        P\left(Q_0 + Q^{(7)}\right) &= (I - R_6)(I+S) = I - R_7
    \end{align*}
    where, by construction of the Neumann series, $R_7 := I - (I - R_0)(I + S)$ is polyhomogeneous with leading orders $(\infty,\infty,-1,\infty,\infty,\infty,\infty,\infty,\infty)$.
\end{proof}

\section{Uniform structure of the resolvent}

With the aid of the stronger parametrix construction in Theorem \ref{theorem_laplacian-parametrix-with-infinite-order-remainder}, we can study the resolvent of the Laplacian. Let $P_\epsilon(\lambda) = x^2 x_\epsilon^{-\alpha} (\Delta_{\epsilon} - \lambda) = P_\epsilon - \lambda x^2 x_\epsilon^{-\alpha}$. Let $P(\lambda)$ be the family $P_\epsilon(\lambda)$ viewed as an element of $\mathrm{Diff}^*_\gas(M)$. (Of course, $P(0) = P$). As before, assume that $\Delta_0 = \Delta_g$ denotes the Laplacian on the domain $x^{2-\alpha} H^2_0(dV_g) \subset L^2(dV_g)$, which is a domain of self-adjointness under the assumption $\alpha > 2/n$ \cite{dHIKM_gas-giants_2024}.

From classical perturbation theory, the Dirichlet spectrum of $\Delta_{\epsilon}$ varies continuously and piecewise smoothly in $\epsilon$. It was also shown in \cite{dHIKM_gas-giants_2024} that, if $\{ \lambda_\epsilon \}_{0 < \epsilon < 1}$ is a continuous family such that $\lambda_\epsilon$ is a Dirichlet eigenvalue of $\Delta_{\epsilon}$, then $\lambda_\epsilon$ limits to an eigenvalue $\lambda_0$ of $\Delta_0$. We show that, if $\lambda_0$ is a simple eigenvalue, then in fact the limiting behavior of $\lambda_\epsilon$ as $\epsilon \to 0$ is polyhomogeneous. To prove this, we study the resolvent $P(\lambda)^{-1}$. (We call it a resolvent, but of course it is actually a multiple of the resolvent of $\Delta_{\epsilon}$.)

\subsection{Away from the spectrum of $\Delta_g$}

For $\epsilon_0 > 0$, define ${}^{\epsilon_0} \sgsp{M}{\gas}$ and ${}^{\epsilon_0} \dbsp{M}{\gas}$ to be the subsets of $\sgsp{M}{\gas}$ and $\dbsp{M}{\gas}$, respectively, where $\epsilon < \epsilon_0$. Define ${}^{\epsilon_0}\Psi^{\ast,\ast}_{\gas}$ to consist of the restrictions of kernels in $\Psi^{\ast,\ast}_{\gas}$ to ${}^{\epsilon_0} \dbsp{M}{\gas}$. The following result characterizes $P(\lambda)^{-1}$ away from the spectrum of $\Delta_0$.

\begin{theorem}
\label{theorem_resolvent-holomorphic-away-from-spectrum}

    Suppose that $\lambda_0 \in \C \setminus \mathrm{spec} (\Delta_g)$. Then there exists $\epsilon_0 > 0$ such that the resolvent $P(\lambda)^{-1}$ is a holomorphic function of $\lambda$ on a neighborhood of $\lambda_0$, taking values in $^{\epsilon_0} \Psi^{\ast,\ast}_{\gas}$.
\end{theorem}

A little care should be taken in defining a holomorphic function valued in a space of polyhomogeneous distributions. Many of the complications are avoided if the values of the function share a common index family, which is always the case for us. If we write that $A(\lambda) \in \Psi^{m,\ast}_{\gas}(M)$ is holomorphic, we mean that there exists an index family $\mathcal{E}$ such that $A(\lambda)$ is a holomorphic map into $\Psi^{m,\mathcal{E}}_{\gas}(M)$, in a sense that is stated carefully in, for example, \cite[8.4]{mazzeo-melrose_analytic-surgery_1995}.

The techniques we use are very similar to those used in \cite[5.2]{mazzeo-melrose_analytic-surgery_1995}. The key step is a parametrix construction for $P(\lambda)$. This operator has an identical principal symbol to $P$ as well as identical normal operators at all front faces (except for $\bf{0}$), so our parametrix construction for $P$ works almost equally well for $P(\lambda)$. However, it is important to perform this parametrix construction so that the parametrix and remainder are holomorphic, at least in a small neighborhood of $\lambda_0$. To prevent issues of convergence, we avoid Borel summations. This inhibits us from proving an exact analogue of Theorem \ref{theorem_laplacian-parametrix-with-infinite-order-remainder}. Fortunately, a weaker result suffices:

\begin{lemma}
\label{lemma_laplacian-minus-lambda-parametrix-with-finite-order-remainder}

    Suppose that $\lambda_0 \in \C \setminus \mathrm{spec} (\Delta_g)$. Let $N > 0$.

    There exist holomorphic $Q_0(\lambda), Q_0'(\lambda) \in \Psi^{-2,\ast}_{\gas}(M), {}^{(N)}Q(\lambda), {}^{(N)}Q'(\lambda) \in \Psi^{-\infty,\ast}_{\gas}(M)$ defined for $\lambda$ in a neighborhood of $\lambda_0$ such that
    \begin{align}
        P(\lambda) \left( Q_0(\lambda) + {}^{(N)}Q(\lambda) \right) &= I - {}^{(N)}R(\lambda) \label{eqn_laplacian-minus-lambda-parametrix-right} \\
        \left( Q_0'(\lambda) + {}^{(N)}Q'(\lambda) \right)P(\lambda) &= I - {}^{(N)}R'(\lambda) \label{eqn_laplacian-minus-lambda-parametrix-left}
    \end{align}
    where ${}^{(N)}R(\lambda)$ (respectively ${}^{(N)}R'(\lambda)$) has leading order greater than $N$ at all boundary hypersurfaces of $\dbsp{M}{\gas}$ except for $\bf{2}$ (respectively $\bf{1}$), where it has index set $(-1,0) + \N_0$ except possibly with some increased logarithmic exponents.
\end{lemma}

\begin{proof}
    It suffices to construct the right parametrix $Q_0(\lambda) + {}^{(N)}Q(\lambda)$; the left parametrix construction follows by taking adjoints.
    
    Write $P(\lambda) = P(\lambda_0) - (\lambda - \lambda_0) x^2 x_\epsilon^{-\alpha}$. By the usual method (inverting the symbol of $P(\lambda_0)$ and taking a Neumann series), obtain $Q_0(\lambda_0) \in \Psi^{-2,\ast}_{\gas}(M)$ such that $P(\lambda_0) Q_0(\lambda_0) = I - R_0(\lambda_0)$, where $R_0(\lambda_0) \in \Psi^{-\infty, \ast}_{\gas}$. The parametrix and remainder term are in the large calculus because $x_\epsilon^{-\alpha}$ contributes a non-integer leading order at the faces $\tff, \btff$, and $\tf[1]$. Then $P(\lambda) Q_0(\lambda_0) = I - R_0(\lambda_0) - (\lambda-\lambda_0) x^2 x_\epsilon^{-\alpha} Q_0(\lambda_0)$. The term $x^2 x_\epsilon^{-\alpha} Q_0(\lambda_0)$ lies in $\Psi^{-2,\ast}_{\gas}$. The Neumann series $I + \sum_{j=1}^\infty (\lambda-\lambda_0)^j (x^2 x_\epsilon^{-\alpha} Q_0(\lambda_0))^j$ is convergent (with respect to seminorms induced by conormal estimates at the diagonal and polyhomogeneous estimates at the boundary) if $|\lambda - \lambda_0|$ is sufficiently small, and so defines a holomorphic function of $\lambda$ in a neighborhood of $\lambda_0 \in \C$. Right-multiply $Q_0(\lambda_0)$ by this Neumann series to obtain $Q_0(\lambda) \in \Psi^{-2,\ast}_{\gas}$, holomorphic in a neighborhood of $\lambda_0$. With careful bookkeeping, the index families of our parametrices and remainders could be determined, but we will not do this in the interest of brevity.
    
    Observe that the term $\lambda x^2 x_\epsilon^{-\alpha}$ vanishes at the front faces $\tff, \btff, \zff$, and $\bff$. Therefore, the model problems for $P(\lambda)$ at these faces are identical to those of $P = P(0)$. The normal operator of $P(\lambda)$ at $\bf{0}$ is precisely $P_0 - \lambda x^{2-\alpha} = x^{2-\alpha}(\Delta_g - \lambda)$, which is an invertible operator $x^{2-\alpha} H^2_0 \to x^{2-\alpha} L^2$ for $\lambda$ near $\lambda_0$.

    Therefore, we can replicate the steps of the proof of Theorem \ref{theorem_laplacian-parametrix-with-compact-remainder-vanishing-at-b0-phg-index-family} to obtain $Q_1(\lambda) \in \Psi^{-\infty,\ast}_{\gas}(M)$ such that
    \begin{align*}
        P(\lambda) \left( Q_0(\lambda) + Q_1(\lambda) \right) &= I - R_1(\lambda)
    \end{align*}
    All steps can be carried out holomorphically in $\lambda$, so $Q_1(\lambda)$ is holomorphic near $\lambda_0$. The index families of $Q_1(\lambda), R_1(\lambda)$ are larger than those of the parametrix and remainder from \ref{theorem_laplacian-parametrix-with-compact-remainder-vanishing-at-b0-phg-index-family}, again due to the $x^2 x_\epsilon^{-\alpha}$ term, but this is not a serious problem. The leading terms of $R_1(\lambda)$ at $\bf{0}, \tff, \btff, \zff$, and $\bff$ may be lower than those of the old remainder, but the difference is strictly less than $1$; this is essentially attributed to the fact that $2 - \alpha > 0$.

    Then we proceed to prove an analogue of Proposition \ref{proposition_laplacian-parametrix-with-compact-remainder-vanishing-at-side-faces-phg-index-family}. Because solving away the error to infinite order at the side faces $\bf{1}, \tf[1]$ requires a Borel sum, we only solve away to a large finite order by taking partial Borel sums. Produce $Q^{(2)}(\lambda) \in \Psi^{-\infty,\ast}_{\gas}(M)$ such that
    \begin{align*}
        P(\lambda) \left( Q_0(\lambda) + Q^{(2)}(\lambda) \right) &= I - R_2(\lambda)
    \end{align*}
    where $R_2(\lambda)$ has a large positive leading order $M = M(N)$ at $\bf{1}$ and $\tf[1]$. The value of $M$ needs to be large enough for the partial Neumann sum argument that follows.

    Now replicate the steps from the proof of Theorem \ref{theorem_laplacian-parametrix-with-infinite-order-remainder}. Assuming $M$ was chosen sufficiently large, there exists $L \in \N$ such that $R_2(\lambda)^L$ has leading order at least $N$ at all faces besides $\bf{2}$. Let $S = \sum_{j=1}^{L-1} R_2(\lambda)^j$, so that
    \begin{align*}
        P(\lambda) \left( Q_0(\lambda) + Q^{(2)}(\lambda) \right) (I + S) &= I - R_2(\lambda)^L
    \end{align*}
    Collecting residual terms, obtain $Q^{(3)}(\lambda) \in \Psi^{-\infty,\ast}_{\gas}(M)$ such that
    \begin{align*}
        P(\lambda) \left(Q_0(\lambda) + Q^{(3)(\lambda)}\right) &= I - {}^{(N)}R(\lambda)
    \end{align*}
    where ${}^{(N)}R(\lambda) = R_2(\lambda)^L$ has the required leading asymptotics.
\end{proof}

Another ingredient in the proof of Theorem \ref{theorem_resolvent-holomorphic-away-from-spectrum} is a semi-ideal property. Recall that the remainder term in Theorem \ref{theorem_laplacian-parametrix-with-infinite-order-remainder} is rapidly vanishing at all boundary hypersurfaces except $\bf{2}$. In particular, this means that it descends to a polyhomogeneous distribution on $[0,1) \times M^2$. Such kernels are called, in the fashion of \cite{melrose_aps-index-theorem_1993}, `very residual'. Very residual kernels enjoy a semi-ideal property: given an operator $B$ satisfying mild assumptions and very residual operators $A, C$, then, assuming the composition $ABC$ is well-defined, it is itself a very residual operator.

However, as Lemma \ref{lemma_laplacian-minus-lambda-parametrix-with-finite-order-remainder} produces a holomorphic remainder term that is vanishing only to finite order at faces other than $\bf{2}$ (or $\bf{1}$), we need the following modification of the semi-ideal property.

\begin{lemma}
\label{lemma_semi-ideal}

    Let $J \in \N$, $K,L \in \R$. Let $F_1, F_2$ be index sets and let $f_1, f_2$ be real-valued constants such that $F_i > f_i$. Let $f_0 \in \R$ be another constant. Then there exists $N \in \R$ such that the following is true.

    Suppose $A \in \Psi^{-\infty,\mathcal{F}_1}_{\gas}(M), C \in \Psi^{-\infty,\mathcal{F}_2}_{\gas}(M)$, where $\mathcal{F}_1$ (resp. $\mathcal{F}_2$) is an index family on $\dbsp{M}{\gas}$ such that $\mathcal{F}_1(\mathrm{h})$ (resp. $\mathcal{F}_2(\mathrm{h})$) has leading order strictly greater than $N$ for all boundary hypersurfaces $\mathrm{h}$ except for $\bf{1}$ (resp. $\bf{2}$) and $\mathcal{F}_i(\bf{i}) = F_i$ for $i = 1, 2$.
    
    Suppose too $B \in \mathcal{A}^{f_0}([0,1)_\epsilon, \mathcal{L}(L^2(M, \Omega^{1/2})))$, meaning $B \in C^\infty((0,1)_\epsilon, \mathcal{L}(L^2(M, \Omega^{1/2})))$ and $B$ satisfies conormal estimates of order $f_0$ at $0$: there exist $\epsilon_0 > 0$ and $C_j > 0$ for each $j \in \N_0$ such that
    \begin{align*}
        \| (\epsilon \partial_\epsilon)^j B_\epsilon \|_{\mathcal{L}(L^2(M))} \leq C_j \epsilon^{f_0}
    \end{align*}
    for $\epsilon \in (0,\epsilon_0)$.

    Then the scalar kernel $\mathcal{K}_T$ of $T:=ABC$, normalized with respect to a standard half-density $(d\epsilon dx_1 dx_2 dy_1 dy_2)^{1/2}$ on $[0,1) \times M^2$ (not $\dbsp{M}{\gas}$), obeys the following estimates of polyhomogeneous type. For $\ell \in \R$, define
    \begin{align*}
        B_{F_i,\ell} (x_i \partial_{x_i}) = \prod_{\substack{(z,p) \in F_i \\ \mathrm{Re} z \leq \ell}} (x_i \partial_{x_i} - z)
    \end{align*}
    Then, for $\ell < L$, $i \in \{ 1, 2 \}$, and $j = 3 - i$,
    \begin{align}
        \left| B_{F_i,\ell}(x_i \partial_{x_i}) \mathcal{K}_T \right| \lesssim \epsilon^K x_i^\ell x_j^{f_j} \label{eqn_semi-ideal-phg-estimate}
    \end{align}
    and this estimate is stable under the action of b-differential operators on $[0,1) \times M^2$, of order up to $J$.
\end{lemma}

\begin{proof}
Let $A, B, C$ satisfy the assumptions above, for $N$ to be determined. Let their kernels be $\mathcal{K}_A, \mathcal{K}_B, \mathcal{K}_C$, normalized with respect to $(d\epsilon dx_1 dx_2 dy_1 dy_2)^{1/2}$, not $\sigma$ as we have previously done. Then, $\mathcal{K}_T$, viewed as a distribution on $[0,1) \times M^2$ is given by
\begin{align}
    \int \mathcal{K}_A(\epsilon,x_1,y_1,x_1',y_1') \mathcal{K}_B(\epsilon,x_1',y_1',x_2',y_2') \mathcal{K}_C(\epsilon,x_2',y_2',x_2,y_2) \,dx_1' dy_1' dx_2' dy_2' \label{eqn_semi-ideal-triple-product}
\end{align}
where the middle kernel is distributional.

We shall show that $B_{F_i,\ell}(x_i \partial_{x_i}) \mathcal{K}_T$ for $\ell < L$ satisfies conormal estimates. We describe only the $i = 1$ case here. Applying $B_{F_1,\ell}(x_1 \partial_{x_1})$ to $\mathcal{K}_T$, we see that it acts only on $\mathcal{K}_A$, so we shall first establish control over $B_{F_1,\ell}(x_1 \partial_{x_1}) \mathcal{K}_A$. Lift $B_{F_1,\ell} \mathcal{K}_A$ from $[0,1) \times M^2$ (parametrized by $\epsilon, x_1, y_1, x_1', y_1'$) to $\dbsp{M}{\gas}$. Let $U$ be a collar neighborhood of $\bf{1}$ in $\dbsp{M}{\gas}$ so that $U \cong \bf{1} \times [0,1)$. We can choose a family of coordinate neighborhoods on $U$, wherein $x_1 \partial_{x_1}$ lifts in each chart to $\varrho_{\bf{1}} \partial_{\varrho_{\bf{1}}}$, where $\varrho_{\bf{1}}$ is a local boundary-defining function of $\bf{1}$. In fact, we can always take $\varrho_{\bf{1}} = s = x_1/x_1'$ or $\varrho_{\bf{1}} = r = x_1/\epsilon$. (See Figures \ref{fig_double-space-upper-region-coords}, \ref{fig_double-space-middle-region-coords}, and \ref{fig_double-space-lower-region-coords}.) Because $\mathcal{K}_A$ is polyhomogeneous on $\dbsp{M}{\gas}$ with index set $F_1$ at $\bf{1}$, a bound of the following form holds on $U$:
\begin{align}
    \left| B_{F_1,\ell}(\varrho_1 \partial_{\varrho_1}) \mathcal{K}_A \right| \lesssim \varrho_{\bf{1}}^\ell \phi
\end{align}
where $\phi$ is polyhomogeneous on $\bf{1}$. The leading order of $\phi$ at each boundary face of $\bf{1}$ is $N$ plus a constant, due to the normalization of $A$. This bound is stable under b-derivatives on $\dbsp{M}{\gas}$.

Now blow down along $\beta^{012}_{0,1,2} : \dbsp{M}{\gas} \to [0,1) \times M^2$. On $\bf{1} \subset [0,1) \times M^2$, $\phi$ is no longer polyhomogeneous, but it still has at most polynomial growth at all boundary faces, albeit its leading orders at $\bf{0} \cap \bf{1}$ and $\bf{2} \cap \bf{1}$ may be worse than on $\dbsp{M}{\gas}$, although they can be bounded by increasing affine functions of $N$. Moreover, the growth estimates are stable under finitely many b-derivatives. Using these properties of $\phi$ and the identity $\varrho_{\bf{1}} = x_1/x_1'$ or $\varrho_{\bf{1}} = x_1/\epsilon$, the following bound holds on $\beta^{012}_{0,1,2}(U)$ for $\ell < L$:
\begin{align}
    \left| B_{F_1,\ell}(x_1 \partial_{x_1}) \mathcal{K}_A \right| \lesssim x_1^\ell (x_1')^{K'} \epsilon^{K'} \label{eqn_semi-ideal-phg-estimate-left}
\end{align}
where $K'$ is a affine function of $L$ and $N$, decreasing in $L$ and increasing in $N$.

To show stability under b-derivatives, suppose $D \in \mathrm{Diff}^J_{\mathrm{b}}([0,1) \times M^2)$. Then $D$ lifts to b-differential operator on $\dbsp{M}{\gas}$, with some coefficients possibly singular at the boundary, although of order no lower than $-J$. By reducing $K'$ if necessary, we can arrange for $\left| D B_{F_1,\ell}(x_1 \partial_{x_1}) \mathcal{K}_A \right| \lesssim x_1^\ell (x_1')^{K'} \epsilon^{K'}$ to hold, where now $K'$ is an affine function that is increasing in $N$ and decreasing in $J, L$.

% By choosing $N$ to be sufficiently large, we can ensure $K' \geq K$.

Now we handle the behavior of $\mathcal{K}_A$ away from $\bf{1}$. In any neighborhood $V$ of the boundary in $\dbsp{M}{\gas} \setminus U$, we can select a coordinate chart wherein the boundary-defining functions of boundary hypersurfaces that intersect $V$ are given by quotients of $\epsilon, x_1, x_1'$, and $|Y| = |y_1-y_1'|$. In $V$, $\mathcal{K}_A$ obeys a bound
\begin{align*}
    |\mathcal{K}_A| \lesssim \left( \prod_i \varrho_{\mathrm{h}_i} \right)^{N'}
\end{align*}
where $N'$ is an increasing affine function of $N$ and $\{ \mathrm{h}_i \}$ are the boundary hypersurfaces of $\dbsp{M}{\gas}$ which intersect $V$. Of course, $\bf{1}$ is not one of these. This bound is stable under b-derivatives on $\dbsp{M}{\gas}$. We can of course weaken this bound: for $N'_i \leq N'$, we have
\begin{align*}
    |\mathcal{K}_A| \lesssim \prod_i \left(\varrho_{\mathrm{h}_i}\right)^{N'_i}
\end{align*}
By choosing $N_i'$ so as to judiciously cancel denominators, we obtain a bound $|\mathcal{K}_A| \lesssim x_1^{K'} (x_1')^{K'} \epsilon^{K'}$, where $K'$ is an increasing affine function of $N$. By selecting $N$ to be sufficiently large, we ensure that $K' \geq L$ and so the bound \eqref{eqn_semi-ideal-phg-estimate-left} also holds on $V$ for $\ell < L$.
By increasing $N$ if necessary, we can arrange for this bound to be stable under $\mathrm{Diff}^J_{\mathrm{b}}([0,1) \times M^2)$. Hence the bound \eqref{eqn_semi-ideal-phg-estimate-left} holds globally on $\dbsp{M}{\gas}$.

Having established control over $B_{F_1,\ell}(x_1 \partial_{x_1}) \mathcal{K}_A$, it remains to control $\mathcal{K}_B$ and $\mathcal{K}_C$. By taking $N$ to be sufficiently large and using the same techniques we used for $\mathcal{K}_A$, one can show that the kernel $\mathcal{K}_C$ satisfies an estimate $|\mathcal{K}_{C}(\epsilon,x_2',y_2',x_2,y_2)| \lesssim \epsilon^{K'} (x_2')^{K'} x_2^{f_2}$, stable under b-derivatives of order up to $J$. Apply the Cauchy-Schwarz inequality to the integral in $dx_1' dy_1'$ in \eqref{eqn_semi-ideal-triple-product} to control it by the product of an $L^2$ norm of $\mathcal{K}_A(\epsilon,x_1,y_1,\cdot,\cdot)$ and an $L^2$ norm of $\int \mathcal{K}_B \mathcal{K}_C \,dx_2' dy_2'$. Control the former using the conormal estimates on $\mathcal{K}_A$ and control the latter using the operator norm of $\mathcal{K}_B$ and conormal estimates on $\mathcal{K}_C$. The integrals in $(x_1',y_1')$ and in $(x_2',y_2')$ converge if $K'$ is large enough, which can be guaranteed by choosing $N$ large enough. Although the operator norm of $\mathcal{K}_B$ may blow up as $\epsilon \to 0$, we can compensate by making $N$ (and hence $K'$) larger and making $\mathcal{K}_A, \mathcal{K}_C$ vanish faster in $\epsilon$. If $K'$ is chosen large enough, we can ensure $\mathcal{K}_T$ vanishes to order $K$ in $\epsilon$.

Altogether, this proves the bound \eqref{eqn_semi-ideal-phg-estimate} for $i = 1$. This bound is stable under b-differential operators of order up to $J$ on $[0,1) \times M^2$. Stability under $x_1 \partial_{x_1}, \partial_{y_1}$ follows from the stability of the bound \eqref{eqn_semi-ideal-phg-estimate-left}. Stability under $x_2 \partial_{x_2}, \partial_{y_2}$ follows from the stability of the conormal bounds on $\mathcal{K}_C$. Stability under $\epsilon \partial_\epsilon$ follows from the stability of the aforementioned bounds in addition to the $\mathcal{L}(L^2)$-conormality of $B$. The $i = 2$ case is similar.
\end{proof}

We need a more general version where $A = A(\lambda), B = B(\lambda), C = C(\lambda)$ depend on $\lambda \in \C$. Suppose that $A(\lambda), C(\lambda)$ are holomorphic in $\Psi^{-\infty,\ast}_{\gas}$, in the sense of \cite{mazzeo-melrose_analytic-surgery_1995}. Suppose $B(\lambda)$ is holomorphic as a function valued in $\mathcal{A}([0,\epsilon_0)_\epsilon,\mathcal{L}(L^2(M, \Omega^{1/2})))$ for some $\epsilon_0 > 0$. If these assumptions hold, then Lemma \ref{lemma_semi-ideal} holds in a holomorphic sense; namely, the distributions $B_{E_i,\ell}(x_i \partial_{x_i}) \mathcal{K}_T $ are holomorphic, valued in $\mathcal{A}([0,\epsilon_0) \times M^2, \Omega^{1/2})$. The proof is almost identical. Equipped with this holomorphic semi-ideal-type property, we can move on to prove Theorem \ref{theorem_resolvent-holomorphic-away-from-spectrum}.

\begin{proof} (Theorem \ref{theorem_resolvent-holomorphic-away-from-spectrum})

Fix $K,L \in \R$ and let $N > 0$ be determined by Lemma \ref{lemma_semi-ideal}. Let $Q_0(\lambda), Q_0'(\lambda) \in \Psi^{-2,\ast}_{\gas}(M), {}^{(N)}Q(\lambda), {}^{(N)}Q'(\lambda), {}^{(N)}R(\lambda), {}^{(N)}R'(\lambda) \in \Psi^{-\infty,\ast}_{\gas}(M)$ be defined for $\lambda$ in a neighborhood of $\lambda_0$, as in Lemma \ref{lemma_laplacian-minus-lambda-parametrix-with-finite-order-remainder}.

By construction, ${}^{(N)}R(\lambda), {}^{(N)}R'(\lambda)$ vanish to leading order $N$ at all boundary hypersurfaces where $\{\epsilon = 0\}$. By choosing $N$ to be large enough, find $\epsilon_0 > 0$ such that, for $\epsilon \in (0,\epsilon_0)$ and $\lambda$ in the neighborhood of holomorphicity, $\|{}^{(N)}R_\epsilon(\lambda)\|_{x_\epsilon^{1-\alpha} x L^2(dV_{g_\epsilon})} < 1/2$. Hence $I - {}^{(N)}R_\epsilon(\lambda)$ is invertible on $x_\epsilon^{1-\alpha} x L^2(dV_g)$. Its inverse can be expressed as a Neumann series, so it is a holomorphic map to $\mathcal{A}([0,\epsilon_0)_\epsilon,\mathcal{L}(x_\epsilon^{1-\alpha} x L^2(dV_{g_\epsilon})))$.

Invoking the same argument for $I - {}^{(N)}R'_\epsilon(\lambda)$ and applying \eqref{eqn_laplacian-minus-lambda-parametrix-right} and \eqref{eqn_laplacian-minus-lambda-parametrix-left} shows that, for $\epsilon \in (0,\epsilon_0)$, $P_\epsilon(\lambda)$ is invertible on $x_\epsilon^{1-\alpha} x L^2(dV_{g_\epsilon})$ and the inverse is a holomorphic map to $\mathcal{A}([0,\epsilon_0)_\epsilon,\mathcal{L}(x_\epsilon^{1-\alpha} x L^2(dV_{g_\epsilon})))$.

In what follows, consider all operators to be restricted to $\epsilon \in [0,\epsilon_0)$. Left-multiply \eqref{eqn_laplacian-minus-lambda-parametrix-left} and right-multiply \eqref{eqn_laplacian-minus-lambda-parametrix-right} by $P(\lambda)^{-1}$ to obtain
\begin{align*}
    P(\lambda)^{-1} &= Q_0(\lambda) + {}^{(N)}Q(\lambda) + \left(P(\lambda)^{-1}\right)\left( {}^{(N)}R(\lambda) \right) \\
    P(\lambda)^{-1} &= Q_0'(\lambda) + {}^{(N)}Q'(\lambda) + {}^{(N)}R'(\lambda) P(\lambda)^{-1}
\end{align*}
Substituting the second identity into the first gives
\begin{align*}
    P(\lambda)^{-1} &= \tilde{Q}_0(\lambda) + {}^{(N)}\tilde{Q}(\lambda) + {}^{(N)}R'(\lambda) P(\lambda)^{-1} \left( {}^{(N)}R(\lambda) \right)
\end{align*}
for $\tilde{Q}_0(\lambda) \in {}^{\epsilon_0} \Psi^{-2,\ast}_{\gas}(M)$, ${}^{(N)}\tilde{Q}(\lambda) \in {}^{\epsilon_0} \Psi^{-\infty,\ast}_{\gas}(M)$, both depending holomorphically on $\lambda$ near $\lambda_0$.

We shall now show that $P(\lambda)^{-1}$ is holomorphic in ${}^{\epsilon_0} \Psi^{-2,\ast}_{\gas}(M)$, or, rather, that $P(\lambda)^{-1} - \tilde{Q}_0(\lambda)$ is holomorphic in ${}^{\epsilon_0} \Psi^{-\infty,\ast}_{\gas}(M)$. Apply Lemma \ref{lemma_semi-ideal} to the triple composition ${}^{(N)}R'(\lambda) P(\lambda)^{-1} \left( {}^{(N)}R(\lambda) \right)$. (As stated, Lemma \ref{lemma_semi-ideal} requires the middle operator to be bounded on $L^2(M, \Omega^{1/2})$, but this can be arranged by conjugating by weights.) By the remarks following Lemma \ref{lemma_semi-ideal}, ${}^{(N)}R'(\lambda) P(\lambda)^{-1} \left( {}^{(N)}R(\lambda) \right)$ satisfies the finite-order polyhomogeneous estimates \eqref{eqn_semi-ideal-phg-estimate} on $[0,1) \times M^2$ for $\ell < L$ and these estimates depend holomorphically on $\lambda$. Thus it also satisfies finite-order polyhomogeneous estimates at all boundary hypersurfaces of $\dbsp{M}{\gas}$. Of course, ${}^{(N)}\tilde{Q}(\lambda)$ is polyhomogeneous on $\dbsp{M}{\gas}$, so it follows that $P(\lambda)^{-1} - \tilde{Q}_0(\lambda)$ satisfies finite-order polyhomogeneous estimates on $\dbsp{M}{\gas}$. But, by choosing $K, L$ to be arbitrarily large, it follows that $P(\lambda)^{-1} - \tilde{Q}_0(\lambda)$ is truly (holomorphically) polyhomogeneous on $\dbsp{M}{\gas}$.
\end{proof}

\subsection{Near the spectrum of $\Delta_g$}

Now we describe the behavior of $P(\lambda)^{-1}$ for $\lambda$ near the spectrum of $\Delta_g$. Certainly, for each fixed $\epsilon$, $P_\epsilon(\lambda)^{-1}$ is meromorphic with a pole of order one at each of its eigenvalues. (Meromorphy in the operator norm topology follows from the analytic Fredholm theorem; meromorphy in a $\Psi$DO sense follows from the work of \cite{mazzeo-melrose_meromorphic-extension-resolvent_1987}.) However, as the eigenvalues vary with $\epsilon$, there is no hope that $P(\lambda)^{-1}$ is meromorphic in $\Psi^{-2,\ast}_{\gas}(M)$. Fortunately, the singularities of $P_\epsilon(\lambda)^{-1}$ are, at least locally, the zeroes of a holomorphic function depending polyhomogeneously in $\epsilon$. Canceling these singularities leaves a holomorphic family.

\begin{theorem}
\label{theorem_resolvent-meromorphic-near-spectrum}

    Suppose that $\lambda_0 \in \mathrm{spec} (\Delta_g)$. Then there exists $\epsilon_0 > 0$ and a holomorphic function $q(\lambda)$ defined in a neighborhood of $\lambda_0$, valued in $\mathcal{A}_{\mathrm{phg}}([0,\epsilon_0)_\epsilon)$, such that the zeroes of $q_\epsilon(\lambda)$ are eigenvalues of $\Delta_\epsilon$ and $q_\epsilon(\lambda) P(\lambda)^{-1}$ extends to  a holomorphic function of $\lambda$ on a neighborhood of $\lambda_0$, taking values in $^{\epsilon_0} \Psi^{-2,\ast}_{\gas}$.
\end{theorem}
\begin{proof}

    When $\lambda_0 \in \mathrm{spec} (\Delta_g)$, of course $N_{\bf{0}}(P(\lambda_0)) = x^{2-\alpha}(\Delta_g - \lambda_0)$ is not invertible and so the parametrix construction of 
    Lemma \ref{lemma_laplacian-minus-lambda-parametrix-with-finite-order-remainder} does not work for $P(\lambda_0)$. We thus modify it with a perturbation. Let $\Pi_{\lambda_0}$ be the orthogonal projection onto the $\lambda_0$-eigenspace of $\Delta_g$ in $x^{2-\alpha} L^2(dV_g)$. Let $\{ \phi^{(j)}_0 \}_{j=1}^d$ be an orthonormal basis for this eigenspace. By the standard regularity theory of elliptic 0-differential operators \cite{mazzeo_differential-edge-operators-I_1991}, each $\phi^{(j)}_0$ is polyhomogeneous on $M$, with index set generated by $(1,0), (\beta+1,0)$, so certainly $\Pi_{\lambda_0}$ has polyhomogeneous kernel on $\dbsp{M}{0'}$. Lift each $\phi^{(j)}_0$ to a polyhomogeneous function $\phi^{(j)} = \phi^{(j)}(\epsilon, x, y)$ on $\sgsp{M}{\gas}$ so that the functions $\phi^{(j)}_\epsilon(\cdot) = \phi^{(j)}(\epsilon, \cdot)$ are linearly independent, at least for $\epsilon < \epsilon_0$ for some $\epsilon_0 > 0$, and lie in $x_\epsilon^{1-\alpha} x L^2(dV_{g_\epsilon})$. (For example, $\phi^{(j)} = (x_\epsilon/x)^{1-\alpha+n\alpha/4} \phi^{(j)}_0$.) Then define $\Pi'_\epsilon : x_\epsilon^{1-\alpha} x L^2(dV_{g_\epsilon}) \to x_\epsilon^{1-\alpha} x L^2(dV_{g_\epsilon})$ by
    \begin{align*}
        \Pi'_\epsilon f = \sum_{j=1}^d \langle  \phi^{(j)}_\epsilon, f \rangle_{x_\epsilon^{1-\alpha} x L^2(dV_{g_\epsilon})}\phi^{(j)}_\epsilon
    \end{align*}
    For $\epsilon < \epsilon_0$, this is an operator of finite rank $d$. As the $\phi^{(j)}$ are polyhomogeneous, $\Pi'$ has polyhomogeneous kernel on $\dbsp{M}{\gas}$ and so $\Pi' \in \Psi^{-\infty,\ast}_{\gas}(M)$. By construction, $N_{\bf{0}}(\Pi') = \Pi_{\lambda_0}$.
    
    Then the normal operator $N_{\bf{0}}(P(\lambda) + x_\epsilon^{-\alpha} x^2  \Pi') = P_0(\lambda) + x^{2-\alpha} \Pi_{\lambda_0} = x^{2-\alpha}(\Delta_0 - \lambda + \Pi_{\lambda_0})$ is certainly invertible for $\lambda$ near $\lambda_0$, and the other normal operators of $P(\lambda) + x_\epsilon^{-\alpha} x^2  \Pi'$ agree with those of $P(\lambda)$ and hence $P(0)$. Then the parametrix construction of Lemma \ref{lemma_laplacian-minus-lambda-parametrix-with-finite-order-remainder} works equally well here and we obtain $Q_0(\lambda) \in {}^{\epsilon_0}\Psi^{-2,\ast}_{\gas}(M)$, ${}^{(N)}Q(\lambda) \in {}^{\epsilon_0}\Psi^{-\infty,\ast}_{\gas}(M)$, ${}^{(N)}R(\lambda) \in {}^{\epsilon_0}\Psi^{-\infty,\ast}_{\gas}(M)$ such that
    \begin{align*}
        \left( P(\lambda) + x_\epsilon^{-\alpha} x^2  \Pi' \right) \left( Q_0(\lambda) + {}^{(N)}Q(\lambda) \right) &= I - {}^{(N)}R(\lambda)
    \end{align*}
    where ${}^{(N)}R(\lambda)$ vanishes to order $N$ at all boundary hypersurfaces other than $\bf{2}$. We also obtain a left parametrix.

    Then, using the same arguments as in the proof of Theorem \ref{theorem_resolvent-holomorphic-away-from-spectrum}, we find that $P(\lambda) + x_\epsilon^{-\alpha} x^2  \Pi'$ is invertible (possibly requiring $\epsilon_0$ to be reduced) and its inverse is holomorphic near $\lambda_0$, valued in $\Psi^{-2,\ast}_{\gas}(M)$. The inverse is a parametrix for $P(\lambda)$:
    \begin{align}
        P(\lambda) \left(P(\lambda) + x_\epsilon^{-\alpha} x^2  \Pi'\right)^{-1} &= I - x_\epsilon^{-\alpha} x^2 \Pi' \left(P(\lambda) + x_\epsilon^{-\alpha} x^2 \Pi'\right)^{-1} =: I - x_\epsilon^{-\alpha} x^2 \Pi''(\lambda) \label{eqn_laplacian-minus-lambda-right-parametrix-with-finite-rank-error}
    \end{align}
    for some $\Pi''(\lambda)$, holomorphic in $^{\epsilon_0}\Psi^{-\infty,\ast}_{\gas}(M)$ for $\lambda$ near $\lambda_0$. For $\epsilon < \epsilon_0$, the operator $\Pi''_\epsilon(\lambda)$ has finite-dimensional range contained in $\Phi_\epsilon = \mathrm{span} \{ \phi^{(j)}_\epsilon \}_{j=1}^d$. Then $I - x_\epsilon^{-\alpha} x^2 \Pi''_\epsilon(\lambda)$ leaves the space $x_\epsilon^{-\alpha} x^2\Phi_\epsilon$ invariant, and it is invertible if and only if its restriction to this subspace is invertible. Let $\pi''_\epsilon(\lambda)$ be the matrix of $x_\epsilon^{-\alpha} x^2 \Pi''_\epsilon(\lambda)|_{x_\epsilon^{-\alpha} x^2 \Phi_\epsilon}$ with respect to the $x_\epsilon^{-\alpha} x^2 \phi^{(j)}_\epsilon$ and let $\pi''(\lambda)$ denote that same matrix viewed as a function of $\epsilon$. The entries of $I_d - \pi''(\lambda)$ are holomorphic in $\lambda$ and polyhomogeneous at $\epsilon = 0$ and hence so is its determinant, which we denote by $q_\epsilon(\lambda)$ or $q(\lambda)$, leaving the $\epsilon$-dependence implicit. For $\epsilon < \epsilon_0$ and $\lambda$ near $\lambda_0$, the function $q_\epsilon(\lambda)$ has zeroes precisely at the eigenvalues of $\Delta_\epsilon$, by \eqref{eqn_laplacian-minus-lambda-right-parametrix-with-finite-rank-error}.
    
    The matrix $q(\lambda) (I_d - \pi''(\lambda))^{-1}$, which is just the adjugate matrix of $I_d - \pi''(\lambda)$, extends continuously to the zeroes of $q(\lambda)$ and is holomorphic with polyhomogeneous entries. If $I_d - \pi''_\epsilon(\lambda)$ and hence $I - x_\epsilon^{-\alpha} x^2\Pi''_\epsilon(\lambda)$ is invertible, then, explicitly, $(I - x_\epsilon^{-\alpha} x^2 \Pi''_\epsilon(\lambda))^{-1} = I + (I-x_\epsilon^{-\alpha} x^2\Pi''_\epsilon(\lambda))|_{x_\epsilon^{-\alpha} x^2 \Phi_\epsilon}^{-1} (x_\epsilon^{-\alpha} x^2 \Pi''_\epsilon(\lambda))$. Then $q(\lambda) (I - x_\epsilon^{-\alpha} x^2\Pi''(\lambda))^{-1}$ extends continuously to the zeroes of $q(\lambda)$ and is holomorphic, valued in $\Psi^{0,\ast}_{\gas}(M)$. Right-multiplying \eqref{eqn_laplacian-minus-lambda-right-parametrix-with-finite-rank-error} by this operator, it follows that $q(\lambda) P(\lambda)^{-1}$ similarly extends to be holomorphic.
\end{proof}

\subsection{Results about the spectrum}

The characterization of $P(\lambda)^{-1}$ essentially reduces the study of the spectrum of $\Delta_\epsilon$ to a problem in classical perturbation theory. The following are straightforward consequences.

First, by integrating the resolvent around a small contour enclosing an eigenvalue of $\Delta_g$, we get the following:
\begin{corollary}
    \label{corollary_laplacian-spectral-projector-is-phg}
     
    Suppose that $\lambda_0 \in \mathrm{spec} (\Delta_g)$. Then there exists $\epsilon_0 > 0$ and $\delta > 0$ such that $\Pi$, the spectral projector of $\Delta_\epsilon$ onto $(\lambda_0 - \delta, \lambda_0 + \delta)$, lies in $^{\epsilon_0} \Psi^{-\infty,\ast}_{\gas}$.
\end{corollary}

Consequently, we can characterize simple eigenvalues and their corresponding eigenfunctions:
\begin{corollary}
    \label{corollary_laplacian-simple-eigenfunctions-eigenvalues-phg}

    Suppose $\lambda_0 \in \mathrm{spec} (\Delta_g)$ is a simple eigenvalue. Then there exists $\epsilon_0 > 0$ and $\delta > 0$ such that, for $\epsilon \in [0,\epsilon_0)$, $\mathrm{spec} (\Delta_\epsilon) \cap (\lambda_0 - \delta, \lambda_0 + \delta) = \{\lambda_\epsilon\}$. The eigenvalue $\lambda_\epsilon$ has multiplicity one and satisfies $\lim_{\epsilon \to 0} \lambda_\epsilon = \lambda_0$. In fact, $\lambda_\epsilon - \lambda_0$ is polyhomogeneous. The real part of its leading index is $\beta+1$.

    There also exists $\psi = \psi(\epsilon,x,y) \in \mathcal{A}_{\mathrm{phg}}({}^{\epsilon_0}\sgsp{M}{\gas})$ such that, for $\epsilon \in (0,\epsilon_0)$, $\psi(\epsilon,\cdot) \in x_\epsilon^{1-\alpha} x L^2(dV_{g_{\epsilon}})$ is a unit eigenfunction of $\Delta_\epsilon$. The real parts of the leading indices of $\psi$ at $\bf{1}, \tf[1], \bf{0}$ are $1,\beta+1,0$, respectively.
\end{corollary}
\begin{proof}

    If $\lambda_0$ is simple, then $\lambda_0$ is a zero of multiplicity one for the holomorphic function $q_0(\lambda)$ from the proof of Theorem \ref{theorem_resolvent-meromorphic-near-spectrum}. By Rouch\'{e}'s theorem, there exists a neighborhood of $\lambda_0$ in which $q_\epsilon(\lambda)$ has a unique zero, of multiplicity one, if $\epsilon$ is sufficiently small. Indeed, for any $\delta > 0$ sufficiently small, there exists $\epsilon_0$ such that $(\lambda_0 - \delta, \lambda_0 + \delta)$ contains a unique zero. This proves the first part, except for the polyhomogeneity of $\lambda_\epsilon - \lambda_0$.

    To construct $\psi$, let $\phi_0 \in x^{2-\alpha} L^2(dV_g)$ be a unit eigenfunction of $\Delta_0$ and extend it to a polyhomogeneous $\phi = \phi(\epsilon,x,y)$ on $\sgsp{M}{\gas}$ so that $\phi(0,\cdot) = \phi_0$ and $\phi_\epsilon(\cdot) := \phi(\epsilon,\cdot) \in x_\epsilon^{1-\alpha} x L^2(dV_{g_\epsilon})$, much as was done in the proof of Theorem \ref{theorem_resolvent-holomorphic-away-from-spectrum}. As before, let $\Pi$ be the spectral projector for $\Delta_\epsilon$ onto $(\lambda_0 - \delta, \lambda_0 + \delta)$. The projection $\Pi_\epsilon \phi_\epsilon$ is polyhomogeneous on $\sgsp{M}{\gas}$ by Corollary \ref{corollary_laplacian-spectral-projector-is-phg} and the mapping properties of \ztrtext{} operators (Proposition \ref{proposition_gas-giant-phg-mapping-properties}). As it limits to $\phi_0$ as $\epsilon \to 0$, it is certainly nonzero for $\epsilon$ close to $0$. Then we can define
    \begin{align*}
        \psi_\epsilon &= \frac{\Pi_\epsilon \phi_\epsilon}{\left\| \Pi_\epsilon \phi_\epsilon \right\|_{x_\epsilon^{1-\alpha} x L^2(dV_{g_\epsilon})}}
    \end{align*}
    This is a unit $\Delta_\epsilon$-eigenfunction for each $\epsilon$. The denominator is polyhomogeneous in $\epsilon$ by the pushforward theorem, so indeed $\psi_\epsilon$ is polyhomogeneous on ${}^{\epsilon_0}\sgsp{M}{\gas}$ if $\epsilon_0 > 0$ is chosen sufficiently small. Let $\psi(\epsilon,x,y) = \psi_\epsilon(x,y)$.

    Writing $\lambda_\epsilon = \langle \Delta_\epsilon \psi_\epsilon, \psi_\epsilon \rangle_{x_\epsilon^{1-\alpha} x L^2(dV_{g_\epsilon})}$, it is then immediate that $\lambda_\epsilon$ and hence $\lambda_0 - \lambda_\epsilon$ are polyhomogeneous.

    To read off the leading asymptotics of $\psi$, one could explicitly compute the index family of the spectral projector, but there is an easier approach via the parametrix construction. Let $\check{P}_\epsilon = P_\epsilon(\lambda_\epsilon)$. This operator is in the \ztrtext{} calculus because $\lambda_\epsilon$ is polyhomogeneous. The parametrix construction from Theorem \ref{theorem_laplacian-parametrix-with-compact-remainder-phg-index-family} works equally well for $\check{P}$ because it has the same normal operators at $\tff, \btff, \zff, \bff$, and $\bf{0}$ as $P(\lambda_0)$. We are unable to replicate the parametrix construction of Theorem \ref{theorem_laplacian-parametrix-with-infinite-order-remainder} exactly, because the normal operator at $\bf{0}$ is not invertible. But, in fact, by following the proof of that theorem (and taking adjoints, because that proof was written for a right parametrix), we can still obtain a left parametrix $\check{Q}$ so that $\check{Q} \check{P} = I - \check{R}$, where the kernel of $\check{R}$ has index set $\emptyset$ at all boundary hypersurfaces except for $\bf{1}, \tf[1]$, and $\bf{0}$, where the real parts of the leading indices are $1, \beta+1, 0$. By construction $\check{P} \psi = 0$, so $\psi = \check{R} \psi$. By mapping properties of the \ztrtext{} calculus (Proposition \ref{proposition_gas-giant-phg-mapping-properties}), $\psi$ has the claimed leading asymptotics.
\end{proof}

This result extends the characterization of eigenvalues and eigenfunctions of $\Delta_\epsilon$ from \cite{dHIKM_gas-giants_2024}. With more careful index bookkeeping, sharper estimates of the index families of the resolvent, eigenvalues, and eigenfunctions of $\Delta_\epsilon$ could be obtained.

The assumption of simplicity is not a fundamental limitation. The point of Theorem \ref{theorem_resolvent-meromorphic-near-spectrum} is to reduce the perturbation theory of the gas giant spectrum to classical perturbation theory. We expect that the non-simple case could be handled largely with tools from classical perturbation theory, without much further input from the \ztrtext{} calculus.

  \printbibliography

\end{document}